\documentclass[preprint,11pt]{elsarticle}
\usepackage[bookmarks=true,colorlinks=true,linkcolor=blue]{hyperref}

\usepackage{amsmath,amssymb,amsthm,mathabx}
\usepackage{graphicx,esint,ulem}

\biboptions{sort&compress}
\usepackage[usenames,dvipsnames,svgnames,table]{xcolor}

\usepackage{color}
\definecolor{jwbGreen}{rgb}{0, .6, 0}
\definecolor{jbaPurple}{HTML}{6600FF}
\definecolor{purple}{rgb}{.7, 0., .8}
\newcommand{\red}{\color{red}}

\newcommand{\blue}{\color{blue}}

\usepackage{calc}
\usepackage[margin=1.in]{geometry}

\usepackage{empheq}
\usepackage[most]{tcolorbox}

\newtcbox{\mymath}[1][]{%
    nobeforeafter, math upper, tcbox raise base,
    enhanced, colframe=blue!60!black,
    colback=blue!20, boxrule=1pt,
    #1}

\usepackage{algorithm,algpseudocode}
\usepackage{algorithmicx}
\algrenewcommand\alglinenumber[1]{\footnotesize #1:} 

\usepackage{tcolorbox}
\usepackage[english]{babel}
\usepackage{amsmath}
\usepackage{amsfonts}
\usepackage{fullpage}
\usepackage{graphicx}
\usepackage{setspace}
\usepackage{float}

\usepackage[intoc]{nomencl}
\makenomenclature

\usepackage{tabu}

\usepackage{multirow}

\usepackage{fancyvrb}

\usepackage{mdframed}

\allowdisplaybreaks

\usepackage{tikz}
\usetikzlibrary{plotmarks}

\newtheorem{theorem}{Theorem}
\newtheorem*{theoremNoNumber}{Theorem}

\newtheorem{definition}{Definition}

\newtheorem{lemma}{Lemma}[section]

\newcommand{\citeCount}[1]{}

\newcommand{\bogus}[1]{{}}
\newcommand{\ssf}{\scriptscriptstyle}

\newcommand{\mni}{\medskip\noindent}

\newcommand{\labelFont}{\footnotesize}

\newlength{\ycbTop}
\newlength{\ycbMid}%

\newcommand{\p}{\partial}

\newcommand{\f}[2]{\frac{#1}{#2}}

\def\ba#1\ea{\begin{align}#1\end{align}}

\def\bas#1\eas{\begin{align*}#1\end{align*}}

\def\bat#1\eat{\begin{alignat}{3}#1\end{alignat}}

\def\bats#1\eats{\begin{alignat*}{3}#1\end{alignat*}}

\newcommand{\bse}{\begin{subequations}}
\newcommand{\ese}{\end{subequations}}

\newcommand{\dpmx}{D_{+x}D_{-x}}

\newcommand{\dpmt}{D_{+t}D_{-t}}

\newcommand{\Dzt}{D_{0t}}
\newcommand{\Dpt}{D_{+t}}
\newcommand{\Dmt}{D_{-t}}
\newcommand{\dt}{\Delta t}

\newcommand{\dr}{{\Delta r}}

\newcommand{\defeq}{\overset{{\rm def}}{=}}
\newcommand{\eqdef}{\overset{{\rm def}}{=}}

\newcommand{\ev}{\mathbf{ e}}

\newcommand{\iv}{\mathbf{ i}}
\newcommand{\jv}{\mathbf{ j}}
\newcommand{\kv}{\mathbf{ k}}

\newcommand{\rv}{\mathbf{ r}}

\newcommand{\xv}{\mathbf{ x}}

\newcommand{\Gv}{\mathbf{ G}}

\newcommand{\Lv}{\mathbf{ L}}

\newcommand{\Qv}{\mathbf{ Q}}
\newcommand{\Rv}{\mathbf{ R}}

\newcommand{\Uv}{\mathbf{ U}}
\newcommand{\Vv}{\mathbf{ V}}

\newcommand{\half}{{1\over2}}

\newcommand{\Real}{{\mathbb R}}

\newcommand{\Complex}{{\mathbb C}}
\newcommand{\Integers}{{\mathbb Z}}

\newcommand{\zerov}{\mathbf{0}}

\newcommand{\Bc}{{\mathcal B}}

\newcommand{\Ec}{{\mathcal E}}

\newcommand{\Gc}{{\mathcal G}}

\newcommand{\Lc}{{\mathcal L}}

\newcommand{\Oc}{{\mathcal O}}

\newcommand{\grad}{\nabla}

\newcommand{\nd}{n_d}

\newcommand{\eps}{\epsilon}

\newcommand{\ds}{{\Delta s}}

\newcommand{\Dmx}{D_{-x}}

\newcommand{\lap}{\Delta}  

\newcommand{\lamHat}{\hat{\lambda}}
\newcommand{\Uhat}{\hat{U}}

\newcommand{\Lss}{{\ssf L}}
\newcommand{\Rss}{{\ssf R}}

\newcommand{\Nuc}{{n_{\rm u}}}
\newcommand{\CFL}{{\rm CFL}}

\newcommand{\safetyFactor}{s_f}
\newcommand{\Rhat}{\hat{R}}

\newlength{\tfwidth}
\newlength{\tfheight}
\newlength{\tfxa}
\newlength{\tfxb}
\newlength{\tfya}
\newlength{\tfyb}
\newcommand{\trimFigWithBox}[6]{%
\setlength\fboxsep{0pt}%
\setlength\fboxrule{1.0pt}
\fbox{\includegraphics[width=#2, clip, trim=#3 #4 #5 #6]{#1}}%
}
\newcommand{\trimFigNoBox}[6]{%
\setlength\fboxsep{1pt}
\setlength\fboxrule{0.0pt}
\fbox{\includegraphics[width=#2, clip, trim=#3 #4 #5 #6]{#1}}%
}
\newcommand{\trimFigHeightWithBox}[6]{%
\setlength\fboxsep{0pt}%
\setlength\fboxrule{1.0pt}
\fbox{\includegraphics[height=#2, clip, trim=#3 #4 #5 #6]{#1}}%
}
\newcommand{\trimFigHeightNoBox}[6]{%
\setlength\fboxsep{1pt}
\setlength\fboxrule{0.0pt}
\fbox{\includegraphics[height=#2, clip, trim=#3 #4 #5 #6]{#1}}%
}

\newsavebox\figBox

\newcommand{\trimw}[6]{%
\sbox\figBox{\includegraphics{#1}}
\setlength{\tfwidth}{\the\wd\figBox}
\setlength{\tfheight}{\the\ht\figBox}
\setlength{\tfxa}{\tfwidth*\real{#3}}%
\setlength{\tfxb}{\tfwidth*\real{#4}}%
\setlength{\tfya}{\tfheight*\real{#5}}%
\setlength{\tfyb}{\tfheight*\real{#6}}%
\trimFigNoBox{#1}{#2}{\tfxa}{\tfya}{\tfxb}{\tfyb}%
}

\newcommand{\trimwb}[6]{%

\sbox\figBox{\includegraphics{#1}}
\setlength{\tfwidth}{\the\wd\figBox}
\setlength{\tfheight}{\the\ht\figBox}
\setlength{\tfxa}{\tfwidth*\real{#3}}%
\setlength{\tfxb}{\tfwidth*\real{#4}}%
\setlength{\tfya}{\tfheight*\real{#5}}%
\setlength{\tfyb}{\tfheight*\real{#6}}%
\trimFigWithBox{#1}{#2}{\tfxa}{\tfya}{\tfxb}{\tfyb}%
}

\newcommand{\trimh}[6]{%
\sbox\figBox{\includegraphics{#1}}
\setlength{\tfwidth}{\the\wd\figBox}
\setlength{\tfheight}{\the\ht\figBox}
\setlength{\tfxa}{\tfwidth*\real{#3}}%
\setlength{\tfxb}{\tfwidth*\real{#4}}%
\setlength{\tfya}{\tfheight*\real{#5}}%
\setlength{\tfyb}{\tfheight*\real{#6}}%
\trimFigHeightNoBox{#1}{#2}{\tfxa}{\tfya}{\tfxb}{\tfyb}%
}

\newcommand{\trimhb}[6]{%

\sbox\figBox{\includegraphics{#1}}
\setlength{\tfwidth}{\the\wd\figBox}
\setlength{\tfheight}{\the\ht\figBox}
\setlength{\tfxa}{\tfwidth*\real{#3}}%
\setlength{\tfxb}{\tfwidth*\real{#4}}%
\setlength{\tfya}{\tfheight*\real{#5}}%
\setlength{\tfyb}{\tfheight*\real{#6}}%
\trimFigHeightWithBox{#1}{#2}{\tfxa}{\tfya}{\tfxb}{\tfyb}%
}
\usepackage{xargs}

\newcommandx{\figByHeight}[9][5=0, 6=0, 7=0, 8=0,9=]{
\draw (#1,#2) node[anchor=south west,xshift=-16pt,yshift=-4pt] {\trimh{#3}{#4}{#5}{#6}{#7}{#8}};}
\newcommandx{\figByHeightb}[9][5=0, 6=0, 7=0, 8=0,9=]{
\draw (#1,#2) node[anchor=south west,xshift=-16pt,yshift=-4pt] {\trimhb{#3}{#4}{#5}{#6}{#7}{#8}};}

\newcommandx{\figByHeightWithLabel}[9][5=0, 6=0, 7=0, 8=0,9=]{
\draw (#1,#2) node[anchor=south west,xshift=-16pt,yshift=-4pt] {\trimh{#3}{#4}{#5}{#6}{#7}{#8}} node[draw=white,fill=white,inner sep=1pt,anchor=south west] {#9};}
\newcommandx{\figByHeightWithLabelb}[9][5=0, 6=0, 7=0, 8=0,9=]{
\draw (#1,#2) node[anchor=south west,xshift=-16pt,yshift=-4pt] {\trimhb{#3}{#4}{#5}{#6}{#7}{#8}} node[draw=white,fill=white,inner sep=1pt,anchor=south west] {#9};}

\newcommandx{\figByWidth}[9][5=0, 6=0, 7=0, 8=0,9=]{
\draw (#1,#2) node[anchor=south west,xshift=-16pt,yshift=-4pt] {\trimw{#3}{#4}{#5}{#6}{#7}{#8}};}
\newcommandx{\figByWidthb}[9][5=0, 6=0, 7=0, 8=0,9=]{
\draw (#1,#2) node[anchor=south west,xshift=-16pt,yshift=-4pt] {\trimwb{#3}{#4}{#5}{#6}{#7}{#8}};}

\newcommandx{\figByWidthWithLabel}[9][5=0, 6=0, 7=0, 8=0,9=]{
\draw (#1,#2) node[anchor=south west,xshift=-16pt,yshift=-4pt] {\trimw{#3}{#4}{#5}{#6}{#7}{#8}} node[draw=white,fill=white,inner sep=1pt,anchor=south west] {#9};}
\newcommandx{\figByWidthWithLabelb}[9][5=0, 6=0, 7=0, 8=0,9=]{
\draw (#1,#2) node[anchor=south west,xshift=-16pt,yshift=-4pt] {\trimwb{#3}{#4}{#5}{#6}{#7}{#8}} node[draw=white,fill=white,inner sep=1pt,anchor=south west] {#9};}

\usepackage{todonotes}

\definecolor{jwbGreen}{rgb}{0, .6, 0}
\definecolor{jbaPurple}{HTML}{6600FF}
\definecolor{purple}{rgb}{.7, 0., .8}
\definecolor{pinegreen}{rgb}{0.0, 0.47, 0.44}

\begin{document}

\begin{frontmatter}
 \title{High-order Accurate Implicit-Explicit Time-Stepping Schemes for Wave Equations on Overset Grids}

\author[rpi]{Allison~M. Carson\fnref{NSFgrants}}
\ead{carsoa@rpi.edu}

\author[rpi]{Jeffrey W.~Banks\fnref{DOE}}
\ead{banksj3@rpi.edu}

\author[rpi]{William D.~Henshaw\corref{cor}\fnref{NSFgrants}}
\ead{henshw@rpi.edu}

\author[rpi]{Donald~W.~Schwendeman\fnref{NSFgrants}}
\ead{schwed@rpi.edu}

\address[rpi]{Department of Mathematical Sciences, Rensselaer Polytechnic Institute, Troy, NY 12180, USA}

\cortext[cor]{Corresponding author}

\fntext[NSFgrants]{Research supported by the National Science Foundation under grants DMS-1519934 and DMS-1818926.}

\fntext[DOE]{This work was partially performed under DOE contracts from the ASCR Applied Math Program.}

\begin{abstract}

  New implicit and implicit-explicit 
time-stepping methods for the wave equation in second-order form are described with application to two and three-dimensional problems discretized on overset grids.
The implicit schemes are single step, 
three levels in time, and based on the modified equation approach.
Second and fourth-order accurate schemes are developed and they incorporate upwind dissipation for stability on overset grids.
The fully implicit schemes are useful for certain applications such as the WaveHoltz algorithm for solving Helmholtz problems where very large time-steps are desired.
Some wave propagation problems are geometrically stiff due to localized regions
of small grid cells, such as grids needed to resolve fine geometric features, and for these situations the
implicit time-stepping scheme is combined with an explicit scheme:
the implicit scheme is used for component grids containing small cells while the
explicit scheme is used on the other grids such as background Cartesian grids.
The resulting partitioned implicit-explicit scheme can be many times faster than
using an explicit scheme everywhere. The accuracy and stability of the schemes are studied
through analysis and numerical computations.

\end{abstract}


\end{frontmatter}

\clearpage
\tableofcontents

\clearpage

\nomenclature[01]{ME:}{Modified Equation}
\nomenclature[02]{EMEp:}{Explicit Modified Equation scheme, $p^{\rm th}$ order accurate}
\nomenclature[03]{IMEp:}{Implicit Modified Equation scheme, $p^{\rm th}$ order accurate}
\nomenclature[04]{SPIEp:}{Spatially Partitioned Implicit Explicit scheme, $p^{\rm th}$ order accurate}
\nomenclature[05]{SCHEMEp-UW-PC:}{SCHEMEp + upwind disssipation + predictor corrector}
\nomenclature[06]{LTS:}{Local Time Stepping}
\nomenclature[07]{LIM:}{Locally Implicit Method}

\printnomenclature

\section{Introduction} \label{sec:introduction}

The wave equation in second-order form is an important model for many applications in science and engineering involving wave propagation.
Example applications include acoustics, electromagnetics and elasticity;
such problems are often posed mathematically as partial differential equations with appropriate initial and boundary conditions.
Wave propagation problems are often solved most efficiently 
using high-order accurate explicit time-stepping schemes. Explicit schemes can be fast and memory efficient. The time-step
for such schemes is limited by the usual CFL stability condition involving the size of grid cells and the wave speed.
Thus, there are some situations, such as with a locally fine mesh or a locally large wave speed, when an explicit scheme
with a global time-step is
inefficient since a small time-step would be required everywhere.
For such situations, we say the problem is \textsl{geometrically stiff} or \textit{materially stiff}.
An example of a geometrically stiff problem is the diffraction of an incident wave from a knife-edge as shown in Figure~\ref{fig:knifeEdge}.
The solution of this problem is computed using an overset grid for which there are small grid cells near the tip
of the knife-edge.  These small cells force the time-step of an explicit method to be reduced by a factor of $20$ from that required by the Cartesian background grid.  (More information concerning overset grids and our numerical schemes for such grids is given later.)
There are two common approaches to overcome this stiffness, local time-stepping (LTS) and locally implicit methods (LIM).
LTS methods use a local time-step dictated by the local time-step restriction. LIM's use an implicit method on only part of the domain, usually where the grid cells are smallest.

{
\newcommand{\drawContour}[7]{%
\begin{scope}[#1]
\draw(0.0,0) node[anchor=south west,xshift=-4pt,yshift=+0pt] {\trimfiga{#2}{\figWidtha}};
\begin{scope}[xshift=0cm,yshift=-2pt]
  \draw (\xcb,\ycb) node[anchor=south west,xshift=0.25cm,yshift=.5cm,rotate=-90] {\trimfigcb{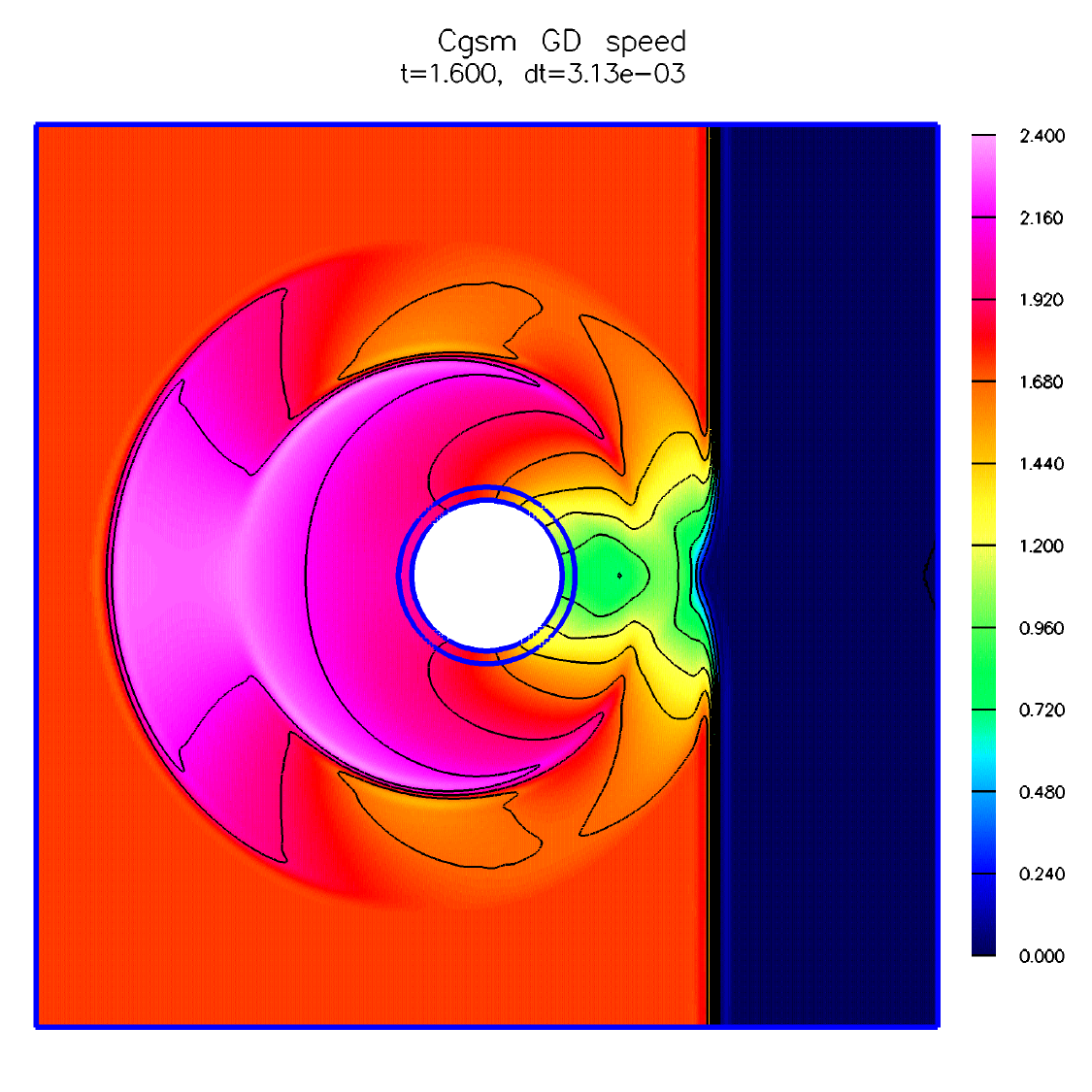}{\cbWidth}{\cbHeight}};
  \draw (.8,0) node[anchor=north,xshift=+3pt,yshift=+2pt] {\scriptsize $#6$};
  \draw (4.8,0) node[anchor=north,xshift=+0pt,yshift=+2pt] {\scriptsize $#7$};
\end{scope}
\end{scope}
}
\newcommand{\cbWidth}{.2cm}
\newcommand{\cbHeight}{4cm}
\newcommand{\xcb}{.5cm}
\newcommand{\ycb}{-.2cm}
\setlength{\ycbTop}{\ycb+\cbHeight}
\setlength{\ycbMid}{\ycb+\cbHeight*\real{.5}}
\newcommand{\trimfigcb}[3]{\includegraphics[width=#2, height=#3, clip, trim=17cm 2.35cm 1.65cm 2.35cm]{#1}}
\newcommand{\figWidtha}{5.15cm}
\newcommand{\trimfiga}[2]{\trimw{#1}{#2}{.05}{.12}{.07}{.1}}
\begin{figure}[htb]
\begin{center}
\begin{tikzpicture}
   \useasboundingbox (0,.3) rectangle (14.5,5);  

   \figByHeight{0.0}{0}{tipGridThine64Zoom1}{5.cm}[0.1][0.2][0.0][0.0]
   \figByHeight{4.5}{0}{tipGridThine64Zoom2}{5.cm}[0.1][0.1][0.0][0.0]
   \figByHeightb{2.5}{1}{tipGridThine64Zoom3}{3.cm}[0.1][0.1][0.1][0.1]

   \begin{scope}[xshift=8.5cm,yshift=-5pt]
     \drawContour{xshift=0.0cm,yshift=0.cm}{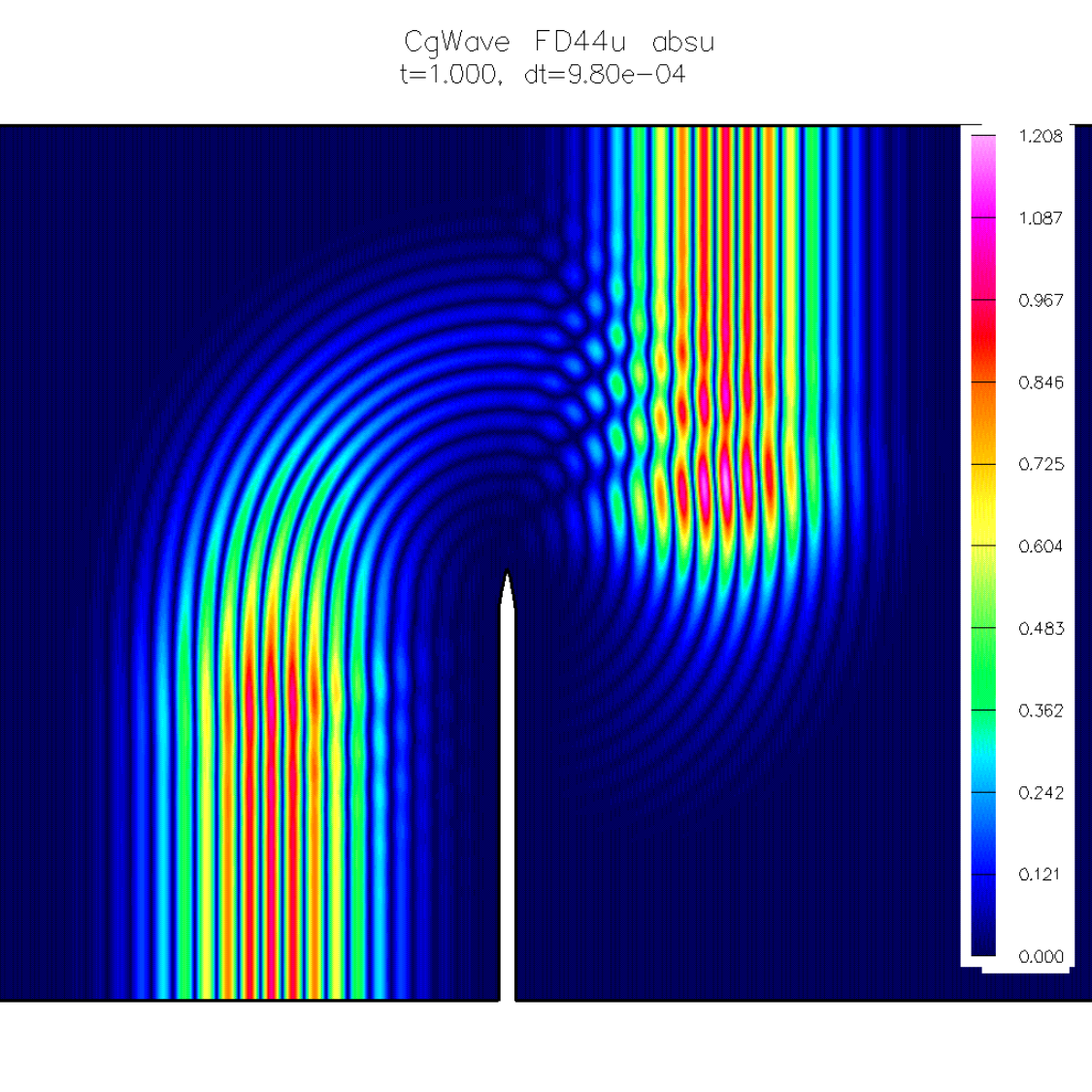}{$|u|$, $t=1.0$}{$v$}{$t=1.0$}{$0.$}{$1.2$}    
   \end{scope}   

\end{tikzpicture}
\end{center}
\caption{Geometrically stiff problem: scattering of a modulated Gaussian plane wave from a knife edge.
    Left: overset grid for the geometry showing magnified views of the tip grid which has very small grid cells.
    Right: contours of $|u|$ computed with the new SPIE scheme; the tip grid was advanced implicitly while other
    grids were advanced explicitly resulting in a time-step that was about $20$ times larger than using an explicit scheme on all grids.
    }
\label{fig:knifeEdge}
\end{figure}
}

In this article we develop new locally implicit time-stepping schemes for the wave equation in second-order form on overset grids based
on the modified equation (ME) approach. These schemes are high-order accurate single-step schemes that use three time-levels and a compact spatial stencil. The schemes depend on one or more parameters that determine the degree of implicitness; the second-order accurate scheme depends on one parameter while the fourth-order accurate scheme depends on two parameters. For certain ranges of these parameters the schemes are unconditionally stable in time.
A small amount of upwind dissipation is added to the schemes for stability on overset grids. The upwind dissipation
can be added in several ways, for example, in a fully implicit manner or in a predictor-corrector fashion where the upwinding is added in a separate explicit step.

Our implicit time-stepping ME scheme, denoted by IME, is combined with a ME-based explicit time-stepping scheme, denoted by EME, in a spatially partitioned manner. 
The EME schemes we use have a compact stencil and have a time-step restriction that is independent 
of the order of accuracy\footnote{Many EME schemes take powers of an matrix operator (leading to a wider stencil) and the time-step restriction depends on the order of accuracy.}.
We say that these compact EME schemes are able to take a \textit{CFL-one} time-step.
This is in contrast to typical linear multi-step methods where explicit higher-order schemes tend to have smaller time-step restrictions,
or to popular Discontinuous Galerkin (DG) methods where the time-step restriction
 typically scales as $1/( 2 P + 1)$ with $P$ being the degree
of the polynomial basis~\cite{TaubeDumbserMunzSchneider2009}.
For an overset grid, a locally implicit scheme can be used, for example,
on boundary-fitted component grids that resolve small geometrical features. The EME scheme can then be used on Cartesian background grids or on curvilinear grids that have similar grid spacings to the background grids. In this way the time-step for the EME scheme is not restricted by the small grid cells used to resolve small geometrical features.
In typical applications, the majority of the grid points belong to background Cartesian grids, and
the solution on these grid points can be advanced very efficiently with a CFL-one time-step. This can make the hybrid IME-EME scheme
much more efficient than using the EME scheme everywhere with a small (global) CFL time-step.  We refer to this hybrid scheme as a Spatially Partitioned Implicit-Explicit (SPIE) scheme.
Note that the EME scheme is more accurate since ME schemes are most accurate for the CFL number close to one: unlike method-of-lines schemes, the accuracy of ME schemes is degraded for small CFL numbers. The implicit matrix formed by the IME schemes is definite, and it is well suited to a solution by modern Krylov-based methods or multigrid.

Normally there is no benefit in using implicit time-stepping and taking a large (greater than one) global CFL time-step for wave propagation problems as the accuracy of the solution is usually degraded.
However, there are applications where implicit time-stepping methods for the wave equation using a large CFL number can be useful.  For example, implicit methods for the wave equation are an attractive option for each iteration step of the
WaveHoltz algorithm~\cite{appelo2020waveholtz,EmWaveHoltzPengAppelo2022,appeloElWaveHoltz2022} which solves for time periodic (Helmholtz) 
solutions\footnote{Note that the dispersion errors due to the large CFL time-step can be eliminated by an adjustment to the forcing frequency.}.
The WaveHoltz algorithm can solve Helmholtz problems for frequencies anywhere in the spectrum without the need to invert an indefinite matrix as is common with many approaches. Each iteration of the WaveHoltz algorithm requires a solution of a wave equation over a given period, and just a few implicit time-steps per period (e.g.~5--10) are needed which leads to very large CFL numbers on fine grids.  This is one of our motivations for developing stable IME schemes for overset grids.

There is a large literature on ME, LTS, and LIM schemes for the wave equation. 
Here we provide a brief synopsis, for further references, see~\cite{GroteMehlinMitkova2015,AlmquistMehlin2017}, for example.
Explicit ME schemes for the wave equation go back, at least, to the work of Dablain~\cite{Dablain1986} and of Shubin and Bell~\cite{ShubinBell1987}. 
Local time-stepping has most often been used for PDEs that are written as first-order systems in time.
LTS has been used for decades with adaptive mesh refinement (AMR) 
since the pioneering work of Berger and Oliger~\cite{Berger3}.
Local time-stepping has also been developed, for example, for Runge-Kutta time-stepping~\cite{GroteMehlinMitkova2015,AlmquistMehlin2017,LiuLiHu2014} and arbitrary high-order ADER schemes~\cite{TaubeDumbserMunzSchneider2009}.
Of note is the ME-based LTS method for the wave equation of Diaz and Grote~\cite{DiazGrote2009}, where it was found
necessary to have a small overlap of one or two cells between the coarse and fine cells in order 
to retain the time-step dictated by the coarse mesh.
Also of note is the locally time-stepped Runge-Kutta finite difference scheme of Liu, Li, and Hu~\cite{LiuLiHu2014}
for the wave equation on block-structured grids.
Another reason for using local time-stepping is to couple two difference schemes together.
For example, Beznsov and Appel\"o~\cite{BeznosovAppelo2021} use local time-stepping for the wave equation to couple a DG scheme (which has a small time-step) with an Hermite scheme (which has a CFL-one time-step). The DG scheme is used on boundary-fitted grids of an overset grid for accurate treatment of the boundary conditions.

Compact implicit difference approximations lead to globally implicit systems (although sometimes with a time-step restriction) and these have been used for 
wave equations by a number of authors~\cite{BrittTurkelTsynkov2018,LiLiaoLin2019,KahanaSmithTurkelTsynkov2022}. 
To overcome the cost of the global
implicit solve, it is common to use locally one-dimensional approximate factorizations such as the alternating-direction-implicit (ADI) scheme~\cite{LimKimDouglas2007,KimLim2007}. 
A variety of locally implicit methods for wave equations have been developed, for example~\cite{Piperno2006,Verwer2011,ChabassierImperiale2016}.
For some formulations care is required in coupling the implicit and explicit schemes 
to avoid an order of accuracy reduction in time.
Of particular note is the 
fourth-order accurate locally implicit method for the wave equation of Chabassier and Imperiale~\cite{ChabassierImperiale2016}.
They use a Finite Element Method (FEM) discretization (with mass lumping to form a diagonal mass matrix) and 
a mortar element method with Lagrange multipliers to couple the implicit and explicit methods.
The implicit ME scheme in~\cite{ChabassierImperiale2016} is similar to our implicit scheme except that we use finite difference approximations and
a more compact approximation (which leads to different stability restrictions). 
Our approach uses a simple coupling between implicit and explicit regions based on overset grid interpolation. The price for this simpler
coupling is that upwind dissipation is needed to ensure stability.

We have been developing high-order accurate algorithms for a variety of wave propagation problems on overset grids. These include the solution to Maxwell's equations of electromagnetics for linear and nonlinear dispersive materials ~\cite{mxsosup2018,adegdm2019,adegdmi2020,maxwellMLA2022,ssmx2023}, the solution of linear and non-linear compressible elasticity~\cite{smog2012,flunsi2016} and the solution of incompressible elasticity~\cite{ism2023}.
A fourth-order accurate ME scheme for Maxwell's equations in second-order form on overset grids was developed in~\cite{max2006b}.
Extensions of the implicit and implicit-explicit time-stepping methods developed in this article will be very useful for these other applications, both to treat geometric stiffness and for solving Helmholtz problems using the WaveHoltz algorithm.

The work presented in the remaining sections of this article are organized as follows.  Explicit and implicit ME schemes for the wave equation are introduced in Section~\ref{sec:implicitExplcitMESchemes}, which also serves to establish some notation.  Details of the second and fourth-order accurate IME schemes for Cartesian grids are given in Section~\ref{sec:IMEonCartesian} where a von Neumann stability analysis is also performed.  Methods of upwind dissipation for IME schemes are described in Section~\ref{sec:upwindDissipation}, and this is followed in Section~\ref{sec:spieSchemes} by a formulation and a GKS stability analysis of our new SPIE schemes.  Section~\ref{sec:oversetGrids} discusses the implementation of the new ME schemes for overset grids, and Section~\ref{sec:matrixStability} provides results of a matrix stability analysis the ME schemes for one-dimensional overset grids.  Numerical results are discussed in Section~\ref{sec:numericalResults} and concluding remarks are offered in Section~\ref{sec:conclusions}.

\section{Three-level explicit and implicit ME schemes for the wave equation} \label{sec:implicitExplcitMESchemes}

We are interested in solving an initial-boundary-value (IBVP) problem for the wave equation in second-order form for a function $u(\xv,t)$ on a domain $\Omega$, with boundary $\Gamma$, in $\nd$ 
space dimensions,
\bse
\label{eq:ibvp}
\begin{alignat}{3}
    \partial_t^2 u      & =  \Lc u, \ \ \ \ &&\xv\in\Omega, \ \ t> 0, \label{eq:secondOrderWaveEquation}\\
    u(\xv,0)            & = u_0(\xv), \ \  \ \ &&\xv\in\Omega,  \label{eq:IC1} \\
    \partial_t u(\xv,0) & = u_1(\xv), \ \ \ \ &&\xv\in\Omega,   \label{eq:IC2} \\
    \Bc u(\xv,t)        & =g(\xv,t),  \ \ \ \ &&\xv\in\Gamma, \ \ t> 0.
\end{alignat}
\ese
Here, $\xv=[x_1,...,x_{\nd}]^T\in\Real^{\nd}$ is the vector of spatial coordinates, $t$ is time, and $\Lc$ is the spatial part of the wave operator, 
\ba  
 \Lc \eqdef c^2 \lap u, \quad  \lap \eqdef \sum_{d=1}^{\nd}\partial_{x_d}^2,
\ea
with wave speed $c>0$.
The initial conditions on $u$ and $\partial_t u$ are specified by the given functions $u_0(\xv)$ and $u_1(\xv)$, respectively, and the boundary conditions, denoted by the boundary condition operator $\Bc$, may be of Dirichlet, Neumann, or Robin type with given boundary data $g(\xv,t)$. 

\medskip
We begin with a description of the three-level ME schemes that ignores
the specifics of the spatial discretizations.
 Details of the grids and
spatial discretization are left to later sections.
The explicit and implicit ME schemes are both based on the standard second-order accurate 
central difference approximation to the second time-derivative of $u$,
\ba
  \Dpt\Dmt\, u & = \f{u(\xv,t+\dt)-2 u(\xv,t)+u(\xv,t-\dt)}{\dt^2} ,  \label{eq:DptDmtu}
\ea
where $\dt$ is the time-step and  
$\Dpt$ and $\Dmt$ are forward and backward divided difference operators
in time given by
\ba
   \Dpt\, u(\xv,t) \eqdef \f{u(\xv,t+\dt) - u(\xv,t)}{\dt} , \quad
   \Dmt\, u(\xv,t) \eqdef \f{u(\xv,t) - u(\xv,t-\dt)}{\dt} ,
\ea
respectively.  Expanding the terms in~\eqref{eq:DptDmtu} using Taylor series gives the following expansion, 
\ba
  \Dpt\Dmt\, u &  = \partial_t^2 u + \frac{\dt^2}{12}\partial_t^4 u + \frac{\dt^4}{360}\partial_t^6 u + \cdots.
  \label{eq:taylorExpansion} 
\ea 
Even time derivatives on the right-hand side of~\eqref{eq:taylorExpansion} are replaced by space derivatives using the governing equation~\eqref{eq:secondOrderWaveEquation} to give
\ba
   \Dpt\Dmt\, u  = \Lc u +  \frac{\dt^2}{12} \Lc^2 u + \frac{\dt^4}{360} \Lc^3 u + \cdots
     . \label{eq:taylorExpansionLc}
\ea
To form a $p^{\rm th}$ order accurate in time scheme, the expansion~\eqref{eq:taylorExpansionLc} 
is truncated to $p/2$ terms, and the spatial operators in the resulting truncated expansion are discretized with a 
time-weighted average of three time levels.  For example, second-order accurate explicit or implicit three-level ME schemes take the form
\ba
   \Dpt\Dmt u(\xv,t) = \Lc_{2,h} \Big(  \alpha_2 u(\xv,t+\dt) +\beta_2  u(\xv,t) + \gamma_2 u(\xv,t-\dt) \Big),
   \label{eq:ME2} 
\ea
where $\Lc_{2,h}$ is a second-order accurate approximation of $\Lc$ on a grid with representative grid spacing~$h$.  The coefficients $(\alpha_2,\beta_2,\gamma_2)$ are the weights in the time-weighted average of $u(\xv,t)$.  The fourth-order accurate scheme has the form
\ba
   \Dpt\Dmt u(\xv,t) 
      & = \Lc_{4,h} \Big(  \alpha_2 u(\xv,t+\dt) +\beta_2  u(\xv,t) + \gamma_2 u(\xv,t-\dt) \Big) \nonumber \\
       & \qquad - \dt^2 \Lc_{2,h}^2 \Big(  \alpha_4 u(\xv,t+\dt) +\beta_4 u(\xv,t) + \gamma_4 u(\xv,t-\dt) \Big) ,
    \label{eq:ME4}  
\ea
where $\Lc_{4,h}$ is a fourth-order accurate approximation of $\Lc$ and $(\alpha_4,\beta_4,\gamma_4)$ are coefficients involved in the time-average of the correction term.
Higher-order accurate schemes for $p=6,8,\ldots$ can be defined in a similar way but for this article we only
consider schemes for $p=2$ and~$4$.
The explicit ME schemes we use have $\alpha_{2m}=\gamma_{2m}=0$ and $\beta_{2m}\ne0$ for $m=1$ and~$2$, while the implicit schemes have $\alpha_{2m}\ne 0$ for $m=1$ or~$2$.

Truncation error analysis can be used to determine the constraints on 
the parameters $\alpha_{2m}$, $\beta_{2m}$ and $\gamma_{2m}$ for $p^{th}$ order accuracy. 
The truncation error of the $p=2$ scheme in~\eqref{eq:ME2}, denoted by $\tau_2(\xv,t)$, 
is 
\ba
  \tau_2(\xv,t) = \big( 1 - (\alpha_2 + \beta_2 + \gamma_2)\big) \,  \Lc u
            - \dt \big(\alpha_2-\gamma_2\big)\, \partial_t \Lc u + O(\dt^2 + h^2).
\ea
For second-order accuracy in $\dt$ and $h$ we take 
\bse
\label{eq:ME2accuracy}
\ba
  \alpha_2 + \beta_2 + \gamma_2 & = 1,\\
  \alpha_2-\gamma_2 & = 0.
\ea
\ese 
These two conditions involving the three parameters $(\alpha_2,\beta_2,\gamma_2)$  give a single-parameter family of second-order accurate schemes discussed further in Section~\ref{sec:IMEonCartesian}.
Note that the condition $\alpha_2=\gamma_2$ implies that the schemes are symmetric in time which implies the
schemes are time reversible.  A similar analysis of the $p=4$ scheme in~\eqref{eq:ME4} leads to the conditions
\bse
\label{eq:ME4accuracy}
\ba
  \alpha_2 + \beta_2 + \gamma_2 & = 1,\\
  \alpha_2-\gamma_2 & = 0,\\
  \alpha_4-\gamma_4 & = 0,\\
   \half(\alpha_2 + \gamma_2)-(\alpha_4 +\gamma_4)-\beta_4 & = \frac{1}{12}.
\ea
\ese
The four constraints in~\eqref{eq:ME4accuracy} involving the six parameters $(\alpha_{2m},\beta_{2m},\gamma_{2m})$, $m=1,2$, implies a two-parameter family of fourth-order accurate schemes as discussed further in Section~\ref{sec:IMEonCartesian}.
Note that these $p=4$ schemes are also symmetric in time.
Choices of the parameters that lead to stable schemes for $p=2$ and~$4$ are discussed in Section~\ref{sec:stabilityAnalysis}.

\section{Implicit modified equation (IME) schemes on Cartesian grids} \label{sec:IMEonCartesian}

In order to analyze the proposed IME schemes in more detail we introduce a spatial approximation 
for Cartesian grids in Section~\ref{sec:spatialApprox}.
This allows us to show the form of the fully discrete schemes and to 
perform a von Neumann stability analysis 
in Section~\ref{sec:stabilityAnalysis}.

\subsection{Spatial approximation on Cartesian grids} \label{sec:spatialApprox}

Let the domain $\Omega=[0,2\pi]^{\nd}$ be a box in $\nd$ dimensions discretized
with a Cartesian grid with $N_d$ grid points in each direction.
Denote the grid points as
\ba
   \xv_{\jv} = [j_1 h_1,..., j_{\nd}h_{\nd}]^T
\ea
for multi-index $\jv\in\Integers^{\nd}$ and grid spacings $h_d=2\pi/N_d$.
Let $U_{\jv}^n\approx u(\xv_{\jv},t^n)$ be elements of a grid function at time $t^n=n\dt$.
Define the usual divided difference operators to be
\ba
  D_{+\xv_d} U_{\jv}^n \defeq \frac{U_{\jv + \ev_d}^n - U_{\jv}^n}{h_d}, \; \; \; D_{-\xv_d} U_{\jv}^n \defeq \frac{U_{\jv}^n - U_{\jv - \ev_d}^n}{h_d}, \; \; \; D_{0\xv_d} U_{\jv}^n \defeq \frac{U_{\jv + \ev_d}^n - U_{\jv - \ev_d}^n}{2h_d}, \label{eq:ddopts}
\ea
where $\ev_d\in\Real^{\nd}$ is the unit vector in direction $d$  (e.g.~$\ev_2=[0,1,0]$). 

The compact $p^{th}$ order accurate discretization of the operator $\Lc$ can be written in the form
\ba
  \Lc_{p,h} \defeq c^2 \,  \sum_{m=0}^{p/2-1} \kappa_m \left[ \sum_{d=1}^{\nd}  h_d^{2m}(D_{+\xv_d}D_{-\xv_d})^{m+1} \right]. \label{eq:LcDiscretization}
\ea
where, for example, $\kappa_0=1$, $\kappa_1=-1/12$,  $\kappa_2=1/90$ and $\kappa_3=-1/560$.
The compact second-order accurate approximation to $\Lc^2$ is just the square of the second-order accurate
approximation $\Lc_{2,h}$
\ba
    \Lc_{2,h}^2 = \Lc_{2,h} \, \Lc_{2,h} .
\ea
Although not used here, note that the compact fourth-order accurate approximation to $\Lc^2$ is not
the square of $\Lc_{4,h}$ as $(\Lc_{4,h})^2$ has a wider stencil than is needed~\cite{lcbc2022}.

Given the accuracy requirements~\eqref{eq:ME2accuracy} 
we write the fully discrete second-order accurate ME scheme (denoted by IME2) 
in terms of a single free parameter $\alpha_2$, 
\ba
  &  \Dpt\Dmt U_\jv^n = \Lc_{2,h} \Big( \alpha_2 \, U_\jv^{n+1} + (1-2\alpha_2) \,U_\jv^{n}  + \alpha_2 \, U_\jv^{n-1} \Big).
         \label{eq:ME2discrete}   
\ea
Note that larger values of $\alpha_2$ correspond to schemes that are \textsl{more implicit} with $\alpha_2=0$ being the
explicit EME2 scheme.
Similarly the fully discrete fourth-order accurate ME scheme (denoted by IME4) 
involves two two free parameters $\alpha_2$ and $\alpha_4$, 
\ba
 \Dpt\Dmt U_\jv^{n} 
      & = \Lc_{4,h} \Big( \alpha_2 \, U_\jv^{n+1} + (1-2\alpha_2) \,U_\jv^{n}  + \alpha_2 \, U_\jv^{n-1} \Big) \nonumber \\
       & \qquad - \dt^2 \Lc_{2,h}^2 
       \Big(  \alpha_4 \, U_\jv^{n+1} + ( \alpha_2 - 2\alpha_4-\f{1}{12} ) \,U_\jv^{n}  + \alpha_4 \, U_\jv^{n-1}    \Big) .
  \label{eq:ME4discrete} 
\ea
Larger values of $\alpha_2$ and $\alpha_4$ correspond to schemes that are more implicit with $\alpha_2=\alpha_4=0$
being the explicit EME4 scheme.

\subsection{Stability analysis of the implicit modified equation (IME) schemes} \label{sec:stabilityAnalysis}

The stability of the IME schemes~\eqref{eq:ME2discrete} and~\eqref{eq:ME4discrete} is now studied
using von Neumann analysis, assuming solutions that are periodic in space.
Von Neumann analysis expands the solutions in a Fourier series and determines conditions so that all Fourier
modes remain stable.
There are numerous definitions for stability, but for our purposes here we make the following definition:
\begin{definition}[Stability]
  \label{def:stability}
  A numerical scheme for the wave equation is stable if there are no Fourier modes with non-zero wave-number, $\kv\ne\zerov$, whose
  magnitude grow in time.  For the zero wave-number, $\kv=\zerov$, case the linear in time mode given by $u = c_0 + c_1 t$ 
  for constants $c_0$ and $c_1$, is permitted
  since this is an exact solution to the wave equation.
\end{definition}
The explicit ME schemes (with $\alpha_2=\alpha_4=0$) are known to be CFL-one stable (at least for $p=2,4,6$), 
meaning stable for 
\ba
  c^2 \dt^2\sum_{d=1}^{n_d} \frac{1}{h_d^2} <1.
\ea

For implicit ME schemes we are generally interested in unconditional stability, that is stability for any $\dt>0$.
The constraint on $\alpha_2$ for unconditional (von Neumann) stability of the second-order accurate IME2 scheme is summarized by the following theorem.
\begin{theorem}[IME2 Stability]
\label{thm:stabilityIME2}
 The IME2 scheme~\eqref{eq:ME2discrete} is unconditionally stable on a periodic domain provided
 \ba
     \alpha_2 \ge \f{1}{4}.       
\ea
\end{theorem}
The constraints on $\alpha_2$ and $\alpha_4$ for unconditional stability of the fourth-order accurate IME4 scheme are summarized by the following theorem.
\begin{theorem}[IME4 Stability]
\label{thm:stabilityIME4}
 The IME4 scheme~\eqref{eq:ME4discrete} is unconditionally stable on a periodic domain provided
\bse
\label{eq:stabilityConditionIME4}
\ba
    & \alpha_2 \ge \f{1}{12}, \\
    & \alpha_4 \ge 
       \begin{cases}
           \f{1}{4} \alpha_2 - \f{1}{48},                                  &   \text{when $\alpha_2\ge \f{1}{4}$}, \\
           \f{1}{4} \alpha_2 - \f{1}{48} + \f{8}{9} (\f{1}{4}-\alpha_2)^2, &  \text{when $\f{1}{12} \le \alpha_2 \le \f{1}{4}$ }.
       \end{cases} 
\ea
\ese
\end{theorem}
The proofs of Theorems~\ref{thm:stabilityIME2} and~\ref{thm:stabilityIME4} 
are given in~\ref{sec:stabilityProof2IME2} and~\ref{sec:stabilityProof2IME4}, respectively.

\section{Upwind dissipation and implicit modified equation (IME-UW) schemes} \label{sec:upwindDissipation}

The EME and IME schemes described in previous sections have no dissipation and are neutrally stable. 
As a result, perturbations to the schemes, such as with variable coefficients or 
when the schemes are applied on overset grids may lead to instabilities. For single curvilinear grids, stable schemes
can be defined using special discretizations such 
as summation by parts (SBP)~\cite{Strand1994,Olsson1995,Mattsson-Nordstrom-2004,AppeloPetersson2009,DuruKreissMattsson2014} 
methods or the schemes described in~\cite{max2006b}.
Overset grids are a greater challenge and in this case we add dissipation for stability.
Upwind dissipation for the wave equation in second-order form was first developed in~\cite{sosup2012} and
applied to Maxwell's equations in~\cite{mxsosup2018}. An optimized version of the upwind dissipation was
developed in~\cite{ssmx2023} and it is this optimized version that we use here as a template for the IME scheme.

We first consider upwind dissipation for the explicit ME scheme to establish our basic approach and to introduce some notation.
For the explicit scheme, dissipation is added in a predictor-corrector fashion,
\bse
\label{eq:EME-UW}
\ba
   & U_\jv^{(0)}  = 2 U_\jv^n - U_\jv^{n-1}  + \dt^2 \, \Lv_p\, U_\jv^n , \\
   & U_\jv^{n+1} = U_\jv^{(0)}  - \nu_p \, \dt^2 \, \Qv_p \left[ \f{U_\jv^{(0)} - U_\jv^{n-1}}{2\dt} \right],  \label{eq:upwindStep}
\ea
\ese
where $\Lv_p$ denotes the (full) spatial operator for the $p^{\rm th}$-order accurate scheme,
$\nu_p$ is an upwind dissipation parameter, and 
$\Qv_p$ is a dissipation operator, which on a Cartesian grid takes the form
\ba
  & \Qv_p \eqdef \sum_{d=1}^{\nd} \f{c}{h_d} \left[ - \Delta_{+x_d}\Delta_{-x_d}  \right]^{p/2+1} , \label{eq:Qdef}
\ea
where $\Delta_{\pm x_d}$ are undivided difference operators corresponding to the divided difference operators defined in~\eqref{eq:ddopts}.  Note that the dissipation operator $\Qv_p$ has a stencil of width $p+3$ compared to the stencil width of $p+1$ for $\Lv_p$.
The wider stencil for the dissipation reflects the upwind character of the operator~\cite{sosup2012}.
Also note that the addition of the upwind dissipation does not change the order of accuracy of the scheme.

The dissipation operator $\Qv_p$ in~\eqref{eq:upwindStep} acts on an approximation of $\partial_t u(\xv_\jv,t^n)$.  The treatment of this approximation in the predictor-corrector scheme in~\eqref{eq:EME-UW} ensures that the scheme, with dissipation, remains explicit and p$^{\rm th}$-order accurate.  For implicit ME schemes, there is more flexibility in the treatment of this approximation.  Two approaches are described in the next subsections.

\subsection{Monolithic upwind dissipation for IME schemes (IME-UW)}

Upwind dissipation for the implicit ME schemes can be added directly into the single step update (denoted as the IME-UW scheme) as
\ba
 & \Dpt\Dmt U_\jv^n = \Lv_{\alpha p}( U_\jv^{n+1}, U_\jv^n,U_\jv^{n-1}) - \nu_p \, \Qv_p \left[ \f{U_\jv^{n+1} - U_\jv^{n-1}}{2\dt} \right].
   \label{eq:IME-UW}
\ea
Here $\Lv_{\alpha p}$ denotes the spatial part of the IME scheme as given in~\eqref{eq:ME2discrete} and~\eqref{eq:ME4discrete} for some choice of the parameters $\alpha_2$ and $\alpha_4$.
A von Neumann stability analysis for a Cartesian grid leads to the following result.
\begin{theorem} \label{th:IME-UW}
  The IME-UW schemes~\eqref{eq:IME-UW} for $p=2,4$ on a periodic Cartesian grid 
  are unconditionally stable for any $\nu_p>0$ provided $\alpha_2$ satisfies the conditions of
  Theorem~\ref{thm:stabilityIME2}, for $p=2$, or $\alpha_2$ and $\alpha_4$ satisfy the conditions for 
  Theorem~\ref{thm:stabilityIME4} for $p=4$.
\end{theorem}
The proof of this theorem is given in~\ref{sec:stabilityIME-UW}.
The monolithic upwind dissipation allows for any $\nu_p>0$ and there are various possible strategies
for choosing this value~\cite{AllisonCarson2023}.

\subsection{Predictor-corrector upwind dissipation for IME schemes (IME-UW-PC)} \label{sec:IME-UW-PC-MUSTA}

One disadvantage of the upwind scheme~\eqref{eq:IME-UW} is that the dissipation operator 
changes the implicit matrix, increasing the stencil size. This may increase the cost of the implicit solve and 
be undesirable, if for example,
one wants to use an existing multigrid solver not designed for this special matrix.
Dissipation can be added to the IME scheme is a separate explicit correction step as in~\eqref{eq:EME-UW}.
Allowing for for multiple corrections leads to the implicit-predictor explicit-corrector upwind scheme (IME-UW-PC)
\bse
\label{eq:IME-UW-PC-MUSTA}
\ba
 & \f{U_\jv^{(0)} -2 U_\jv^n + U_\jv^{n-1}}{\dt^2 } = \Lv_{\alpha p}\, (U_\jv^{(0)}, U_\jv^n, U_\jv^{n-1}) , \\
 & U_\jv^{(k)} = U_\jv^{(k-1)}  - \nu_p \, \dt^2 \, \Qv_p \left[ \f{U_\jv^{(k-1)} - U_\jv^{n-1}}{2\dt} \right] , \qquad k=1,2,\ldots,\Nuc , \\ 
 & U_\jv^{n+1} = U_\jv^{(\Nuc)} .
\ea
\ese
where $\Nuc$ denotes the number of upwind correction steps.
The sequence of corrections in~\eqref{eq:IME-UW-PC-MUSTA} can be combined and written succinctly as 
\ba
   U_\jv^{n+1} &=  \Rv_p^{\Nuc} U_\jv^{(0)} + (I- \Rv_p^{\Nuc} ) U_\jv^{n-1}, 
\ea
where
\ba
    \Rv_p \eqdef I - \f{\nu_p \dt}{2} \Qv_p.
\ea
The conditions on $\nu_p$ for stability are specified in the following theorem.
The theorem covers the cases of using an implicit or an explicit ME predictor.
\begin{theorem} \label{th:IME-UW-PC-MUSTA}
   The upwind predictor-corrector 
   scheme~\eqref{eq:IME-UW-PC-MUSTA} is stable on the periodic domain 
  provided the non-dissipative predictor scheme (explicit or implicit) is stable and provided
  \ba
    0 \le \nu_p < \f{\sigma_\Nuc}{ 2^{p+1} \sum_{d=1}^{\nd} \lambda_{x_d}}.  \label{eq:nupIME-UW-PC}
  \ea
  where $\sigma_\Nuc=2$ for $\Nuc$ even and $\sigma_\Nuc=1$ for $\Nuc$ odd, and 
  where $\lambda_{x_d}$ is the CFL parameter in coordinate direction $d$, 
  \ba
     \lambda_{x_d} \eqdef \f{c\dt}{h_d}.
  \ea
\end{theorem}
The proof of Theorem~\ref{th:IME-UW-PC-MUSTA} is given in~\ref{sec:stabilityIME-UW-PC}.
\mni
In practice a reasonable choice might be 
\ba
  \nu_p  = \f{\safetyFactor \, \sigma_\Nuc }{ 2^{p+1} \sum_{d=1}^{\nd} \lambda_{x_d}} , \label{eq:nuprecommended}
\ea
where $\safetyFactor\in(0,1)$ is a safety factor.

Note from~\eqref{eq:nupIME-UW-PC} that the coefficient of dissipation, $\nu_p$,
decreases as the CFL parameter increases, and thus less dissipation is added as the CFL number increases.
Thus, for large CFL it may become necessary to
use more than one correction step and in general taking $\Nuc$ to be the integer larger than the CFL, $\Nuc = \lceil \CFL \rceil$,
seems appropriate for difficult cases
such as when solving on an overset grid.


\section{Spatially partitioned implicit-explicit (SPIE) ME schemes} \label{sec:spieSchemes}

The IME and EME schemes can be combined in a spatially partitioned manner.
For overset grids, the IME scheme is used on certain components grids while the EME scheme
is applied on all other component grids. A typical strategy is to employ the EME scheme on background Cartesian grids and any curvilinear grids with grids spacings close to a nominal Cartesian value, and then use the IME scheme on any curvilinear component grids that have a minimum grid spacing that is relatively small as compared to the nominal value.

\newcommand{\Gce}{\Gc_{\rm e}}
\newcommand{\Gci}{\Gc_{\rm i}}
\subsection{Formulation of the SPIE scheme}

Let $\Gce$ denote the set of grids that use explicit time-stepping (e.g.~Cartesian grids), and
let $\Gci$ denote the set of grids that use implicit time-stepping (e.g.~curvilinear grids).
The SPIE algorithm consists  of the following three stages:

\mni
\bse
\label{eq:SPIEscheme}
Stage 1. Update explicit grids,
\ba
 & \f{U_{g,\,\jv}^{(0)} -2 U_{g,\,\jv}^n + U_{g,\,\jv}^{n-1}}{\dt^2 } = \Lv_{p}\, U_{g,\,\jv}^n ,   \quad g\in\Gc_{\rm e} , \label{eq:SPIEschemeExplicit}
\ea
Stage 2. Update implicit grids, interpolating from the solution on explicit grids from Stage 1,
\ba
 & \f{U_{g,\,\jv}^{(0)} -2 U_{g,\,\jv}^n + U_{g,\,\jv}^{n-1}}{\dt^2 } = \Lv_{\alpha p}\, (U_{g,\,\jv}^{(0)}, U_{g,\,\jv}^n, U_{g,\,\jv}^{n-1}) ,    \quad g\in\Gc_{\rm i}  \label{eq:SPIEschemeImplicit}
\ea
Stage 3. Add dissipation to all grids. For example, if all grids use a predictor-corrector upwind formulation, then use
\ba
 & U_{g,\,\jv}^{n+1} = U_{g,\,\jv}^{(0)}  + \nu_p \, \dt^2 \, \Qv_p \left[ \f{U_{g,\,\jv}^{(0)} - U_{g,\,\jv}^{n-1}}{2\dt} \right]. 
    \quad g\in \Gc_{\rm i} \cup \Gc_{\rm e}.  \label{eq:SPIEschemeDissipation}
\ea
\ese
Multiple upwind correction steps can also be used as discussed in Section~\ref{sec:IME-UW-PC-MUSTA}.


\subsection{GKS stability analysis of a model problem for the SPIE scheme}\label{sec:gksanalysis}

This section investigates the stability of the SPIE scheme for a one-dimensional overset grid.
Normal mode (GKS) analysis~\cite{GKSI,GKSII} is used to show that the second-order accurate SPIE scheme 
is stable when solving the wave equation on an overset grid for a one-dimensional infinite domain.
Matrix stability analysis on a finite domain is presented later in Section~\ref{sec:matrixStability} and the result
of the matrix analysis supports the conclusions of the GKS analysis discussed in this section.

{

\newcommand{\tikzcircle}[2][red,fill=red]{\tikz[baseline=-0.5ex]\draw[#1,radius=#2] (0,0) circle ;}%
\begin{figure}[hbt!]
  \centering
  \begin{tikzpicture}
    \useasboundingbox (0,.7) rectangle (11,2.8);

     \begin{scope}[yshift=-0.7cm]
     \draw[ultra thick,red,<-] (0,3) -- (6,3);
     \foreach \x in {1,2,3,4,5}
       \draw[ultra thick,red] (\x,3cm-4pt) -- (\x,3cm +4pt); 
 
     \draw (6,3) node {\tikzcircle[red,very thick, fill=red!20]{3pt}};
 
     \draw (3,3) node[anchor=south,yshift=2pt] {\labelFont$x_{\Lss,-3}$} node[anchor=north,yshift=-2pt]{\labelFont$U_{\Lss,-3}$};
     \draw (4,3) node[anchor=south,yshift=2pt] {\labelFont$x_{\Lss,-2}$} node[anchor=north,yshift=-2pt]{\labelFont$U_{\Lss,-2}$};
     \draw (5,3) node[anchor=south,yshift=2pt] {\labelFont$x_{\Lss,-1}$} node[anchor=north,yshift=-2pt]{\labelFont$U_{\Lss,-1}$};
     \draw (6,3) node[anchor=south,yshift=2pt] {\labelFont$x_{\Lss, 0}$} node[anchor=north,yshift=-2pt]{\labelFont$U_{\Lss, 0}$};
     \draw (2,3.35) node[]{$\dots$};

    \end{scope}

    \begin{scope}[xshift=0cm]
      \draw[ultra thick,blue,->] (5,1) -- (11,1);
       \foreach \x in {6,7,8,9,10}
        \draw[ultra thick,blue] (\x,1cm-4pt) -- (\x,1cm +4pt); 
  
      \draw (5,1) node {\tikzcircle[blue,very thick, fill=blue!20]{3pt}};
      \draw (5,1) node[anchor=north,yshift=-2pt] {\labelFont$x_{\Rss, 0 }$} node[anchor=south,yshift=2pt]{\labelFont$U_{\Rss, 0}$};
      \draw (6,1) node[anchor=north,yshift=-2pt] {\labelFont$x_{\Rss, 1 }$} node[anchor=south,yshift=2pt]{\labelFont$U_{\Rss, 1}$};
      \draw (7,1) node[anchor=north,yshift=-2pt] {\labelFont$x_{\Rss, 2 }$} node[anchor=south,yshift=2pt]{\labelFont$U_{\Rss, 2}$};
      \draw (8,1) node[anchor=north,yshift=-2pt] {\labelFont$x_{\Rss, 3 }$} node[anchor=south,yshift=2pt]{\labelFont$U_{\Rss, 3}$};

      \draw (9,1) node[yshift=-9pt]{$\dots$};
    \end{scope}
%
  \end{tikzpicture}
  \caption{One-dimensional overset grid used to assess the stability of the SPIE scheme.
    The explicit scheme is used on the left grid and the implicit scheme is used on the right.
  Interpolation points are marked as circles.}
  \label{fig:oversetGrid1d}
\end{figure}
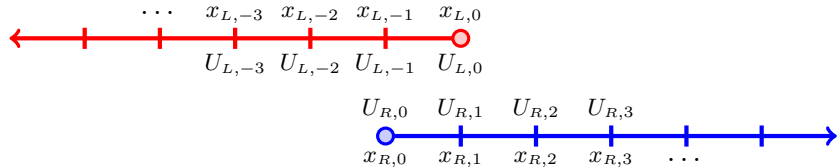 
}

Consider the one-dimensional overlapping grid for $\Omega=(-\infty,\infty)$ shown in Figure~\ref{fig:oversetGrid1d}.
Let $x_{\Lss,j} = (j+1) h$ and $x_{\Rss,j}=j h$ 
denote the grid points for the left and right grids, respectively.
Let $U_{s,j}$, for $s=L,R$, denote the discrete solutions on the two grids.
The grids overlap by a distance $h$ and the solution is interpolated at the interpolation points shown in Figure~\ref{fig:oversetGrid1d}.
The second-order accurate SPIE scheme is used. 
The explicit (EME2) scheme is applied on the left grid and the implicit (IME2) scheme is applied on the right grid.
The discrete equations are 
\bse
\label{eq:pieGKS}
\bat
  & \dpmt U_{\Lss,j}^n =  c^2\dpmx U_{\Lss,j}^n,           
          && \quad j<0,                                \label{eq:pieGKSa} \\
  & \dpmt U_{\Rss,j}^n = c^2\dpmx \left[ \alpha_2 \, U_{\Rss,j}^{n+1} + (1-2\alpha_2) \, U_{\Rss,j}^n + \alpha_2 \, U_{\Rss,j}^{n-1}\right], 
         &&\quad j>0,         \label{eq:pieGKSb} \\
  & |U_{\Lss,j}^n|<\infty , && \quad j\to-\infty,  \label{eq:pieGKSc}\\
  & |U_{\Rss,j}^n|<\infty , && \quad j\to\infty, \label{eq:pieGKSd}\\
  & U_{\Lss,0}^n = U_{\Rss,1}^n, \label{eq:pieGKSe}\\
  & U_{\Rss,0}^n = U_{\Lss,-1}^n,\label{eq:pieGKSf}
\eat
\ese
where we have assumed that the solution is bounded as $|j|\rightarrow \infty$.
Each individual scheme is assumed to be stable and so we take $\lambda=c\dt/h <1$ and
 $\alpha_2 \ge 0$.  Note that the IME2 scheme is unconditionally stable for $\alpha_2 \ge 1/4$, but $\alpha_2\ge0$ is sufficient when $\lambda<1$.

Before proceeding with the stability analysis, it is useful to first state the following lemma related to the stability 
of the EME2 and IME2 schemes for the Cauchy problem. 
\begin{lemma}
\label{thm:gksLemmaIME2}
   Suppose $a$ is a root of  quadratic equation,
   \bse
    \label{eq:lemmaQuadratic}
   \ba
      a^{2} -2 b \, a + 1 = 0, \label{eq:lemmaIME2quadratic}
   \ea
   where $b$ is defined in terms of some $\kappa\in\Complex$ by 
   \ba
     & b \eqdef \frac{1 + (\half-\alpha_2)\, \lambda^2 \,(\kappa -2 + \kappa^{-1})}
                                          {1 - \alpha_2 \, \lambda^2 \,(\kappa -2 + \kappa^{-1})} .
   \ea
   \ese
   Then,  $|\kappa|=1$ implies $|a|=1$,  when (i) $\lambda<1$ and $\alpha_2 \ge0$ or when (ii) $\lambda>0$ and $\alpha_2 \ge 1/4$.
   
\end{lemma}
\noindent
The proof of Lemma~\ref{thm:gksLemmaIME2} is given in~\ref{sec:GKSlemma}.


\medskip
We are now ready to prove the main theorem of this section.
\begin{theorem}
\label{thm:gksSPIE}
   The SPIE scheme in~\eqref{eq:pieGKS}, for the one-dimensional infinite domain overset grid, 
   has no unstable solutions provided $\lambda<1$ and $\alpha_2\ge 0$. 
\end{theorem}
\begin{proof}
  We look for unstable mode solutions of the form
   \bse
   \label{eq:SPIEansatz}
   \ba   
      U_{\Lss,j}^n & = a^n \kappa_\Lss^j,\\
      U_{\Rss,j}^n & = a^n \kappa_\Rss^j , 
   \ea
   \ese
   for some $a\in\Complex$ with $|a|>1$. Note that the same amplification factor $a$ must appear in both the left and right grid functions in order to match the interpolation equations~\eqref{eq:pieGKSe} and~\eqref{eq:pieGKSf}.
   Substituting the ansatz~\eqref{eq:SPIEansatz} into~\eqref{eq:pieGKSa} and~\eqref{eq:pieGKSb} implies,
   \bse
   \label{eq:SPIEquadratics}
   \ba
     & a^{1} -2 a + a^{-1} = \lambda^2 (\kappa_\Lss -2 + \kappa_\Lss^{-1} ) , \\
     & a^{1} -2 a + a^{-1} = \lambda^2 (\kappa_\Rss -2 + \kappa_\Rss^{-1} ) \big( \alpha_2 a +  (1-2\alpha_2) + \alpha_2 a^{-1} \big).
   \ea
   \ese
   The equations~\eqref{eq:SPIEquadratics}  can be written as quadratics for $\kappa_s$, 
   \bse
   \ba
      \kappa_s^{2} -2 b_s \kappa_s + 1  = 0, 
   \ea
   \ese
   for some $b_s$, $s=L,R$,
   with roots denoted by $\kappa_{s,\pm}$. 
   The general solutions for the left and right sides then take the form
   \bse
   \ba
    &  U_{\Lss,j}^n = a^n (c_{+}\kappa_{\Lss+}^j + c_{-}\kappa_{\Lss+}^{-j}),\\
    &  U_{\Rss,j}^n = a^n (d_{+}\kappa_{\Rss+}^j + d_{-}\kappa_{\Rss+}^{-j}),
   \ea
   \ese
   where $c_{\pm}$ and $d_{\pm}$ are constants.
   Note that both equations in~\eqref{eq:SPIEquadratics} are of the form~\eqref{eq:lemmaQuadratic} of Lemma~\ref{thm:gksLemmaIME2}.
   Since we have assumed $|a|>1$, it follows from Lemma~\ref{thm:gksLemmaIME2} that the roots $\kappa_{s,\pm}$ cannot have magnitude equal to one.
   Since the product of the roots $\kappa_{s,\pm}$ is one, we can 
   therefore, without loss of generality, take $|\kappa_{s,+}|<1$, $s=L,R$.

   The boundedness conditions~\eqref{eq:pieGKSc} and~\eqref{eq:pieGKSd} at infinity imply $c_+=0$ and $d_-=0$, reducing the solutions to 
   \bse
   \ba
     &  U_{\Lss,j}^n  = a^n c_{-}\kappa_{\Lss+}^{-j} ,\\
     &  U_{\Rss,j}^n  = a^n  d_{+}\kappa_{\Rss+}^{j}.
   \ea
   \ese
   Applying the interpolation conditions~\eqref{eq:pieGKSe} and~\eqref{eq:pieGKSf} gives
   \bse
   \ba
    &  c_- = d_{+}\kappa_{\Rss+}, \\
    &  d_+ = c_{-}\kappa_{\Lss+},
   \ea
   \ese
   which implies, assuming $c_{-} \ne 0$ and $d_+\ne 0$, that 
   \ba
      \kappa_{\Lss+} \kappa_{\Rss+} =1.
   \ea
   This last equation cannot hold since $|\kappa_{\Lss+} \kappa_{\Rss+}|<1$.
   Therefore, only the trivial solution remains, thus yielding no unstable solutions with $|a|>1$.
\end{proof}

\section{Overset grids, implicit first step, and implicit solvers} \label{sec:oversetGrids}

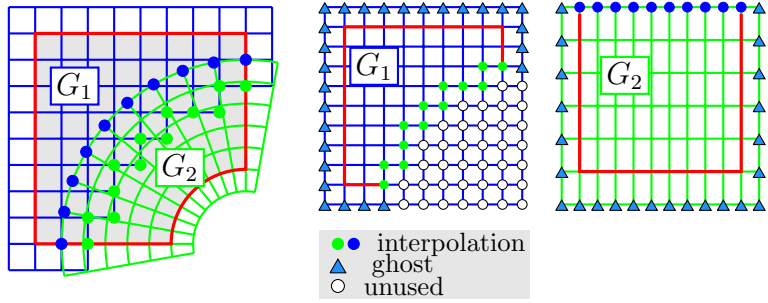
\begin{figure}[hbt]
\begin{center}
\begin{tikzpicture}[scale=.7]
\useasboundingbox (.75,1.25) rectangle (15.5,6);  
%
\begin{scope}[xshift=1cm,yshift=1cm]
\fill[black!10!white,xshift=.5cm,yshift=.5cm] (0,0) -- (2.583333,0) arc (180:90:1.416667) -- (4.,4.) -- (0,4.) -- (0,0);
\draw[-,thick,blue,yshift=.0 cm] 
   \foreach \x/\y in {1.5/0,1.5/.5,2/1,2/1.5,2.5/2,3/2.5,4/3,5/3.5,5/4,5/4.5,5/5}{ (0,\y) -- (\x,\y) }
   \foreach \x/\y in {0/0,.5/0,1/0,1.5/0,2/1,2.5/2,3/2.5,3.5/3,4/3,4.5/3.5,5/3.5}{ (\x,\y) -- (\x,5) };
  \begin{scope}[xshift=4.5cm,yshift=0.5cm]
    \draw[thick,green] \foreach \r in {1.000000,1.416667,1.833333,2.250000,2.666667,3.083333,3.500000}{ (0,\r) arc (90:190:\r)  (0,\r) arc (90:80:\r) };
    \draw[thick,green]
     (0.173648,0.984808)  -- (0.607769,3.446827)
     (0.000000,1.000000)  -- (0.000000,3.500000)
     (-0.173648,0.984808) -- (-0.607769,3.446827)
     (-0.342020,0.939693) -- (-1.197071,3.288924)
     (-0.500000,0.866025) -- (-1.750000,3.031089)
     (-0.642788,0.766044) -- (-2.249757,2.681156)
     (-0.766044,0.642788) -- (-2.681156,2.249757)
     (-0.866025,0.500000) -- (-3.031089,1.750000)
     (-0.939693,0.342020) -- (-3.288924,1.197071)
     (-0.984808,0.173648) -- (-3.446827,0.607769)
     (-1.000000,0.000000) -- (-3.500000,0.000000)
     (-0.984808,-0.173648) -- (-3.446827,-0.607769);
  \end{scope}
  \draw[very thick,red,xshift=.5cm,yshift=.5cm] (0,0) -- (2.583333,0) arc (180:90:1.416667) -- (4.,4.) -- (0,4.) -- (0,0);
%
   \filldraw[green] (1.5,.5)  circle (3pt)
                 (1.5,1 )  circle (3pt)
                 (2  ,1 )  circle (3pt)
                 (2 ,1.5)  circle (3pt)
                 (2 , 2 )  circle (3pt)
                 (2.5,2 )  circle (3pt)
                 (2.5,2.5) circle (3pt)
                 (3 , 2.5) circle (3pt)
                 (3  ,3 )  circle (3pt)
                 (3.5,3 )  circle (3pt)
                 (4  ,3. ) circle (3pt)
                 (4  ,3.5) circle (3pt)
                 (4.5,3.5) circle (3pt);
%
  \begin{scope}[xshift=4.5cm,yshift=0.5cm]
      \filldraw[blue]
       (0.000000,3.500000)    circle (3pt)
       (-0.607769,3.446827)   circle (3pt)
       (-1.197071,3.288924)  circle (3pt) 
       (-1.750000,3.031089)  circle (3pt) 
       (-2.249757,2.681156)  circle (3pt) 
       (-2.681156,2.249757)  circle (3pt) 
       (-3.031089,1.750000)  circle (3pt) 
       (-3.288924,1.197071)  circle (3pt) 
       (-3.446827,0.607769)  circle (3pt) 
       (-3.500000,0.000000)  circle (3pt);
  \end{scope}
   \draw (1.25,3.5) node[thick,draw=blue,fill=white,inner sep=2pt] {\large$G_1$};
   \draw (3.25,1.95) node[thick,draw=green,fill=white,inner sep=2pt] {\large$G_2$};
\end{scope}
%
\definecolor{ghostColour}{named}{DodgerBlue}
\newcommand{\mytrix}{(\x,-.15) -- ++(.3,0) -- ++(-.15,.26) -- (\x,-.15)}
\newcommand{\mytriy}{(-.15,\y) -- ++(.3,0) -- ++(-.15,.26) -- (-.15,\y)}
\begin{scope}[xshift=7cm,yshift=2.25cm,scale=.75]
\draw[-,thick,blue,yshift=.0 cm] 
   \foreach \x in {0,.5,...,5}{ (\x,0) -- (\x,5) }
   \foreach \y in {0,.5,...,5}{ (0,\y) -- (5,\y) };
  \draw[very thick,red,xshift=.5cm,yshift=.5cm] (1.,0) -- (.0,0) -- (.0,4.) -- (4.,4.) -- (4.,3.);
   \filldraw[green] (1.5,.5)  circle (3pt)
                 (1.5,1 )  circle (3pt)
                 (2  ,1 )  circle (3pt)
                 (2 ,1.5)  circle (3pt)
                 (2 , 2 )  circle (3pt)
                 (2.5,2 )  circle (3pt)
                 (2.5,2.5) circle (3pt)
                 (3 , 2.5) circle (3pt)
                 (3  ,3 )  circle (3pt)
                 (3.5,3 )  circle (3pt)
                 (4  ,3. ) circle (3pt)
                 (4  ,3.5) circle (3pt)
                 (4.5,3.5) circle (3pt);
  \filldraw[fill=white,draw=black]  \foreach \x in {2,2.5,...,5}{ (\x,.0) circle (3.5pt) };
  \filldraw[fill=white,draw=black]  \foreach \x in {2,2.5,...,5}{ (\x,.5) circle (3.5pt) };
  \filldraw[fill=white,draw=black]  \foreach \x in {2.5,3,...,5}{ (\x,1.) circle (3.5pt) };
  \filldraw[fill=white,draw=black]  \foreach \x in {2.5,3,...,5}{ (\x,1.5) circle (3.5pt) };
  \filldraw[fill=white,draw=black]  \foreach \x in {3,3.5,...,5}{ (\x,2.0) circle (3.5pt) };
  \filldraw[fill=white,draw=black]  \foreach \x in {3.5,4,...,5}{ (\x,2.5) circle (3.5pt) };
  \filldraw[fill=white,draw=black]  \foreach \x in {4.5,5}      { (\x,3.0) circle (3.5pt) };
  \draw[fill=ghostColour,xshift=-.15cm,yshift=0cm]  \foreach \x in {.5,1.,1.5}{ \mytrix };  
  \draw[fill=ghostColour,xshift=-.15cm,yshift=5cm]  \foreach \x in {.5,1.,...,5}{ \mytrix };  
  \draw[fill=ghostColour,xshift=0cm,yshift=-.15cm]  \foreach \y in {0,.5,...,5}{ \mytriy };
  \draw[fill=ghostColour,xshift=5cm,yshift=-.15cm]  \foreach \y in {3.5,4,4.5}{ \mytriy };
   \draw (1.25,3.5) node[thick,draw=blue,fill=white,inner sep=2pt] {\large$G_1$};
\end{scope}
\begin{scope}[xshift=11.5cm,yshift=2.25cm,scale=.75]
\draw[-,thick,green,yshift=.0 cm] 
   \foreach \x in {0,.454545,...,5}{ (\x,0) -- (\x,5) }
   \foreach \y in {0,.833333,...,5}{ (0,\y) -- (5,\y) };
 \draw[very thick,red,xshift=.454545cm,yshift=.833333cm] (0.,4) -- (.0,0) -- (4.0909,0.) -- (4.0909,4);
 \filldraw[blue]  \foreach \x in {.454545,.909090,...,4.545454}{ (\x,5) circle (3.5pt) };
 \draw[fill=ghostColour,xshift=-.15cm]  \foreach \x in {.454545,.909090,...,4.545454}{ \mytrix };
 \draw[fill=ghostColour,yshift=-.15cm]  \foreach \y in {0,.833333,...,5}{ \mytriy };
 \draw[fill=ghostColour,xshift=5cm,yshift=-.15cm]  \foreach \y in {0,.833333,...,5}{ \mytriy };
\end{scope}
\begin{scope}[xshift=7cm,yshift=.7cm]
  \fill[black!10!white,xshift=-.1cm,yshift=-.25cm] (0,0) -- (4,0) -- (4.,1.3) -- (0,1.3) -- (0,0);
  \filldraw[green,xshift=.0cm,yshift=.8cm] (.25,.0)  circle (3pt);
  \filldraw[blue,xshift=.3cm,yshift=.8cm] (.25,.0)  circle (3pt);
  \draw[xshift=.0cm,yshift=.8cm] (.5,0) node[anchor=west,xshift=6] {\small interpolation};
  \draw[fill=ghostColour,xshift=.0cm,yshift=.4cm] (.35,0) \foreach \x in {.1}{ \mytrix } node[anchor=west,xshift=12,yshift=3] {\small ghost};
  \draw[fill=white,draw=black,xshift=.0cm,yshift=.0cm] (.25,0) circle (3.5pt) node[anchor=west,xshift=6] {\small unused};
\end{scope}
\begin{scope}[xshift=11.5cm,yshift=2.25cm,scale=.75]
   \draw (1.6,3.27) node[thick,draw=green,fill=white,inner sep=2pt] {\large$G_2$};
\end{scope}
\end{tikzpicture}
\end{center}
\caption{Left: composite of a background grid ($G_1$, blue) 
  and a boundary-fitted grid ($G_2$, green) in physical space for the domain defined by
  the interior of the red boundary. The grid points on $G_1$ with green dots interpolate from
  $G_2$ and the grid points on $G_2$ with blue dots interpolate from $G_1.$
  Middle: Plot of $G_1$ showing interpolation points, ghost points (grid points which exist 
  outside the physical boundary), and unused points (grid points which do not affect the computation).
  Right: The green boundary fitted grid, $G_2,$ is mapped to a unit square. The plot shows
  interpolation points and ghost points.
\label{fig:overlappingGrid2dCartoon}}
\end{figure}

The new ME schemes have been implemented for complex geometry using overset grids, which are also known as composite overlapping grids
or Chimera grids. 
As shown in Figure~\ref{fig:overlappingGrid2dCartoon}, an overset grid, denoted as~$\Gc$, 
consists of a set of component grids $\{G_g\}$, $g=1,\ldots,{\mathcal N}$, that cover the entire domain~$\Omega$.
Solutions on the component grids are matched by interpolation~\cite{CGNS}.
Overset grids enable the use of efficient finite difference schemes on structured grids, while simultaneously treating complex geometry with high-order accuracy up to and including boundaries.
Each component grid, $G_g$, is a logically rectangular, curvilinear grid
defined by a smooth mapping from a unit cube in $\nd$ dimensions (called the parameter space
with coordinates~$\rv$) to physical
space~$\xv$,
\begin{equation}
  \xv = \Gv_g(\rv),\qquad \rv\in[0,1]^{\nd},\qquad \xv\in\Real^{\nd}.
\label{eq:gridMapping}
\end{equation}
All grid points in $\Gc$ are classified as discretization, interpolation or unused points~\cite{CGNS}.
The overlapping grid generator {\bf Ogen}~\cite{ogen} from the {\it Overture} framework is used to
construct the overlapping grid information.

\subsection{Discrete approximations on curvilinear grids}

Approximations to derivatives on a curvilinear grid can be formed using the \textit{mapping method}.
Given a mapping $\xv = \Gv_g(\rv)$ and its metric derivatives, $\p r_\ell/\p x_m$, $\ell,m=1,\ldots,\nd$, 
the derivatives of a function $u(\xv)=U(\rv)$ are first written in parameter space using the chain rule, for example,
\ba
   \f{\p u}{\p x_m } = \sum_{\ell=1}^{\nd} \f{\p r_\ell}{\p x_m} \f{\p U}{\p r_\ell } .
\ea
Derivatives of $U$ with respect to $r_\ell$ are then approximated with standard finite differences.
Let $\rv_\iv$ denote grid points on the unit cube, where $i_k=0,1,\ldots,N_k$.
Let $\dr_k=1/N_k$ denote the grid spacing on the unit cube with $\rv_{\iv} = (i_1 \dr_1, i_2\dr_2, i_3\dr_3)$. 
Let $U_\iv \approx U(\rv_\iv)$ and define the difference operators,
\ba
   D_{+ r_\ell} U_{\iv} \eqdef \f{U_{\iv+\ev_\ell} - U_\iv}{\dr_\ell},        \qquad 
   D_{- r_\ell} U_{\iv} \eqdef \f{U_{\iv}       - U_{\iv-\ev_\ell}}{\dr_\ell}, \quad 
   D_{0 r_\ell} U_{\iv} \eqdef \f{U_{\iv+\ev_\ell} - U_{\iv-\ev_\ell}}{2\dr_\ell},
   \label{eq:dividedDiffr}
\ea
where $\ev_\ell$ is the unit vector in direction $\ell$. 
Second-order accurate approximations to the first derivatives, for example, are
\ba
   D_{x_m,h} U_\iv \eqdef \sum_{\ell=1}^{\nd} \left.\f{\p r_\ell}{\p x_m}\right\vert_{\iv} D_{0,r_\ell} U_\iv,
\ea
where we assume the metric terms $\p r_\ell/\p x_m$ are known at grid points from the mapping. 
We do not, however, assume the second derivatives of the mapping are known (to avoid the extra storage) and these 
are computed using finite differences of the metrics.
Using the chain rule, the second derivatives are
\ba
  \f{\p^2 u}{\p x_m \p x_n } = \sum_{k=1}^{\nd} \sum_{l=1}^{\nd} \f{\p r_k}{\p x_m} \f{\p r_l}{\p x_n} \f{\p^2  U}{\p r_k\p r_l }
     +  \sum_{k=1}^{\nd} \left\{ \sum_{l=1}^{\nd} \f{\p r_l}{\p x_n} \f{\p}{\p r_l}  \f{\p r_k}{\p x_m} \right\}  \f{\p U}{\p r_k } .
\ea
The second derivatives are approximated to second-order accuracy using approximations such as 
\bat
    \left.\f{\p^2  U}{\p r_k\p r_l }\right\vert_{\rv_\iv} & \approx  D_{+r_k} D_{-r_l} U_\iv , \qquad&& \text{for $k =  l$}, \\
    \left.\f{\p^2  U}{\p r_k\p r_l }\right\vert_{\rv_\iv} & \approx  D_{0r_k} D_{0r_l} U_\iv , \qquad&& \text{for $k\ne l$}, \\
    \left.\f{\p}{\p r_l}\left( \f{\p r_k}{\p x_m}\right)\right\vert_{\rv_\iv} & \approx D_{0r_l} \left( \left.\f{\p r_k}{\p x_m}\right\vert_\iv \right) .
\eat
Fourth and higher-order accurate approximations are straightforward to form using similar techniques.

\subsection{Boundary conditions and upwind dissipation}

Careful attention to the discrete boundary conditions is important for accuracy and stability, especially for the wave equation which has no natural dissipation.
We use compatibility boundary conditions (CBCs) which are generally more accurate and stable than
one-sided approximations~\cite{lcbc2022}.
A simple CBC uses the governing equation on the boundary. More generally, CBCs for the case of the wave equation are formed by taking even time-derivatives
of the boundary condition and then using the governing equation to replace $\p_t^2$ by $\Lc$. 
For flat boundaries with homogeneous Dirichlet or Neumann boundary conditions CBCs lead to odd or even reflection
conditions, respectively.
For more details on CBCs see~\cite{max2006b,lcbc2022} for example.

The upwind dissipation operator $\Qv_p$ was introduced in Section~\ref{sec:upwindDissipation} and defined in~\eqref{eq:Qdef} for
the case of a Cartesian grid.  More generally for a curvilinear grid the upwind dissipation operator is taken as 
\bse
\label{eq:upwindDissSimpleCurvilinear}
\ba
 &  \Qv_p =  \sum_{\ell=1}^{\nd} \f{c\, \| \grad_{\xv} r_\ell \| }{\dr_\ell} 
       \big( -\Delta_{+r_\ell}\Delta_{-r_\ell} \big)^{p/2+1} ,   \\
\intertext{where}
  & \| \grad_{\xv} r_\ell \|^2  = \sum_{m=1}^{\nd} \left[ \f{\p r_\ell }{\p x_m} \right]^2 .        
\ea
\ese
Here, $\Delta_{\pm r_\ell}$ are undivided difference operators in the $\ell$~coordinate direction of the parameter space~$\rv$ corresponding to the divided difference operators defined in~\eqref{eq:dividedDiffr}.

\subsection{Implicit first time-step} \label{sec:implicitFirstStep}

Implicit three-level ME schemes require two time levels to initiate the time stepping.
The solution at $t=0$ is found directly from the initial condition~\eqref{eq:IC1}.
The solution at the first time-step $t=\dt$ could be found from a Taylor series in time 
using both initial conditions~\eqref{eq:IC1}, and \eqref{eq:IC2}, together with the governing equation~\eqref{eq:secondOrderWaveEquation}.
It would be convienent to use an explicit version of this Taylor series approximation
to obtain the solution at the first time-step;
formally this would not change the stability of the scheme.
In practice, however, very large errors can be introduced in the first explicit time-step at $t=\dt$
when the CFL number is large, often rendering the full computation useless.
Thus the first-time step should be taken implicitly when the CFL number is large.
Further, it would be convenient if this implicit first time-step utilized the
same implicit matrix as subsequent time steps of the full three-level scheme.
In this section we show how this can be accomplished.

\subsubsection{Implicit first time-step: second-order accuracy}

Consider the second-order accurate implicit IME2 scheme, 
re-written here for clarity,  
\ba
 \f{U_\iv^{n+1} - 2 U_\iv^n + U_\iv^{n-1}}{\dt^2} = 
       \Lc_{2,h}\Big[ \alpha_2 U_\iv^{n+1} + \beta_2 U_\iv^n + \alpha_2 U_\iv^{n-1}  \Big] , \label{eq:implicitOrder2}
\ea 
which has the implicit operator
\ba
  A_2 \eqdef  I - \alpha_2 \dt^2 \Lc_{2,h} .  \label{eq:a2op}
\ea
Given the initial conditions,
\bse
\ba
    u(\xv,0)      & = u_0(\xv), \\
    \p_t u(\xv,0) & = u_1(\xv),  \label{eq:initialCondition2}
\ea
\ese
approximate~\eqref{eq:initialCondition2} to second-order accuracy using
\ba
   \Dzt U_\iv^n = \f{U_\iv^{n+1} - U_\iv^{n-1}}{2\dt} = u_{1,\iv} ,
\ea
with $n=0$. Solving for $ U_\iv^{n-1}$ gives 
\ba
   U_\iv^{n-1} = U_\iv^{n+1} - (2\dt)u_{1,\iv}  . \label{eq:firstStepD0t}
\ea
Substituting~\eqref{eq:firstStepD0t} into~\eqref{eq:implicitOrder2} and dividing by $2$
gives the following implicit scheme for the first step ($n=0$),
\ba
 \f{ U_\iv^{n+1} - U_\iv^n -  \dt u_{1,\iv} }{\dt^2} = 
      \Lc_{2,h}\Big[ \alpha_2 U_\iv^{n+1} + \half \beta_2 U_\iv^n - \alpha_2 \dt\, u_{1,\iv}  \Big]  .
\ea
This gives the following update for the first time-step ($n=0$)
\ba
  A_2 U_\iv^{n+1}  = U_\iv^n +  \dt\, u_{1,\iv} + \dt^2 \Lc_{2,h}\Big[  \half \beta_2 U_\iv^n - \alpha_2 \dt\, u_{1,\iv}  \Big] ,
  \label{eq:impFirstStepOrder2}
\ea
that uses the same implicit operator $A_2$ as the later time steps.

\subsubsection{Implicit first time-step: fourth-order accuracy}

\mni
Consider the fourth-order accurate implicit IME4 scheme, re-written here for clarity,
\ba
 \f{U_\iv^{n+1} - 2 U_\iv^n + U_\iv^{n-1}}{\dt^2} = 
      \Lc_{4,h}\Big[ \alpha_2 U_\iv^{n+1} + \beta_2 U_\iv^n + \alpha_2 U_\iv^{n-1}  \Big]  
    -\dt^2 \Lc_{2,h}^2 \Big[ \alpha_4 U_\iv^{n+1} + \beta_4 U_\iv^n + \alpha_4 U_\iv^{n-1}  \Big] .
      \label{eq:implicitOrder4}
\ea 
Scheme~\eqref{eq:implicitOrder4} has the implicit operator
\ba
  A_4 \eqdef  I - \alpha_2 \dt^2 \Lc_{4,h} + \alpha_4 \dt^4 \Lc_{2,h}^2 . \label{eq:a4op}
\ea
To approximate the second initial condition~\eqref{eq:initialCondition2} to fourth-order accuracy, 
we use
\ba
  & \p_t u = \Dzt u  - \f{\dt^2}{6} \p_t^3 u + O(\dt^4).  \label{eq:utOrder4}
\ea
To treat the $\p_t^3 u$ term in~\eqref{eq:utOrder4}, 
take the first time derivative of the governing equation and write it in the form
\ba
  \p_t^3 u = \p_t^2 (\p_t u) = \Lc( \p_t u) . 
\ea
Using this expression for $\p_t^3 u$ in~\eqref{eq:utOrder4} gives the approximation
\ba
   \Dzt U_\iv^n - \f{\dt^2}{6}\, \Lc_{2,h} u_{1,\iv}   = u_{1,\iv} , \qquad (n=0). \label{eq:utOrder4b}
\ea
Solving~\eqref{eq:utOrder4b} for $ U_\iv^{n-1}$ gives
\bse
\label{eq:unm1Order4}
\ba
   U_\iv^{n-1} &= U_\iv^{n+1} - (2\dt)u_{1,\iv} - \f{\dt^3}{3} \, \Lc_{2,h} u_{1,\iv} ,  \label{eq:unm1Order4a} \\
               &= U_\iv^{n+1} - W_\iv,  \label{eq:unm1Order4b}
\ea
where
\ba
  W_\iv &\eqdef  (2\dt)u_{1,\iv} + \f{\dt^3}{3} \, \Lc_{2,h} u_{1,\iv}  .
\ea
\ese
Substituting~\eqref{eq:unm1Order4b}  into~\eqref{eq:implicitOrder4}  and dividing by $2$
gives the following fourth-order accurate implicit scheme for the first step ($n=0$),
\ba
\f{ U_\iv^{n+1} - U_\iv^n - \half W_\iv}{\dt^2} & = 
      \Lc_{4,h}\Big[  \alpha_2 U_\iv^{n+1} + \half \beta_2 U_\iv^n - \half\alpha_2 W_\iv  \Big] \nonumber \\
    &  \qquad - \dt^2  \Lc_{2,h}^2\Big[  \alpha_4 U_\iv^{n+1} + \half\beta_4 U_\iv^n - \half\alpha_4 W_\iv  \Big] .
    \label{eq:impFirstStepOrder4A}
\ea
Rearranging~\eqref{eq:impFirstStepOrder4A} gives 
\ba
A_4 U_\iv^{n+1} &= U_\iv^n + \half W_\iv + \dt^2 \Lc_{4,h}\Big[  \half \beta_2 U_\iv^n - \half\alpha_2 W_\iv  \Big]
                 - \dt^4  \Lc_{2,h}^2\Big[  \half\beta_4 U_\iv^n - \half\alpha_4 W_\iv  \Big]  \qquad (n=0),
\ea     
and substituting for $W_\iv$ leads to
\bse
\label{eq:impFirstStepOrder4}
\ba
A_4 U_\iv^{n+1} &= U_\iv^n + \dt\, u_{1,\iv} + \f{\dt^3}{6} {\red \Lc_{2,h} u_{1,\iv}}  \\
                & \qquad + \dt^2 \Lc_{4,h}\Big(  \half \beta_2 U_\iv^n - \alpha_2(\dt\, u_{1,\iv}  {\blue + \f{\dt^3}{6} \Lc_{2,h} u_{1,\iv}})  \Big) \\
                & \qquad\qquad- \dt^4  \Lc_{2,h}^2\Big(  \half\beta_4 U_\iv^n - \alpha_4 (\dt\, u_{1,\iv} {\blue + \f{\dt^3}{6} \Lc_{2,h} u_{1,\iv}})  \Big) .
\ea
\ese
The terms in {\blue blue} are dropped, based on accuracy, and this keeps the stencil compact.
Note that scheme~\eqref{eq:impFirstStepOrder4} has the same implicit operator $A_4$ as scheme~\eqref{eq:implicitOrder4}.
If the red term ${\red \Lc_{2,h} u_{1,\iv}}$ in~\eqref{eq:impFirstStepOrder4} is replaced
by $\Lc_{4,h} u_{1,\iv}$, then the form of~\eqref{eq:impFirstStepOrder4} is similar to the usual interior update and the
same code can be used for both the first step and later steps with the appropriate choice of coefficients.

\subsection{Solution of the implicit time-stepping equations} \label{sec:implicitSolvers}

For the overset grid results presented in this article the implicit equations are solved
either with a direct sparse solver (for problems that are not too large) or with a
Krylov space method (bi-CG-stab with an ILU preconditioner). 
The approach we currently use is not efficient since an implicit matrix is formed for all grids points, 
on both explicit and
implicit grids. The implicit matrix entries corresponding to the equations at points treated explicitly are simply the identity.
This approach was done so that existing software could be used. A more efficient approach would be to only form
a system for the implicit points. Moreover, it is often the case that the implicit points on different
component grids are not coupled and in this case multiple smaller implicit systems could be formed.

\section{Matrix stability analysis on one-dimensional overset grids} \label{sec:matrixStability}

In this section, matrix stability analysis is used to study the behavior of the new schemes on
a collection of one-dimensional overset grids.
A large number of overset grids with different grid spacings are considered to determine the stability behavior for a wide range of grid configurations.
A scaled upwind dissipation coefficient $\nu_\gamma = \gamma \nu_p$ is used with $\gamma\in[0,1]$ to show how the stability of the scheme depends on the amount of dissipation ranging from no dissipation $\gamma=0$ to {\it full} dissipation $\gamma=1$.  In particular, the results for $\gamma=0$ show the necessity of upwind
dissipation. For some cases, the number of upwind corrections, $\Nuc$, must also be chosen greater than one for stability.

The time-stepping update on an overset grid is written in the form of a single vector update for the 
\textit{active unknowns} (i.e.~unknowns corresponding to the interior equations)
excluding the \textit{constraint unknowns} (i.e.~unknowns associated 
with the boundary conditions and interpolation equations). For homogeneous boundary conditions, this update is written
in terms of the matrices $B_1$ and $B_2$ and the vector $\Vv^n$ of active unknowns at time $t^n$, 
\ba
   & \Vv^{n+1} = B_1  \Vv^n  + B_2 \Vv^{n-1} . \label{eq:matrixStabilityEquationsB}
\ea
For a given overset grid, the associated eigenvalues and eigenvectors can be determined and this
shows whether discrete solutions grow in time or not.

\newcommand{\nGhost}{n_{\rm ghost}}
\newcommand{\tola}{{\rm tol}_a}
\subsection{Matrix stability formulation} \label{sec:matrixStabilityFormulation}

Some details on the construction of the matrix stability equation~\eqref{eq:matrixStabilityEquationsB} is now described.
We assume the problem domain is $\Omega=[-1,1]$ and let it be discretized with an overset grid  as shown in Figure~\ref{fig:finite1doverlapgrid}.
The left domain is $\Omega_\Lss=[-1,b_\Lss]$, where $b_\Lss$ may vary, and the right domain is fixed at $\Omega_\Rss=[0.5,1]$.
Let $U_{\Lss,j}^n$ and  $U_{\Rss,j}^n$ denote grid functions on the left and right grids, respectively, at time~$t^n$.
The active points on the left grid are $j=1,2,\ldots,N_{\Lss}$, while those on the
right grid are $j=0,1,2,\ldots,N_{\Rss}-1$. Dirichlet boundary conditions
are given at $j=0$ (left grid) and at $j=N_{\Rss}$ (right grid). 
The interpolation points on the left grid are at $j=N_{\Lss}+1, \ldots,N_{\Lss}+\nGhost$, where
$\nGhost=p/2+1$ is the number of ghost points. 
The interpolation points on the right grid are at $j=-\nGhost,\ldots,-1$.  The grid spacings for the left and right grids are uniform with $h_L=(b_L+1)/N_L$ and $h_R=0.5/N_R$.  In our study of stability for a collection of overset grids, $N_R$ is held fixed and thus $h_R$ is also fixed.  We then select a ratio of the grid spacings $\delta=h_L/h_R$ which implies $b_L=-1+\delta h_RN_L$.  The value for $N_L$, and the corresponding value for $b_L$, is chosen to minimize the grid overlap while maintaining the assumption of an explicit interpolation as discussed in more detail below.  The sampling of overset grids is performed for grids for a range of grid spacing ratios
$\delta\in[\delta_{{\rm min}},\delta_{{\rm max}}]$ as noted below.

{
\newcommand{\tikzcircle}[2][red,fill=red]{\tikz[baseline=-0.5ex]\draw[#1,radius=#2] (0,0) circle ;}%
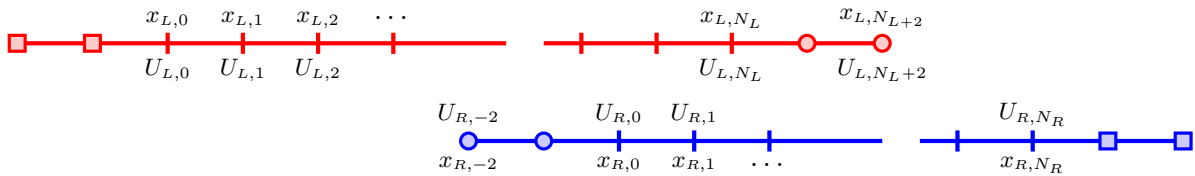
\begin{figure}[hbt!]
  \centering
  \begin{tikzpicture}
    \useasboundingbox (0,.7) rectangle (12,2.8);

     \begin{scope}[yshift=-0.7cm,xshift=.25cm]
     \draw[ultra thick,red] (-2,3) -- (4.5,3);
     \draw[ultra thick,red] (5,3) -- (9.5,3);
     \foreach \x in {0,1,2,3,5.5,6.5,7.5}
       \draw[ultra thick,red] (\x,3cm-4pt) -- (\x,3cm +4pt); 
 
     \draw (8.5,3) node {\tikzcircle[red,very thick, fill=red!20]{3pt}};
     \draw (9.5,3) node {\tikzcircle[red,very thick, fill=red!20]{3pt}};

     \draw [red, very thick,fill=red!20] plot [only marks, mark size=3pt, mark=square*] coordinates {(-2,3) (-1,3)};
 
     \draw (0.0,3) node[anchor=south,yshift=2pt] {\labelFont$x_{\Lss,0}$} node[anchor=north,yshift=-2pt]{\labelFont$U_{\Lss,0}$};
     \draw (1.0,3) node[anchor=south,yshift=2pt] {\labelFont$x_{\Lss,1}$} node[anchor=north,yshift=-2pt]{\labelFont$U_{\Lss,1 }$};
     \draw (2.0,3) node[anchor=south,yshift=2pt] {\labelFont$x_{\Lss,2}$} node[anchor=north,yshift=-2pt]{\labelFont$U_{\Lss,2}$};
     \draw (3,3.35) node[]{$\dots$};
     \draw (7.5,3) node[anchor=south,yshift=2pt] {\labelFont$x_{\Lss,N_{\Lss}}$} node[anchor=north,yshift=-2pt]{\labelFont$U_{\Lss,N_{\Lss}}$};
     \draw (9.5,3) node[anchor=south,yshift=2pt] {\labelFont$x_{\Lss,N_{\Lss+2}}$} node[anchor=north,yshift=-2pt]{\labelFont$U_{\Lss,N_{\Lss}+2}$};    
    \end{scope}

    \begin{scope}[xshift=.75cm]
      \draw[ultra thick,blue] (3.5,1) -- (9,1);
      \draw[ultra thick,blue] (9.5,1) -- (13,1);
       \foreach \x in {5.5,6.5,7.5,10,11}
        \draw[ultra thick,blue] (\x,1cm-4pt) -- (\x,1cm +4pt); 
  
      \draw (3.5,1) node {\tikzcircle[blue,very thick, fill=blue!20]{3pt}};
      \draw (4.5,1) node {\tikzcircle[blue,very thick, fill=blue!20]{3pt}};

     \draw [blue, very thick,fill=blue!20] plot [only marks, mark size=3pt, mark=square*] coordinates {(12,1) (13,1)};

      \draw (3.5,1) node[anchor=north,yshift=-2pt] {\labelFont$x_{\Rss,-2 }$} node[anchor=south,yshift=2pt]{\labelFont$U_{\Rss,-2}$};
      \draw (5.5,1) node[anchor=north,yshift=-2pt] {\labelFont$x_{\Rss, 0 }$} node[anchor=south,yshift=2pt]{\labelFont$U_{\Rss, 0}$};
      \draw (6.5,1) node[anchor=north,yshift=-2pt] {\labelFont$x_{\Rss, 1 }$} node[anchor=south,yshift=2pt]{\labelFont$U_{\Rss, 1}$};
      \draw (11.,1) node[anchor=north,yshift=-2pt] {\labelFont$x_{\Rss,N_{\Rss}}$} node[anchor=south,yshift=2pt]{\labelFont$U_{\Rss,N_{\Rss}}$};

      \draw (7.5,1) node[yshift=-9pt]{$\dots$};
    \end{scope}
%
  \end{tikzpicture}
  \caption{One-dimensional overset grid used for the matrix stability analyses. Interpolation points are marked as circles, ghost points are marked as squares.}
  \label{fig:finite1doverlapgrid}
\end{figure} 
}

In Stages 1 and 2 of the SPIE scheme in~\eqref{eq:SPIEscheme}, the solutions for the left and right grids are advanced one time step (without dissipation) and the boundary/interpolation conditions are applied.  This step can be written in matrix form as
\ba
   Q_0 \Uv^{(0)} = Q_1 \Uv^n + Q_2 \Uv^{n-1}, \label{eq:stage1}
\ea 
where $\Uv^n$ is a vector holding all of the unknowns on the two grids (including ghost points and interpolation points).
The equations in~\eqref{eq:stage1} include the interior equations, boundary conditions and interpolation equations.
The matrix $Q_0$ is equal to the identity matrix at active points using an explicit scheme, while the matrix has values corresponding to the implicit operators $A_2$ in~\eqref{eq:a2op} and $A_4$ in~\eqref{eq:a4op} at active points corresponding to the second and fourth-order accurate implicit schemes, respectively.  The matrix $Q_0$ also includes the boundary conditions and interpolation equations; 
the corresponding rows in $Q_1$ and $Q_2$ are zero.
For example, the boundary conditions on the left grid are
\ba
   & U_{\Lss,0}^{(0)} = 0 , \\
   & U_{\Lss,-j}^{(0)} = - U_{\Lss,j}^{(0)}, \qquad j=1,\ldots,\nGhost,
\ea
where the odd symmetry conditions on the ghost points are determined from compatibility conditions.
The boundary conditions on the right grid are similar.
The values at interpolation points are found using Lagrange interpolation for a stencil of $p+1$ points. For example, an interpolation point on the left grid is
found using a formula of the form 
\ba
  U_{\Lss,k}^{(0)} = \sum_{j=1}\sp{p+1} w_{k,j}^{\Lss} \, U_{\Rss,m_{k}+j}^{(0)} . \label{eq:interpEquations}
\ea
where $m_{k}$ denotes the left index of the interpolation stencil, chosen to make the interpolation
as centered as possible, and $w_{k,j}^{\Lss}$ are interpolation weights. The values on the right hand side of~\eqref{eq:interpEquations} are known as \textit{donor} points.
The interpolation is taken to be explicit so that none of the donor points for one grid are interpolation points for the other grid.

Explicit upwind dissipation is incorporated in Stage 3 of the SPIE scheme. Assuming $\Nuc$ applications of the dissipation, the updates of the solution at this stage have the matrix form
\ba
  &  P_0 \Uv^{(k)} = P_1 \Uv^{(k-1)} + P_2 \Uv^{n-1},   \qquad k=1,2,\ldots,\Nuc, \label{eq:stage2}
\ea 
where $P_0$, like $Q_0$, includes the boundary/interpolation equations.
Finally, the solution at the new time is 
\ba
   \Uv^{n+1} = \Uv^{(\Nuc)}.
\ea
Combining the time step in Stages 1/2 and the corrections in Stage 3 leads to a three-level matrix equation of the form
\ba
   & \Uv^{n+1} = A_1  \Uv^n  + A_2 \Uv^{n-1} . \label{eq:matrixStabilityEquationsA}
\ea
where $A_1$ and $A_2$ are coefficient matrices generated from the ones in~\eqref{eq:stage1} and~\eqref{eq:stage2}.
The constraint unknowns in~\eqref{eq:matrixStabilityEquationsA} can be eliminated
by row operations and this leads to the compressed form in~\eqref{eq:matrixStabilityEquationsB}.
The correctness of the matrices $B_1$ and $B_2$ in~\eqref{eq:matrixStabilityEquationsB} is checked
by comparing, at each time-step, the solution computed using the SPIE scheme in~\eqref{eq:SPIEscheme} with the solution
arising from the compressed form~\eqref{eq:matrixStabilityEquationsB}.

To investigate the growth of solutions to the discrete problem~\eqref{eq:matrixStabilityEquationsB}
we look for solutions of the form $\Vv^n = a^n \, \Vv^0$ which leads to a
quadratic eigenvalue problem for $a$ given by
\ba  
   (a^2 I - a B_1 - B_2) \Vv^0 = 0 .
\ea
This quadratic form can also be written as a regular eigenvalue problem of twice the dimension,
\ba
  \begin{bmatrix}
      O   &  I \\
      B_2 & B_1 
  \end{bmatrix}
  \begin{bmatrix}
     \Vv^0 \\
     \Vv^1
  \end{bmatrix}
  = 
   a 
  \begin{bmatrix}
     \Vv^0 \\
     \Vv^1
  \end{bmatrix}   .  \label{eq:matrixStabEigProblem}
\ea
The eigenvalue problem in~\eqref{eq:matrixStabEigProblem} is solved easily with standard software.

For the stability studies, we set $N_R=10$ for the right grid with fixed domain $\Omega_{\Rss}=[0.5,1]$. This right grid represents
the local boundary grid in a general overset grid. The grid spacing on the left grid is determined by $h_L=\delta h_R$, where $\delta$ is the ratio of grid spacings.  This ratio is varied from $\half$ to $2$ to represent typical overset grids
where the grid spacings in the overlap are chosen to be nearly the same.
For each value of $\delta$, the parameters $b_L$ and $N_L$ for the left grid are determined based on the grid overlap as discussed above.
Finally, the parameter $\gamma$ for the scaled upwind dissipation coefficient $\nu_\gamma = \gamma \nu_p$ is varied between $0$ and $1$ to study how much dissipation is needed to stabilize the SPIE scheme for the different cases considered.

{
\newcommand{\figw}{7cm}
\newcommand{\figh}{6cm}
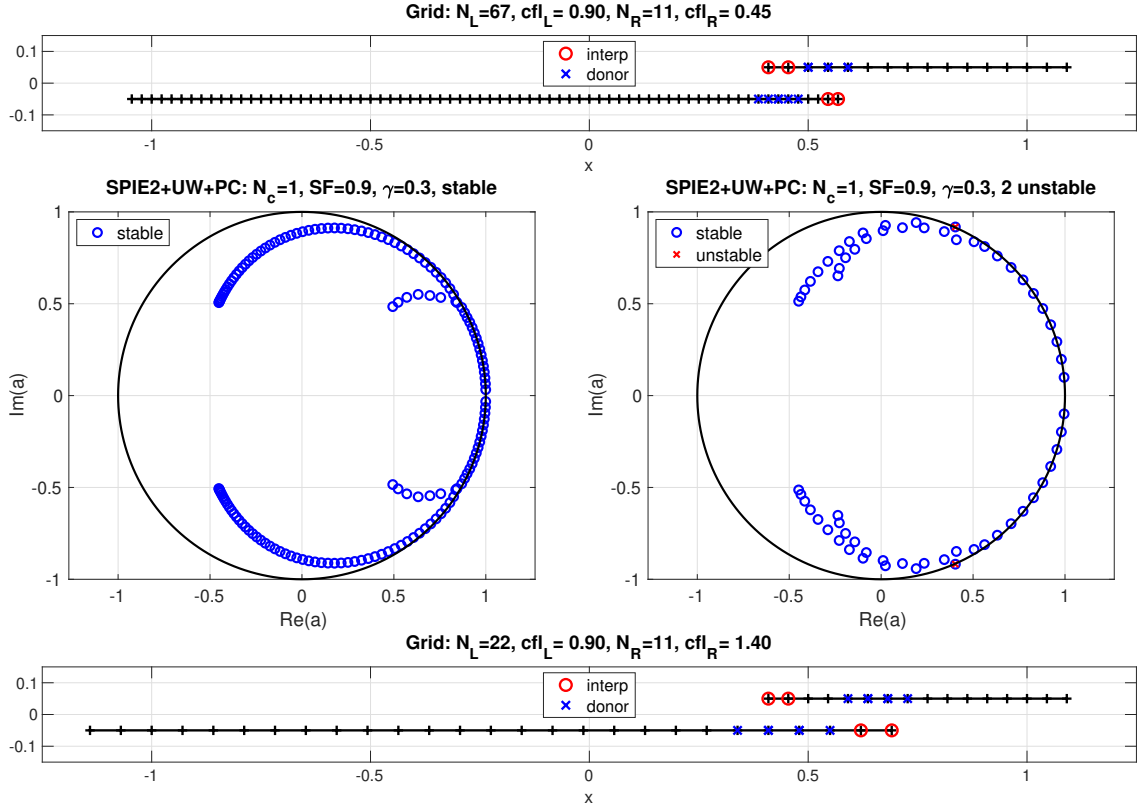
\begin{figure}[htb]
\begin{center}
\begin{tikzpicture}
   \useasboundingbox (0,.5) rectangle (15,1.8*\figh);  
  \figByWidth{  0}{1.4*\figh}{gridSPIE2_UW_PCStableCaseCase1}{15cm}[0.][0.][0.][0.]
  \begin{scope}[yshift=2.3cm]
  \figByWidth{        0}{0*\figh}{eigsComplexPlaneSPIE2_UW_PCStableCaseCase1}{\figw}[0.][0.][0.][0.]
  \figByWidth{1.1*\figw}{0*\figh}{eigsComplexPlaneSPIE2_UW_PCUnstableCase8}{\figw}[0.][0.][0.][0.]
  \end{scope}
  \figByWidth{  0}{0.0*\figh}{gridSPIE2_UW_PCUnstableCase8}{15cm}[0.][0.][0.][0.]

\end{tikzpicture}
\end{center}
\caption{Example stable and unstable cases for the SPIE2-UW-PC scheme with $\gamma=0.3$. 
   Middle left: amplification factors $a$ for the stable case corresponding to the grid on the top, grid-ratio $\delta=.5$.
   Middle right: amplification factors $a$ for the unstable case corresponding to the grid on the bottom, grid-ratio $\delta=1.55$.
   }
\label{fig:stableAndUnstableModesSPIE2}
\end{figure}
}

Figure~\ref{fig:stableAndUnstableModesSPIE2} shows some sample results for the SPIE2-UW-PC scheme.
The left grid uses the EME2 explicit scheme with $\CFL=0.9$, while the right grid uses the IME2 implicit scheme.
Explicit upwind dissipation is added in a corrector step with $\nu_p$ given by~\eqref{eq:nuprecommended} and $\Nuc=1$ corrections.
The safety factor for $\nu_p$ is chosen as $s_f=.9$. 
Two grid cases are shown for $\gamma=0.3$; 
one with grid-ratio $\delta=.5$ and one with $\delta=1.55$.  
The grid plotted on the top of the figure is stable as shown in the middle left plot; all eigenvalues of $a$ satisfy $|a|\le 1 + \tola$,
where $\tola=10^{-8}$.
The grid on the bottom has two unstable modes, as illustrated on the middle right plot.
The conclusion for this representative case is that there is insufficient dissipation for the SPIE2-UW-PC scheme with $\gamma=0.3$ and $\Nuc=1$.

\subsection{Matrix stability numerical results} \label{sec:matrixStabilityNumericalResults}

Results are now presented using different combinations of explicit and implicit schemes;
different orders of accuracy ($p=2$ and $p=4$), and different numbers of upwind corrections.
The following cases are considered 
\begin{enumerate}
   \item EMEp: explicit p$^{\rm th}$-order accurate ME schemes are used on the left and right grids.
   \item IMEp: implicit p$^{\rm th}$-order accurate ME schemes are used on the left and right grids.
   \item SPIEp : An EMEp scheme is used on the left grid and an IMEp scheme is used on the right grid.
\end{enumerate}
Explicit upwind dissipation is used in all cases. 
Unless otherwise specified, the IME schemes use the parameters,  $\alpha_2=1/4$ and $\alpha_4=1/12$. 
For each value of upwind scaling factor $\gamma$, the grid-ratio $\delta$ is varied form $.25$ to $2$ using $N_\delta=101$ different values (i.e. $N_\delta$ different
overset grids).

{
\newcommand{\figw}{7cm}
\newcommand{\figh}{6cm}
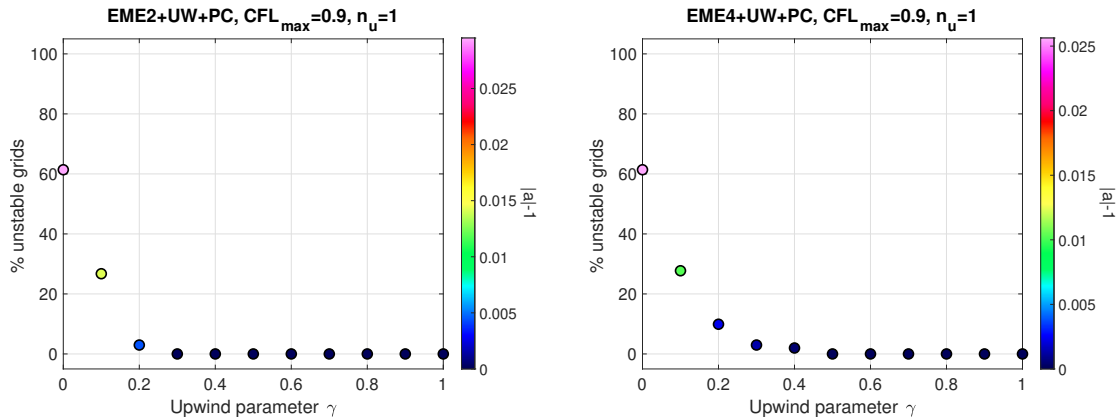
\begin{figure}[htb]
\begin{center}
\begin{tikzpicture}
   \useasboundingbox (0,.5) rectangle (15,.95*\figh);  
  \figByWidth{        0}{0*\figh}{numUnstableModesEME2NC1}{\figw}[0.][0.][0.][0.]
  \figByWidth{1.1*\figw}{0*\figh}{numUnstableModesEME4NC1}{\figw}[0.][0.][0.][0.]
\end{tikzpicture}
\end{center}
\caption{
Overset grid EMEp: Fraction of unstable grids versus dissipation parameter $\gamma$ for the EMEp scheme on an overset grid.
  The colour of each dot corresponds to the value of $|a|-1$ for the most unstable mode.
}
\label{fig:unstableModesEME}
\end{figure}
}

Figure~\ref{fig:unstableModesEME} shows results for the EMEp scheme.
The time-step is chosen so that the CFL number is $0.9$ on the side with the smallest grid spacing.
The number of unstable grids is plotted versus $\gamma$, the scaling factor of the dissipation coefficient $\nu_p$.
For $p=2$ and $\gamma=0$ (no dissipation), roughly $60\%$ of the grids tested are unstable.
This number drops to about $25\%$ when $\gamma=0.1$, and there are no unstable grids for $\gamma\ge 0.3$.
For $p=4$ a value of about $\gamma=0.5$ is sufficient to stabilize all the grids tested.

{
\newcommand{\figw}{7cm}
\newcommand{\figh}{6cm}
\begin{figure}[htb]
\begin{center}
\begin{tikzpicture}
   \useasboundingbox (0,.5) rectangle (15,.97*\figh);  
  \figByWidth{        0}{0*\figh}{numUnstableModesIME2CFL4NC4}{\figw}[0.][0.][0.][0.]
  \figByWidth{1.1*\figw}{0*\figh}{numUnstableModesIME4CFL5NC5}{\figw}[0.][0.][0.][0.]
\end{tikzpicture}
\end{center}
\caption{
Overset grid IMEp: Fraction of unstable grids versus dissipation parameter $\gamma$ for the IMEp scheme on an overset grid.
  The colour of each dot corresponds to the value of $|a|-1$ for the most unstable mode.
  }
\label{fig:unstableModesIME}
\end{figure}
}

Figure~\ref{fig:unstableModesIME} shows results for IMEp schemes.
Note that these schemes use a single implicit solve over both grids (i.e.~the implicit solves are coupled, not partitioned).
For $p=2$ ($p=4$), the time-step is chosen so that CFL number is $4.0$ ($5.0$) on the side with the smallest grid spacing.
The number of explicit upwind corrections is set to $\Nuc=4$ for $p=2$ and $\Nuc=5$ for $p=5$.
This choice is made since the value of $\nu_p$ for the IME schemes scales with the inverse of the CFL.
The results in Figure~\ref{fig:unstableModesIME} show that the schemes have no unstable grids for $\gamma \ge 0.1$.

{
\newcommand{\figw}{7cm}
\newcommand{\figh}{6cm}
\begin{figure}[htb]
\begin{center}
\begin{tikzpicture}
   \useasboundingbox (0,.5) rectangle (15,1.95*\figh);  
  \figByWidth{        0}{0*\figh}{numUnstableModesSPIE2FWNC2}{\figw}[0.][0.][0.][0.]
  \figByWidth{1.1*\figw}{0*\figh}{numUnstableModesSPIE4FWNC2}{\figw}[0.][0.][0.][0.]
  \figByWidth{        0}{1*\figh}{numUnstableModesSPIE2FWNC1}{\figw}[0.][0.][0.][0.]
  \figByWidth{1.1*\figw}{1*\figh}{numUnstableModesSPIE4FWNC1}{\figw}[0.][0.][0.][0.]  

\end{tikzpicture}
\end{center}
\caption{Overset grid SPIE: 
Fraction of unstable grids versus dissipation parameter $\gamma$ for the SPIE scheme on an overset grid 
with weights $\alpha_2=1/4$, and $\alpha_4=1/12$.
  The colour of each dot corresponds to the value of $|a|-1$ for the most unstable mode.
}
\label{fig:unstableModesSPIE_FW}
\end{figure}
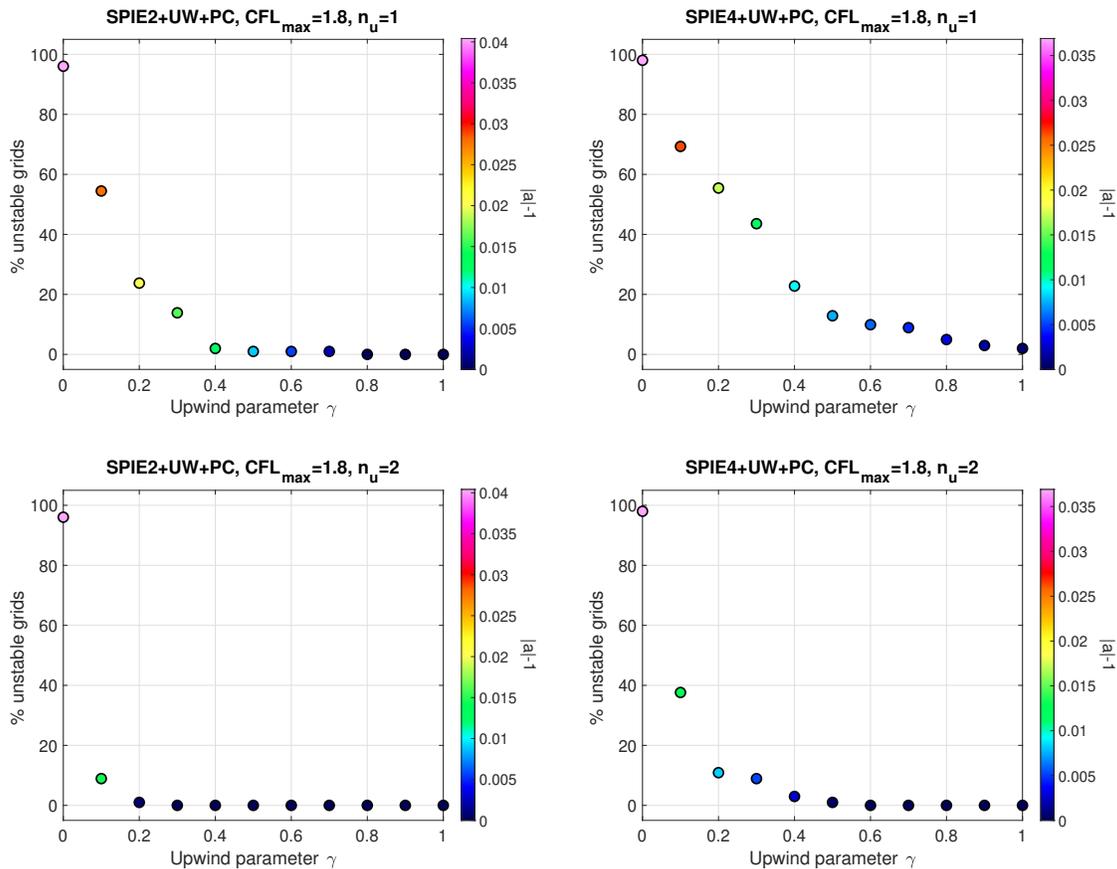
}

Figure~\ref{fig:unstableModesSPIE_FW} shows results for the implicit-explicit SPIE scheme.
The left grid uses an explicit solver and the overall time-step is chosen to match a CFL
number of $0.9$ on this grid. The CFL number on the implicit grid varies between grids and
reaches a maximum of $1.8$.
The left column of plots show results for the second-order accurate SPIE2+UW+PC scheme using $\Nuc=1$ ($s_f=0.9$) and $\Nuc=2$ ($s_f=1.9$)
upwind corrections. With $\Nuc=1$ there are no unstable grids for $\gamma\ge 0.8$.
For $\Nuc=2$, which incorporates more dissipation, there are no unstable grids for $\gamma\ge 0.3$.
The right column of plots show corresponding results for the fourth-order accurate SPIE4+UW+PC scheme using $\Nuc=1$ ($s_f=0.9$) and $\Nuc=2$ ($s_f=1.9$). 
This fourth-order accurate SPIE scheme presents a more difficult case to keep stable. 
With $\Nuc=1$, there are some unstable grids even for $\gamma=1$, while the grids are stable for $\gamma\ge 0.6$ using $\Nuc=2$.

{
\newcommand{\figw}{7cm}
\newcommand{\figh}{6cm}
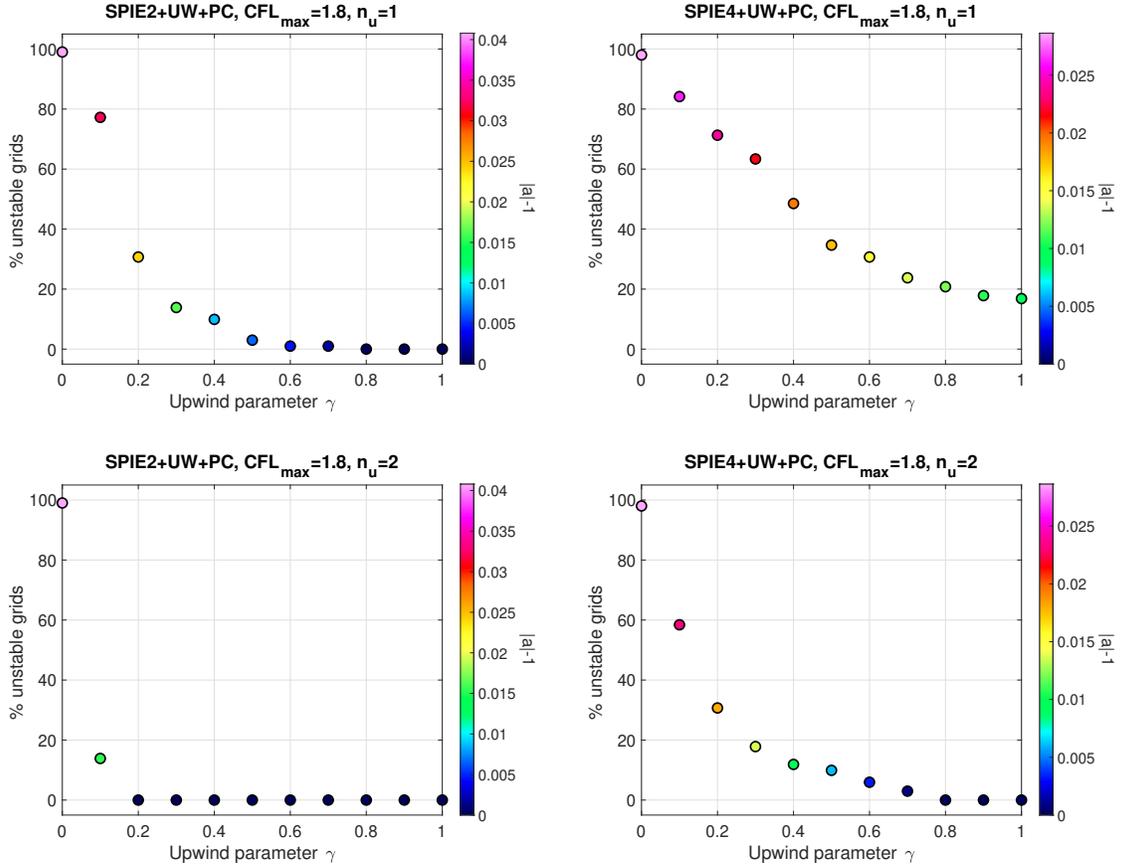
\begin{figure}[htb]
\begin{center}
\begin{tikzpicture}
   \useasboundingbox (0,.5) rectangle (15,1.95*\figh);  
  \figByWidth{        0}{0*\figh}{numUnstableModesSPIE2TrapNC2}{\figw}[0.][0.][0.][0.]
  \figByWidth{1.1*\figw}{0*\figh}{numUnstableModesSPIE4TrapNC2}{\figw}[0.][0.][0.][0.]
  \figByWidth{        0}{1*\figh}{numUnstableModesSPIE2TrapNC1}{\figw}[0.][0.][0.][0.]
  \figByWidth{1.1*\figw}{1*\figh}{numUnstableModesSPIE4TrapNC1}{\figw}[0.][0.][0.][0.]  
\end{tikzpicture}
\end{center}
\caption{Overset grid SPIE + Trap: 
Fraction of unstable grids versus dissipation parameter $\gamma$ for the SPIE scheme on a overset
grid with trapezoidal weights $\alpha_2=1/2$, and $\alpha_4=5/24$.
  The colour of each dot corresponds to the value of $|a|-1$ for the most unstable mode.
}
\label{fig:unstableModesSPIE_TRAP}
\end{figure}
}

A possible reason the SPIE schemes may require more dissipation to remain stable compared to the other cases is that there is a mismatch
in the truncation errors between the left and right grids, the IME schemes generally having larger
errors than the corresponding EME schemes.
Some support for this hypothesis comes from results shown in Figure~\ref{fig:unstableModesSPIE_FW} when using
the \textit{trapezoidal} IME scheme
with $\alpha_2=1/2$ and $\alpha_4=5/24$.
The \textit{Trapezoidal} scheme has a larger truncation error compared to the default scheme.
As seen in Figure~\ref{fig:unstableModesSPIE_TRAP} the \textit{trapezoidal} scheme is more difficult to
stabilize.

\section{Numerical results}  \label{sec:numericalResults}

In this section we present numerical results to demonstrate the
accuracy, stability, and efficiency 
of the proposed new implicit modified equation schemes.  The results are organized into two groups.  The results
in first group are aimed at demonstrating the accuracy and stability of the schemes.  For this group, numerical solutions are
computed for several problems in two and three dimensions where exact solutions are available.  In addition, long-time
simulations are performed for problems in two and three dimensions with random initial conditions as a demonstration of the stability of the schemes.
The results in the second group are used to illustrate the performance of the SPIE schemes for problems with geometric stiffness.

In all examples the wave speed $c$ is taken to be one. For overset grid problems using the SPIE scheme,
 the curvilinear grids are taken to be implicit and the Cartesian grids are taken to be explicit.
Unless otherwise specified the IME schemes use the implicit weighting parameters
\ba
   \alpha_2 = \f{1}{4}, \qquad \alpha_4 = \f{1}{12} .
\ea
For grids $g$ using an implicit method, the coefficient of upwind dissipation  is chosen as 
\ba
    \nu_p = \f{ s_f }{2^{p+1} \, \sqrt{\nd} } \, \f{1}{\lambda_g},   
\ea
where $s_f$ is the safety-factor and 
where $\lambda_g$ is the CFL number for grid $g$, which on a Cartesian grid is given by
\ba
  \lambda_g  \eqdef  c \, \dt \, \sqrt{ \sum_{d=1}^{\nd}\f{1}{h_d^2}  } = \sqrt{ \sum_{d=1}^{\nd}\lambda_{x_d}^2 }. 
\ea
For grids using an explicit method we take
\ba
   \nu_p = \f{ s_f }{2^{p+1} \, \sqrt{\nd} }   .  
\ea
since $\lambda_g\approx1$ for such grids according to the CFL condition.

\subsection{Accuracy and stability of the IME and SPIE schemes}

We begin with numerical results illustrating the accuracy and stability of the second and fourth-order accurate IME and SPIE schemes.

\newcommand{\Gcd}{\Gc_{\rm disk}}
\newcommand{\mTheta}{{m_{\theta}}}
\newcommand{\mr}{{m_r}}
\subsubsection{Eigenmodes on a disk} \label{sec:diskEigenModes}

In this section, eigenmodes of the unit disk in two dimensions are computed.
We look for time-periodic solutions to the wave equation.
In polar coordinates $(r,\theta)$, these solutions have the form
\ba
     u_{\mTheta,\mr}(r,\theta,t) = J_{\mTheta}( k_{\mTheta,\mr} r ) \,  e^{i \mTheta \theta} \, e^{i \omega t}
\ea
where $J_\mTheta$ is the Bessel function of the first kind of (integer) order $\mTheta$,
 $k_{\mTheta,\mr}$, $\mr=1,2,\ldots$ are the positive zeros of $J_\mTheta$ (for the case of Dirichlet boundary conditions)
or $J_\mTheta'$ (for the case of Neumann boundary conditions).  The frequency of vibration for a particular eigenmode is given by 
\ba
   \omega = c \, k_{\mTheta,\mr} .
\ea

Numerical solutions are computed using the fully implicit IME-UW-PC scheme and the mixed explicit/implicit SPIE-UW-PC scheme using the overset grid for the unit disk consisting of a background Cartesian grid (blue) and an annular boundary-fitted grid (green) as shown in Figure~\ref{fig:diskEigenfunctionSolution}.  The grid, denoted by $\Gcd^{(j)}$, has a target grid spacing of $\ds^{(j)}=1/(10 j)$, where the index $j$ determines the size of the grid spacing.  The figure also shows a representative solution at $t=1$  
computed using the IME scheme and the grid $\Gcd^{(4)}$ for the case $(\mTheta,\mr)=(2,2)$ and Dirichlet boundary conditions.  For this Dirichlet case, $k_{\mTheta,\mr}\approx 8.41724414$, while $k_{\mTheta,\mr}\approx 6.70613319$ for Neumann case (not shown).  The rightmost plot in the figure shows the (signed) max-error in the solution.  The error is seen to be smooth with negligible artifacts due to the interpolation at the grid overlap.

{
\newcommand{\drawContour}[7]{%
\begin{scope}[#1]
\draw(0.0,0) node[anchor=south west,xshift=-4pt,yshift=+0pt] {\trimfiga{#2}{\figWidtha}};
  \draw(.5,.5) node[draw,fill=white,anchor=west,xshift=2pt,yshift=1pt,inner sep=2pt] {\scriptsize #3};
\begin{scope}[xshift=0cm,yshift=-2pt]
  \draw (\xcb,\ycb) node[anchor=south west,xshift=0.25cm,yshift=.5cm,rotate=-90] {\trimfigcb{colourBarLines}{\cbWidth}{\cbHeight}};
  \draw (.8,0) node[anchor=north,xshift=+3pt,yshift=+2pt] {\scriptsize $#6$};
  \draw (4.8,0) node[anchor=north,xshift=+0pt,yshift=+2pt] {\scriptsize $#7$};
\end{scope}
\end{scope}
}
\newcommand{\cbWidth}{.2cm}
\newcommand{\cbHeight}{4cm}
\newcommand{\xcb}{.5cm}
\newcommand{\ycb}{-.2cm}
\setlength{\ycbTop}{\ycb+\cbHeight}
\setlength{\ycbMid}{\ycb+\cbHeight*\real{.5}}
\newcommand{\trimfigcb}[3]{\includegraphics[width=#2, height=#3, clip, trim=17cm 2.35cm 1.65cm 2.35cm]{#1}}
\newcommand{\figWidtha}{5cm}
\newcommand{\trimfiga}[2]{\trimw{#1}{#2}{.03}{.12}{.03}{.11}}
\begin{figure}[htb]
\begin{center}
\begin{tikzpicture}
   \useasboundingbox (0,.3) rectangle (16,5);  

  \figByWidth{-.2}{0.15}{sicGridG2}{5.1cm}[0.1][0.1][0.1][0.1]

   \begin{scope}[xshift=4.75cm]
     \drawContour{xshift=0.cm,yshift=0.00cm}{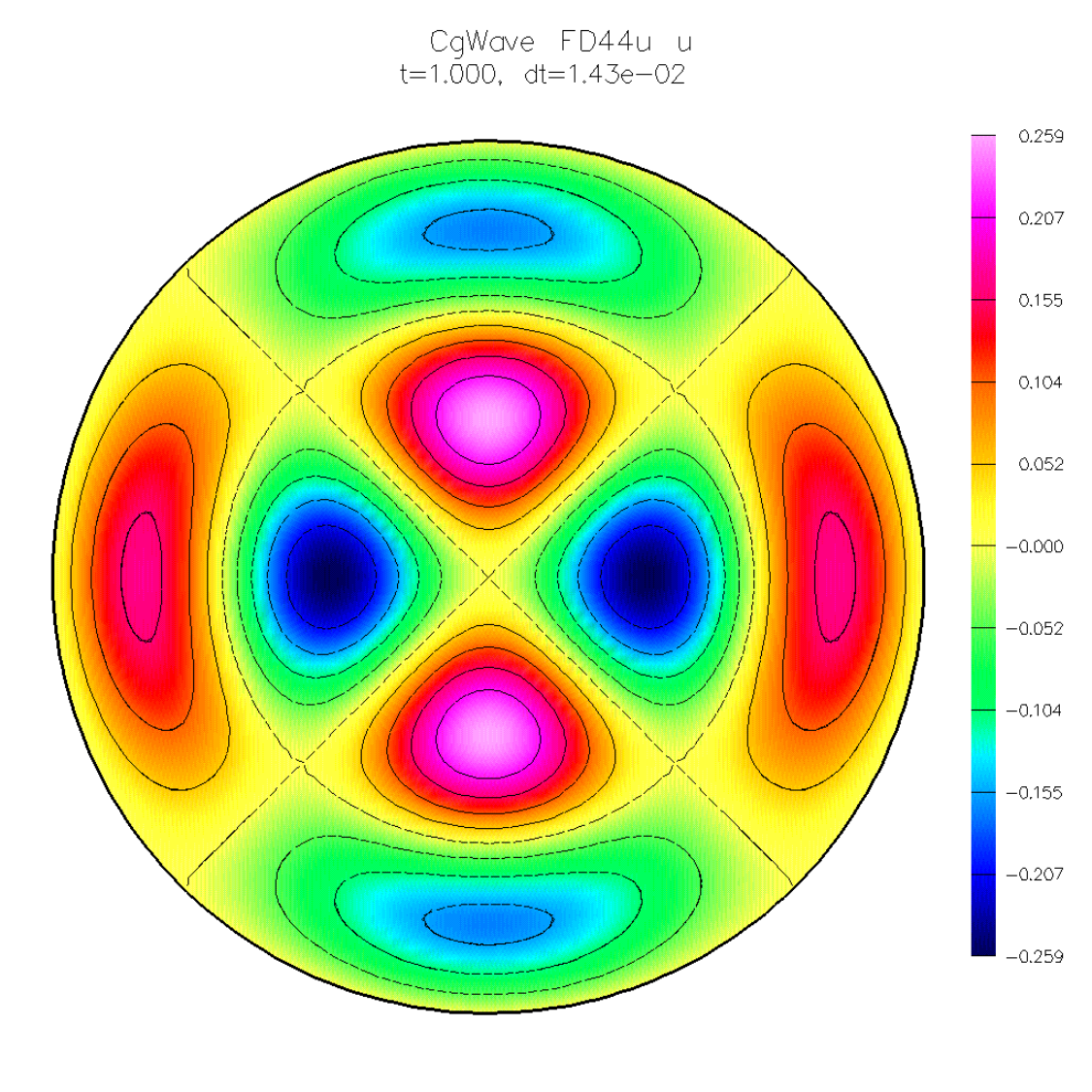}{u, $t=1$}{$v$}{$t=1.0$}{$-.26$}{$.26$}     
     \drawContour{xshift=5.5cm,yshift=0.00cm}{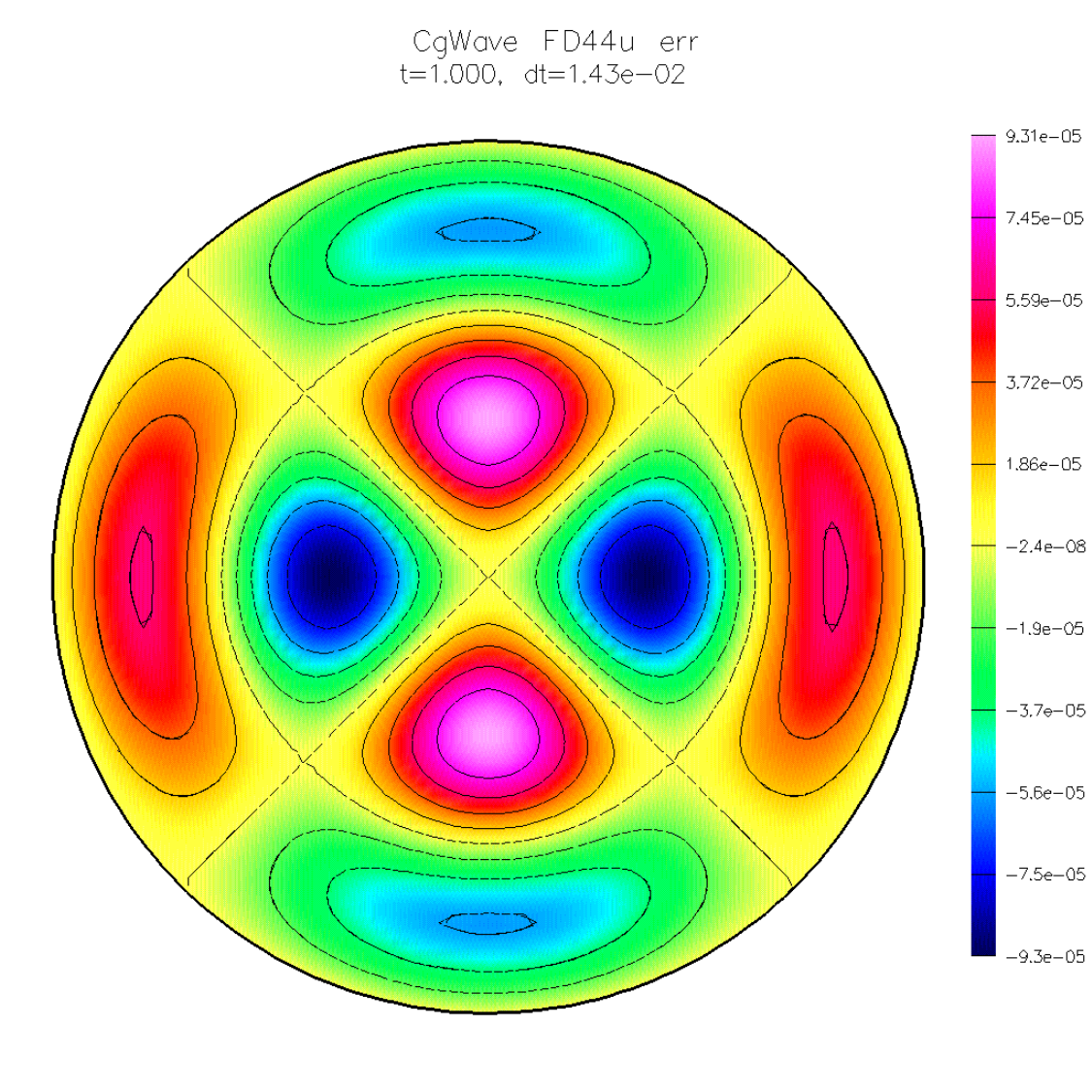}{err, $t=1$}{$v$}{$t=1.0$}{$-9.3e-5$}{$9.3e-5$}     
   \end{scope}

\end{tikzpicture}
\end{center}
\caption{Results for the disk.
    Left: overset grid $\Gcd^{(2)}$ for a disk. 
    Middle: computed eigenfunction $(\mTheta,\mr)=(2,1)$.
   Right: error. Implicit time-stepping, order four, grid $\Gcd^{(4)}$.
    }
\label{fig:diskEigenfunctionSolution}
\end{figure}
}

Figure~\ref{fig:diskFig} shows grid convergence results. Numerical solutions are computed using a time-step $\dt=.04/j$ for grids $\Gcd^{(j)}$, $j=2,4,8,16$.  Max-norm errors at $t=0.7$ are plotted as a function of the grid spacing.  The left-plot in the figure shows results for the eigenmode $(\mTheta,\mr)=(2,2)$ using the IME-UW-PC scheme for both Dirichlet and Neumann cases.  The results show that the numerical solutions are converging at close to the expected rates (as indicated by the reference lines in the log-log plots).  The right-plot in the figure shows results for the same eigenmode, but using the SPIE-UW-PC scheme with time-step determined by the explicit grid.  
As with the case of the fully implicit scheme, the results show that the numerical solutions are converging at close to the expected rates.

{
\newcommand{\figw}{5.5cm}
\newcommand{\figh}{5.5cm}
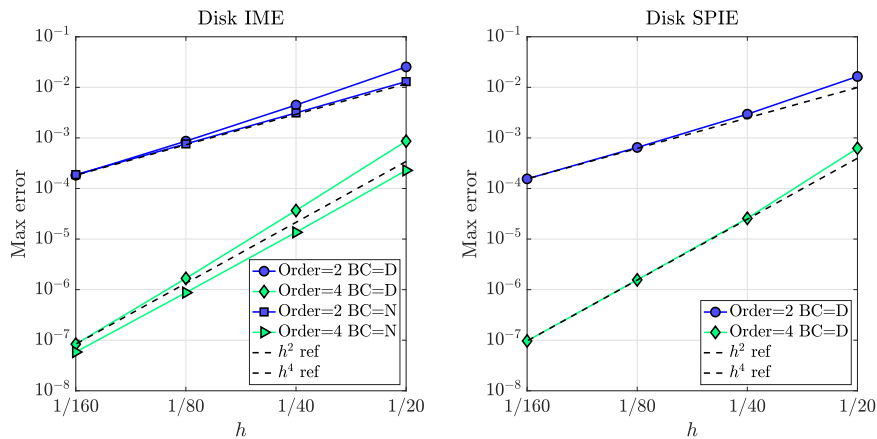
\begin{figure}[htb]
\begin{center}
\begin{tikzpicture}
  \useasboundingbox (0,.75) rectangle (11.5,5.5);  

  \begin{scope}[xshift=0cm]
     \figByWidth{0}{0}{diskIME}{\figw}[0][0][0][0]
     \figByWidth{6}{0}{diskSPIE}{\figw}[0][0][0][0]
  \end{scope} 
\end{tikzpicture}
\end{center}
\caption{
Left: Grid convergence for the IME-UW-PC scheme. Right: Grid convergence for the SPIE-UW-PC scheme.
 }
\label{fig:diskFig}
\end{figure}
}

\newcommand{\Gcscat}{\Gc_{\rm scat}}
\subsubsection{Scattering from a 2D cylinder} \label{sec:cylinderScattering}

We consider the scattering of a plane wave from a cylinder of radius $a$ in two dimensions.
The incident field is taken to be
\ba
   u_{\rm inc}(\xv,t) =  e^{i k (x-c t)} ,
\ea
where $k$ is the wave number of the incident field in the reference direction given by $x$.
The exact solution is written in polar coordinates $(r,\theta)$ with the usual assignment $x=r\cos\theta$.
A homogeneous Dirichlet boundary condition on the cylinder is assumed so that the total field (incident plus scattered) is given by
\bse
\label{eq:scatteringSoln}
\ba
  u(r,\theta,t) &= e^{- i k c t} \, \sum_{m=0}^\infty \eps_m \, i^m \, \left[  J_m(k r) - \f{J_m(ka)}{H_m^{(1)}( ka ) } H_m^{(1)}( kr ) \right] \cos(m \theta) , \\
    &= e^{ i k (x-c t) }   -  e^{- i k c t} \, \sum_{m=0}^\infty \eps_m \, i^m \, \left[ \f{J_m(ka)}{H_m^{(1)}( ka ) } H_m^{(1)}( kr ) \right] \cos(m \theta) , 
\ea
\ese
where $\eps_0=1$ and $\eps_m=2$ for $m>0$, and $H_m^{(1)}(z) = J_m(z) + i Y_m(z)$ is the Hankel function of the first kind of order $m$ defined in terms of the Bessel functions of the first and second kind.
Real-valued solutions are obtained by using either the real or imaginary parts of the solutions in~\eqref{eq:scatteringSoln}.

Numerical solutions are computed using an overset grid, denoted by $\Gcscat^{(j)}$, consisting of two component grids,
a background Cartesian grid covering $[-2,2]^2$ and an annular grid with inner radius $a=0.5$ and outer radius $b=0.8$.
The inner radius represents the cylindrical scatterer with a homogeneous Dirichlet boundary condition applied there, and the boundary conditions on the outer boundaries of the Cartesian grid are set to the exact solution.  The target grid spacing is approximately equal to $\ds^{(j)}=1/(10 j)$ in all directions.

{
\newcommand{\drawContour}[7]{%
\begin{scope}[#1]
\draw(0.0,0) node[anchor=south west,xshift=-4pt,yshift=+0pt] {\trimfiga{#2}{\figWidtha}};
  \draw(.75,.5) node[draw,fill=white,anchor=west,xshift=2pt,yshift=1pt,inner sep=2pt] {\scriptsize #3};
\begin{scope}[xshift=+.2cm,yshift=-2pt]
  \draw (\xcb,\ycb) node[anchor=south west,xshift=0.25cm,yshift=.5cm,rotate=-90] {\trimfigcb{colourBarLines}{\cbWidth}{\cbHeight}};
  \draw (.8,0) node[anchor=north,xshift=+3pt,yshift=+2pt] {\scriptsize $#6$};
  \draw (4.8,0) node[anchor=north,xshift=+0pt,yshift=+2pt] {\scriptsize $#7$};
\end{scope}
\end{scope}
}
\newcommand{\cbWidth}{.2cm}
\newcommand{\cbHeight}{4cm}
\newcommand{\xcb}{.5cm}
\newcommand{\ycb}{-.2cm}
\setlength{\ycbTop}{\ycb+\cbHeight}
\setlength{\ycbMid}{\ycb+\cbHeight*\real{.5}}
\newcommand{\trimfigcb}[3]{\includegraphics[width=#2, height=#3, clip, trim=17cm 2.35cm 1.65cm 2.35cm]{#1}}
\newcommand{\figWidtha}{5cm}
\newcommand{\trimfiga}[2]{\trimw{#1}{#2}{.0}{.117}{.09}{.09}}
\begin{figure}[htb]
\begin{center}
\begin{tikzpicture}
   \useasboundingbox (0,.3) rectangle (16,4.5);  

   \figByWidth{0}{-.05}{cicGridCoarse}{4.75cm}[0.1][0.1][0.1][0.1]

   \begin{scope}[xshift=4.75cm]
    \drawContour{xshift= 0.cm,yshift=0.00cm}{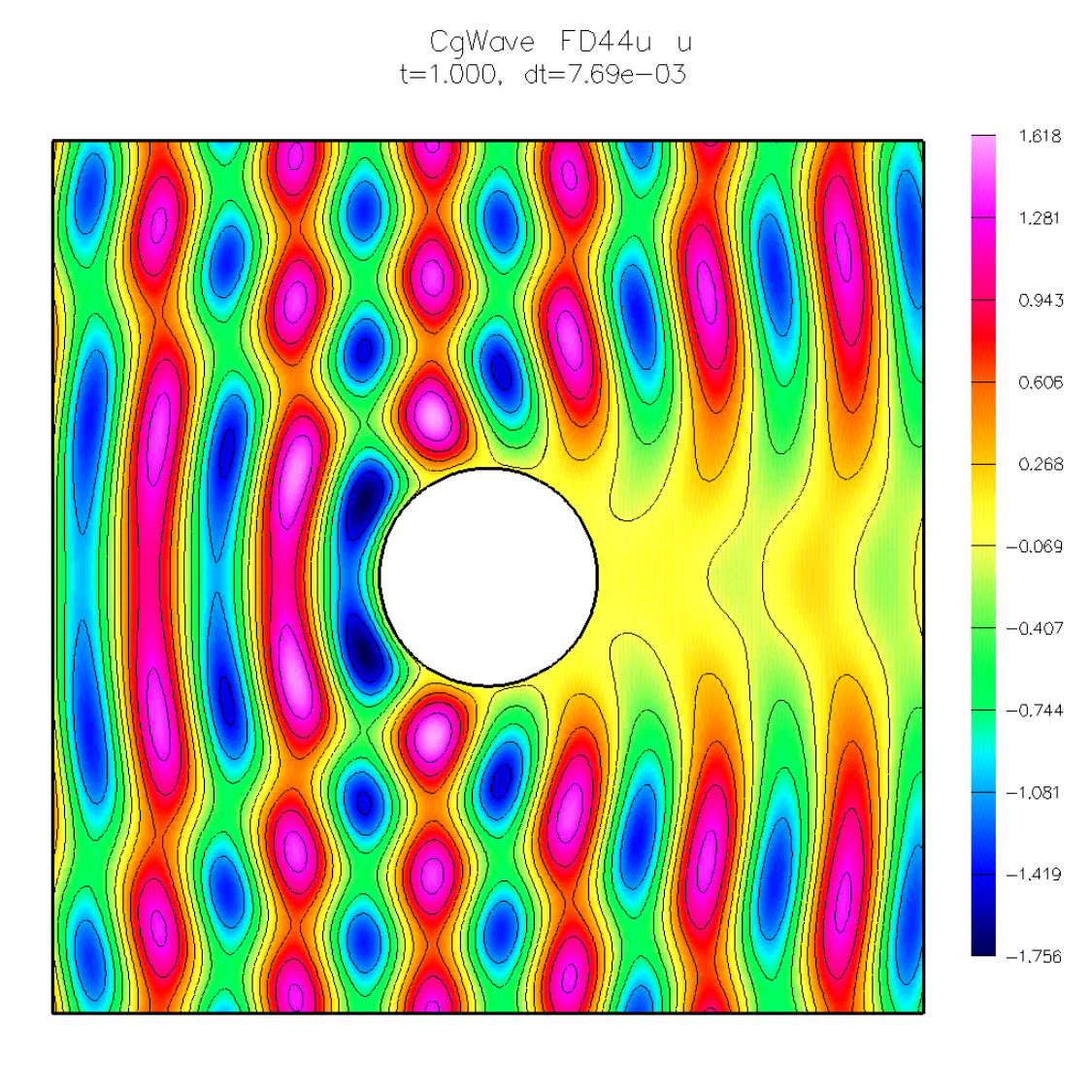}{u, $t=1.0$}{$v$}{$t=0.3$}{$-1.76$}{$1.62$}     
    \drawContour{xshift=5.5cm,yshift=0.00cm}{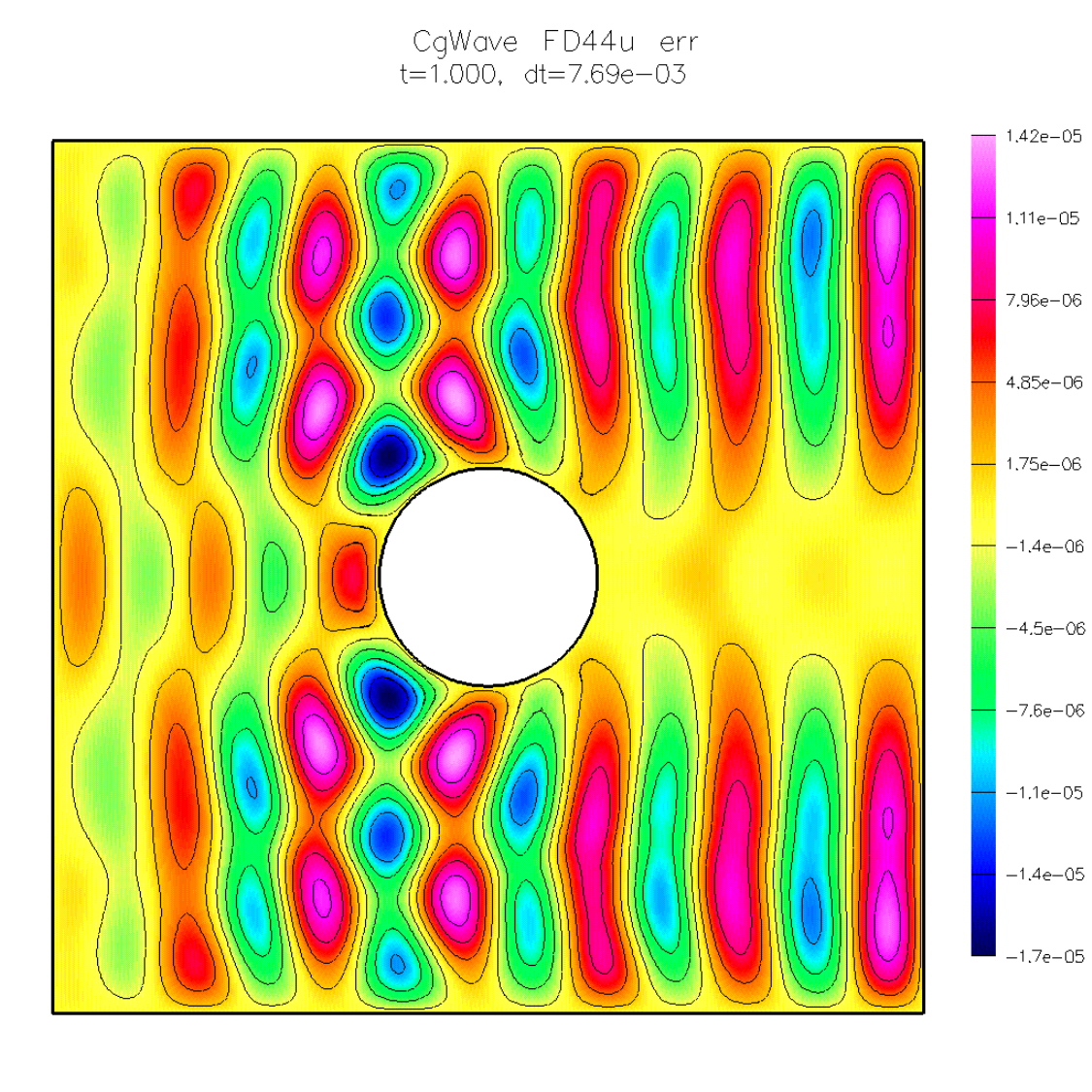}{err, $t=1.0$}{$v$}{$t=0.3$}{$-1.7e-5$}{$1.4e-5$}     
   \end{scope}

\end{tikzpicture}
\end{center}
\caption{Scattering from a cylinder.
  Left: overset grid for scattering from a cylinder. Middle and right: solution and errors for SPIE-UW-PC, order four, grid $\Gcscat^{(8)}$, $k=10$.
    }
\label{fig:cylScatContours}
\end{figure}
}

Figure~\ref{fig:cylScatContours} shows the overset grid $\Gcscat^{(2)}$ and contours of the computed solution and errors at $t=1$ using the fourth-order accurate SPIE-UW-PC scheme. The grid $\Gcscat^{(8)}$ is used for this calculation with $k=10$.
The errors are seen to be smooth.

{
\newcommand{\figw}{5.5cm}
\newcommand{\figh}{5.5cm}
\begin{figure}[htb]
\begin{center}
\begin{tikzpicture}
  \useasboundingbox (0,.75) rectangle (11.5,5.5);  

  \begin{scope}[xshift=0cm]
     \figByWidth{0}{0}{cylScatIME}{\figw}[0][0][0][0]
     \figByWidth{6}{0}{cylScatSPIE}{\figw}[0][0][0][0]
  \end{scope} 
\end{tikzpicture}
\end{center}
\caption{
Left: Grid convergence for the IME-UW-PC scheme. Right: Grid convergence for the SPIE-UW-PC scheme.
 }
\label{fig:cyScatFig}
\end{figure}
}

Figure~\ref{fig:cyScatFig} shows grid convergence results at $t=0.4$ for an incident field with
$k=2$ using the second and fourth-order accurate IME-UW-PC and SPIE-UW-PC schemes.
Max-norm errors are plotted as a function
of the grid spacing and the solutions are seen to be converging at close to the expected rates.

\newcommand{\Gcc}{\Gc_c}
\newcommand{\mz}{{m_z}}
\subsection{Cylinder eigenmodes with implicit time-stepping} \label{sec:cylEigenModes}

{
\newcommand{\figw}{5.5cm}
\newcommand{\figh}{5.5cm}
\begin{figure}[htb]
\begin{center}
\begin{tikzpicture}
  \useasboundingbox (0,.2) rectangle (11,5.25);  

  \begin{scope}[xshift=0cm,yshift=-.5cm]
     \figByWidth{0}{0}{cylinderIME}{\figw}[0][0][0][0]
     \figByWidth{6}{0}{cylinderSPIE}{\figw}[0][0][0][0]
  \end{scope} 
\end{tikzpicture}
\end{center}
\caption{
Left: Grid convergence for the IME-UW-PC scheme. 
Right: Grid convergence for the SPIE-UW-PC scheme.
 }
\label{fig:cylinderFig}
\end{figure}
}

In this section eigenmodes of a three-dimensional cylindrical solid pipe are computed.
The pipe has a radius of 0.5 and extends in the axial direction $z$ from $0$ to $1$.
In cylindrical coordinates $(r,\theta,z)$, the eigenmodes, for Dirichlet boundary conditions, take the form 
\ba
     u_{\mTheta,\mr,\mz} (r,\theta,t) = J_{\mTheta}( k_{\mTheta,\mr} r ) \,  e^{i \mTheta \theta} \, \sin( \mz \pi z ) \, e^{i \omega t} .
\ea

The composite grid for the solid cylinder, denoted by $\Gcc^{(j)}$, consists of two component grids,
each with grid spacings approximately equal to $\ds^{(j)}=1/(10 j)$ in all directions.
One component grid is a boundary-fitted cylindrical shell, while the other component grid is a background Cartesian grid covering the interior of the cylindrical domain (see Figure~\ref{fig:cylinderFig}).

{
\newcommand{\drawContour}[7]{%
\begin{scope}[#1]
\draw(0.0,0) node[anchor=south west,xshift=-4pt,yshift=+0pt] {\trimfiga{#2}{\figWidtha}};
  \draw(.75,.5) node[draw,fill=white,anchor=west,xshift=2pt,yshift=1pt,inner sep=2pt] {\scriptsize #3};
\begin{scope}[xshift=+.2cm,yshift=-2pt]
  \draw (\xcb,\ycb) node[anchor=south west,xshift=0.25cm,yshift=.5cm,rotate=-90] {\trimfigcb{colourBarLines}{\cbWidth}{\cbHeight}};
  \draw (.8,0) node[anchor=north,xshift=+3pt,yshift=+2pt] {\scriptsize $#6$};
  \draw (4.8,0) node[anchor=north,xshift=+0pt,yshift=+2pt] {\scriptsize $#7$};
\end{scope}
\end{scope}
}
\newcommand{\cbWidth}{.2cm}
\newcommand{\cbHeight}{4cm}
\newcommand{\xcb}{.5cm}
\newcommand{\ycb}{-.2cm}
\setlength{\ycbTop}{\ycb+\cbHeight}
\setlength{\ycbMid}{\ycb+\cbHeight*\real{.5}}
\newcommand{\trimfigcb}[3]{\includegraphics[width=#2, height=#3, clip, trim=17cm 2.35cm 1.65cm 2.35cm]{#1}}
\newcommand{\figWidtha}{5.5cm}
\newcommand{\trimfiga}[2]{\trimw{#1}{#2}{.0}{.117}{.09}{.09}}
\begin{figure}[htb]
\begin{center}
\begin{tikzpicture}
   \useasboundingbox (0.2,.3) rectangle (16.5,5);  

   \figByWidth{0.0}{0.15}{pipeGridG2}{5.6cm}[0.][0.][0.05][0.05]

   \begin{scope}[xshift=5.1cm]
    \drawContour{xshift= 0.cm,yshift=0.00cm}{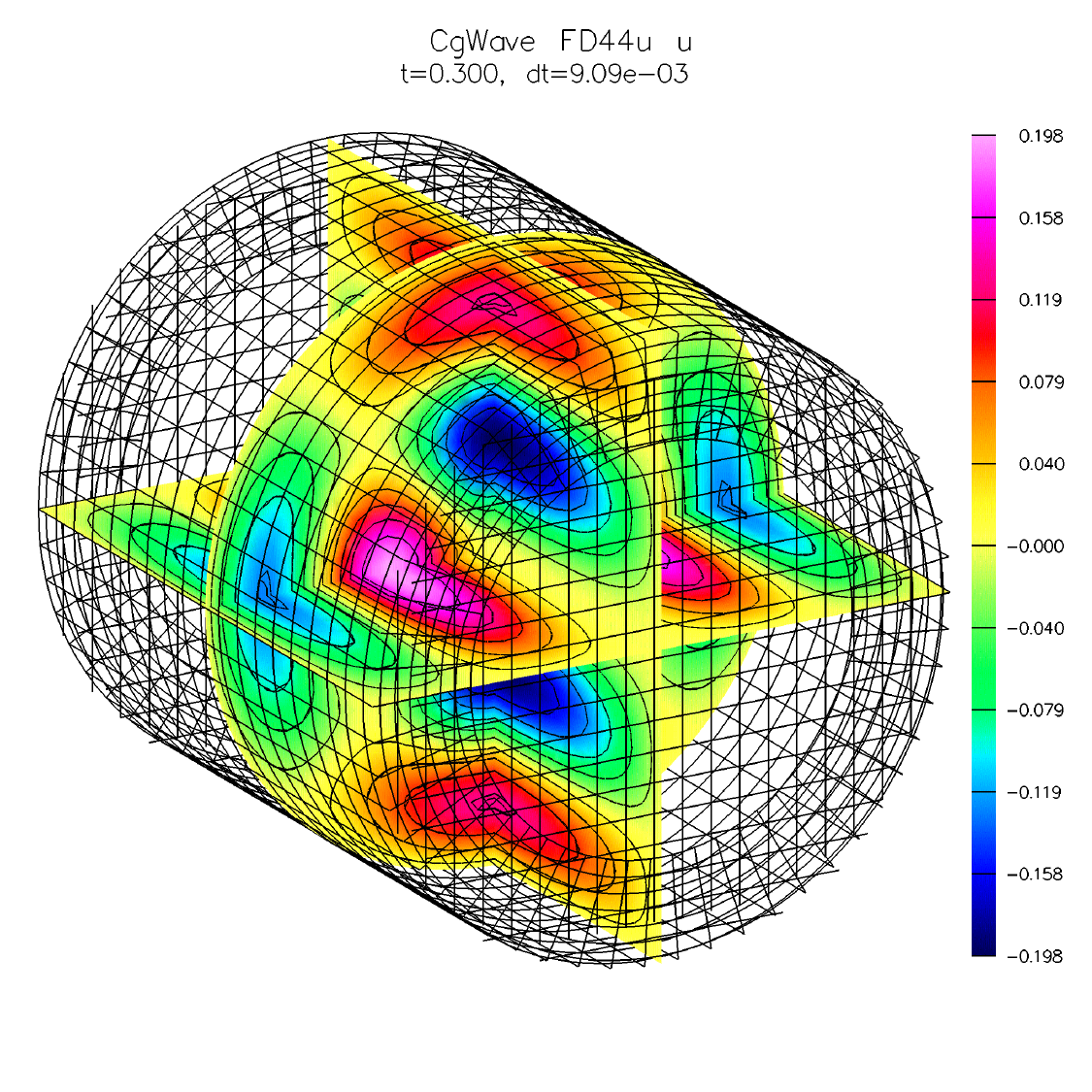}{u, $t=0.4$}{$v$}{$t=0.3$}{$-.41$}{$.41$}    
    \drawContour{xshift= 5.5cm,yshift=0.00cm}{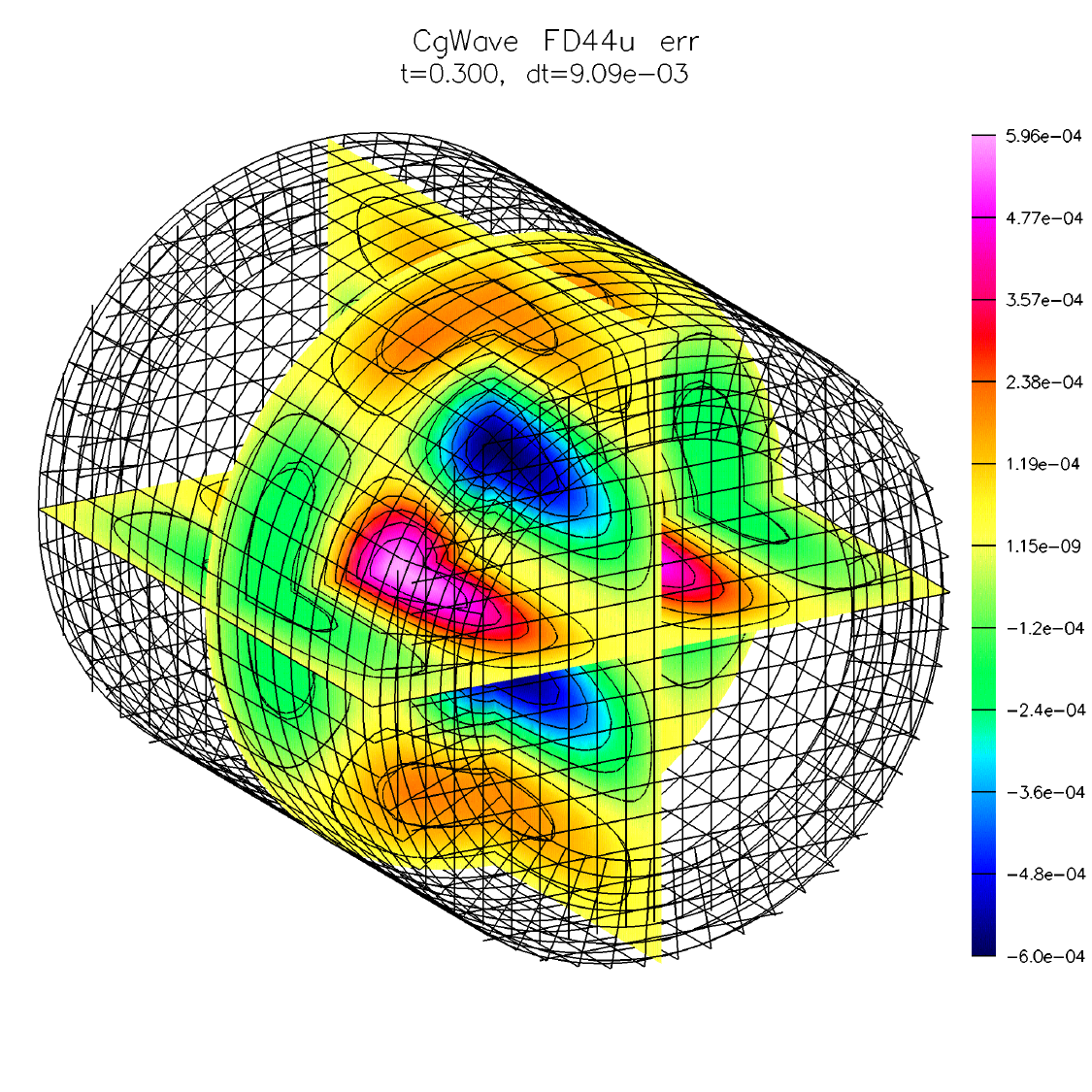}{err, $t=0.4$}{$v$}{$t=0.3$}{$-2.3e-4$}{$2.3e-4$}     
   \end{scope}

\end{tikzpicture}
\end{center}
\caption{Cylinder.
  Left: overset grid $\Gcc^{(2)}$ for a solid cylinder.
  Middle and right: Solution and errors for IME-UW-PC, order four, grid $\Gcc^{(4)}$, eigenmode $(\mTheta,\mr,\mz)=(2,1,1)$. A coarsened version of the grid is shown.
    }
\label{fig:cylEigs}
\end{figure}
}

Figure~\ref{fig:cylinderFig} shows grid convergence results. Max-norm errors are plotted as a function
of the grid spacing. 
Results are shown for IME-UW-PC and SPIE-UW-PC schemes for Dirichlet at $t=0.6$.
The eigenmode was $(\mTheta,\mr,\mz)=(1,2,2)$. For Dirichlet boundary conditions this
corresponded to a value of $k_{\mTheta,\mr}=8.41724414$.
The time-step for the SPIE scheme was chosen as $\dt=.04/j$ for grid $\Gcd^{(j)}$, $j=2,4,8$.
The results show that the computed solution is converging at close to the expected rates.
Figure~\ref{fig:cylEigs} shows contours of the solution and errors for $(\mTheta,\mr,\mz)=(2,2,2)$.

\newcommand{\Gcs}{\Gc_s}
\newcommand{\mPhi}{{m_{\phi}}}
\subsubsection{Eigenmodes on a sphere} \label{sec:sphereEigenModes}

We now consider eigenmode solutions of a solid unit sphere assuming a homogeneous Dirichlet boundary condition on the surface of the sphere.
Introduce spherical polar coordinates $(r,\theta,\phi)$, where $r$ is the radius, $\theta\in[0,2\pi]$ is the angle in the $x$-$y$ plane and $\phi\in[0,\pi]$ the angle from
the $z$-axis. 
We assume time-periodic eigenmodes with frequency $\omega$ having the well known form 
\ba
    u_{\mr,\mTheta,\mPhi} (r,\theta,\phi,t) 
       = r^{-1/2}\, J_{\mPhi +\half}( \lambda_{\mPhi,\mr} r)  \, P_{\mPhi}^{\mTheta}(\cos\phi) \, e^{i \mTheta \theta} \, e^{i\omega t}, 
\ea
where $J_{\mPhi +\half}$, $\mPhi=0,1,2,\ldots$, are Bessel functions of fractional order, $P_{\mPhi}^{\mTheta}$, $\mPhi\ge\mTheta$, are associated Legendre functions, and $\lambda_{\mPhi,\mr}$, $\mr=1,2,\ldots$, are zeros of $J_{\mPhi +\half}$.  The frequency of vibration is given by $\omega=c \, \lambda_{\mPhi,\mr}$.

\begin{figure}[htb]
\begin{center}
\begin{tikzpicture}
   \useasboundingbox (0,.75) rectangle (11.5,5.5);  
  \figByHeight{0.00}{0}{sphereIME}{5.5cm}[0][0][0.][0]
  \figByHeight{6.25}{0}{sphereSPIE}{5.5cm}[0][0][0.][0]
\end{tikzpicture}
\end{center}
\caption{
Left: grid convergence results for the IME-UW-PC scheme.
Right: grid convergence results for the SPIE-UW-PC scheme.
}
\label{fig:solidSphereConvergence}
\end{figure}
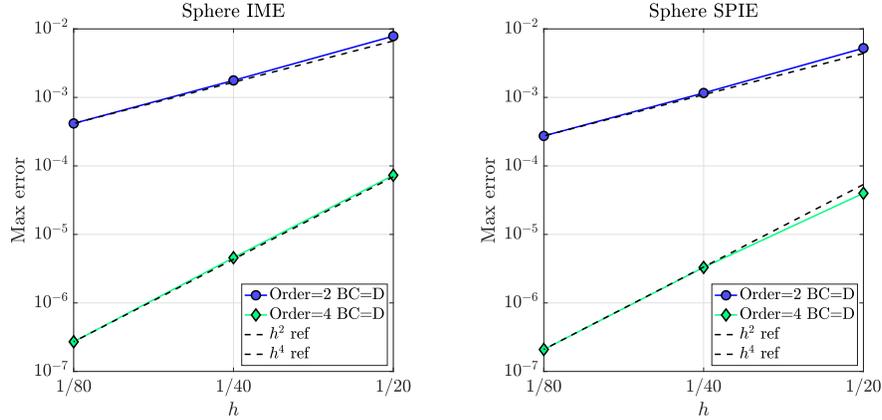

The composite grids for the solid sphere of radius one, denoted by $\Gcs^{(j)}$, consist of four component grids,
each with grid spacing approximately equal to $\ds^{(j)}=1/(10 j)$.  
The sphere is covered with three boundary-fitted patches near the surface as shown on the left in Figure~\ref{fig:solidSphereConvergence}.
There is one patch specified using spherical polar coordinates that covers much of the sphere except
near the poles. To remove the polar singularities there are two patches
that cover the north and south poles, defined by orthographic mappings.
A background Cartesian grid (not shown) covers the interior of the sphere.  The middle image in the figure shows the solution at $t=0.5$  for the eigenmode with 
$(m_\phi,m_\theta,m_r)=(2,1,1)$ 
and $\lambda_{\mPhi,\mr}\approx 5.7634591968945$.  This solution is computed using the fourth-order accurate SPIE-UW-PC scheme and the grid $\Gcs^{(4)}$.  The right image shows the max errors which are smooth.

{
\newcommand{\drawContour}[7]{%
\begin{scope}[#1]
\draw(0.0,0) node[anchor=south west,xshift=-4pt,yshift=+0pt] {\trimfiga{#2}{\figWidtha}};
  \draw(.75,.5) node[draw,fill=white,anchor=west,xshift=2pt,yshift=1pt,inner sep=2pt] {\scriptsize #3};
\begin{scope}[xshift=-.1cm,yshift=-2pt]
  \draw (\xcb,\ycb) node[anchor=south west,xshift=0.25cm,yshift=.5cm,rotate=-90] {\trimfigcb{colourBarLines}{\cbWidth}{\cbHeight}};
  \draw (.8,0) node[anchor=north,xshift=+3pt,yshift=+2pt] {\scriptsize $#6$};
  \draw (4.8,0) node[anchor=north,xshift=+0pt,yshift=+2pt] {\scriptsize $#7$};
\end{scope}
\end{scope}
}
\newcommand{\cbWidth}{.2cm}
\newcommand{\cbHeight}{4cm}
\newcommand{\xcb}{.5cm}
\newcommand{\ycb}{-.2cm}
\setlength{\ycbTop}{\ycb+\cbHeight}
\setlength{\ycbMid}{\ycb+\cbHeight*\real{.5}}
\newcommand{\trimfigcb}[3]{\includegraphics[width=#2, height=#3, clip, trim=17cm 2.35cm 1.65cm 2.35cm]{#1}}
\newcommand{\figWidtha}{4.75cm}
\newcommand{\trimfiga}[2]{\trimw{#1}{#2}{.04}{.14}{.06}{.09}}
\begin{figure}[htb]
\begin{center}
\begin{tikzpicture}
   \useasboundingbox (0,.3) rectangle (16,5.);  

   \figByWidth{0}{0.1}{sphereGridsExploded}{5.9cm}[0][0][0.1][0.1]

   \begin{scope}[xshift=5.75cm]
    \drawContour{xshift= 0.cm,yshift=0.00cm}{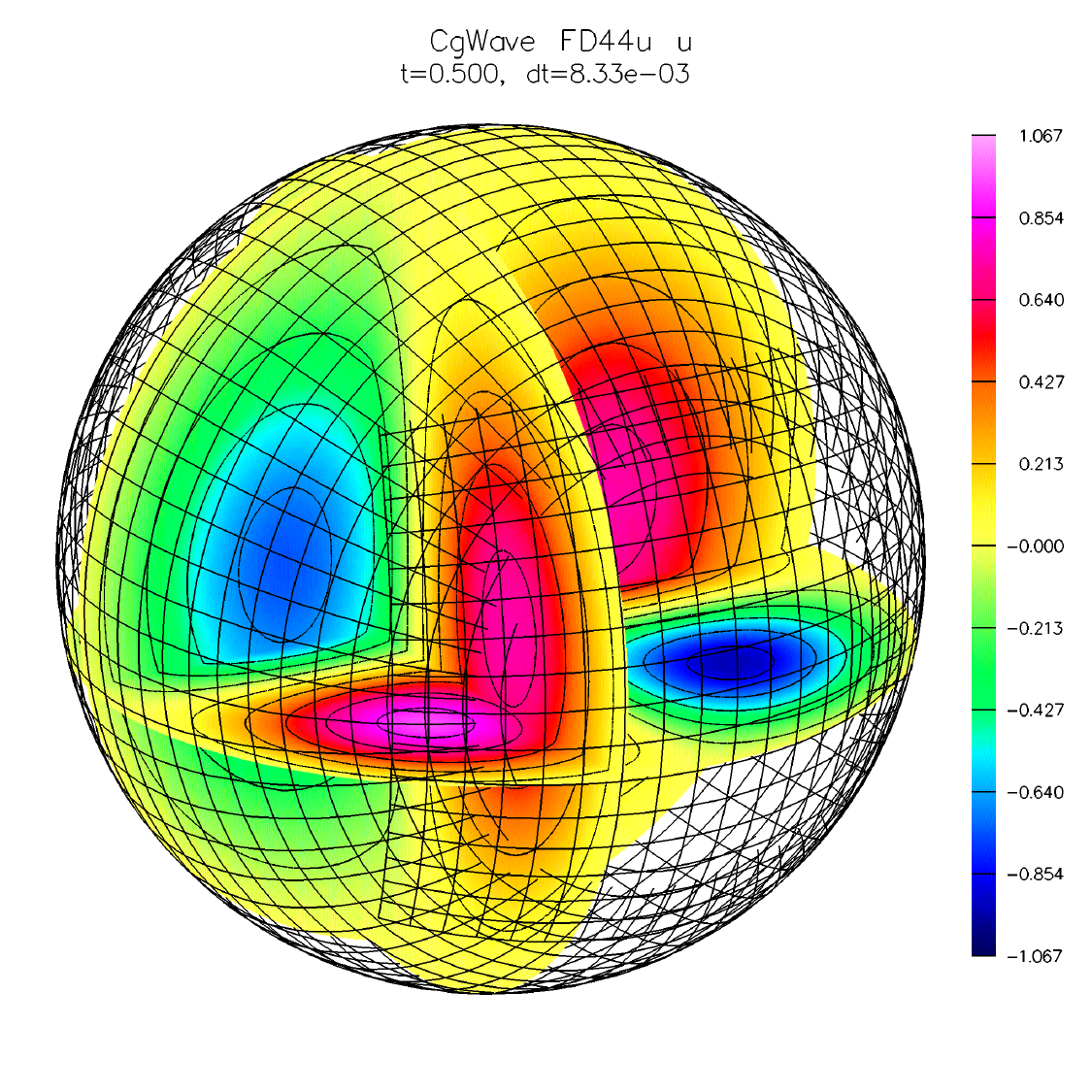}{u, $t=0.5$}{$v$}{$t=0.3$}{$-1.1$}{$1.1$}    
    \drawContour{xshift=5.cm,yshift=0.00cm}{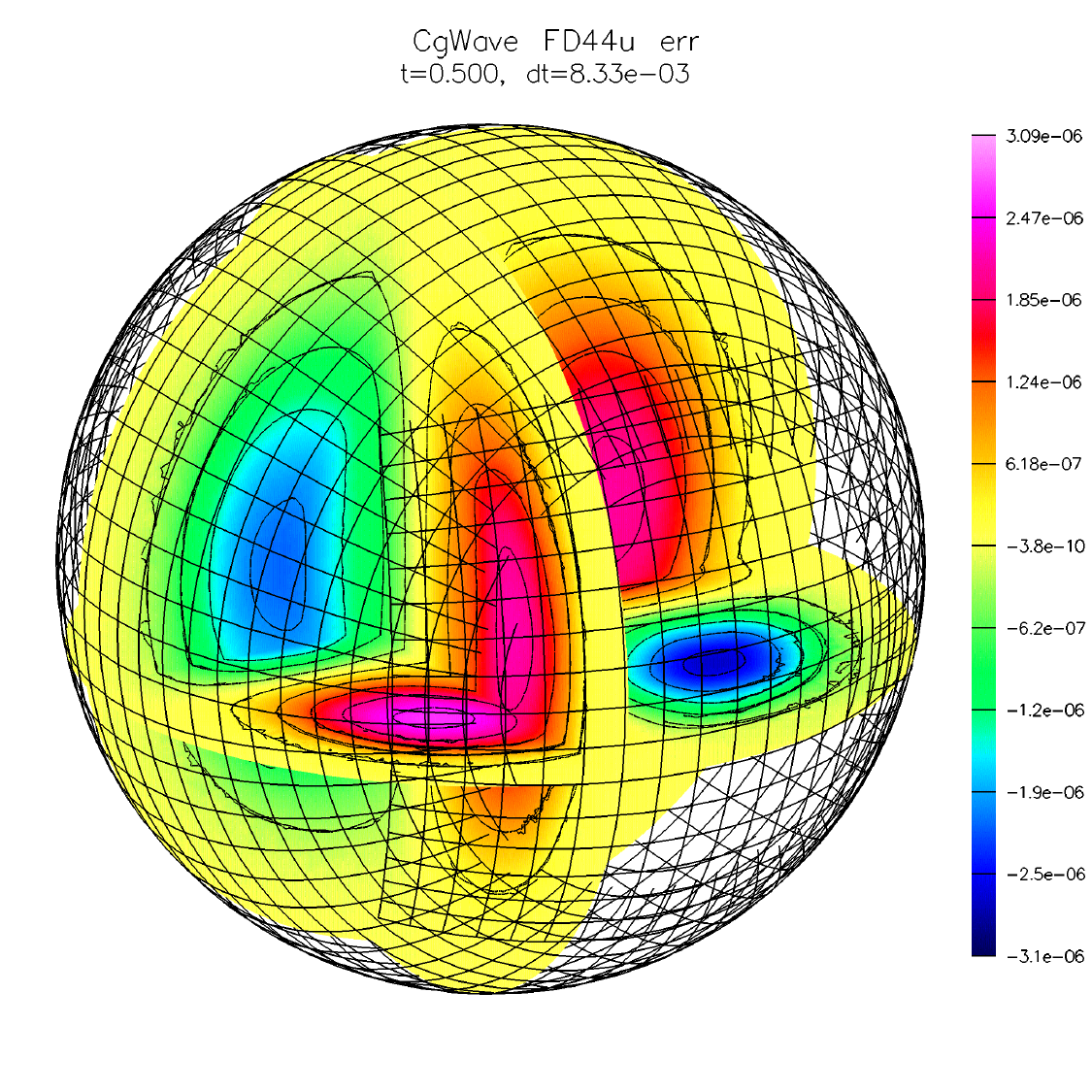}{err, $t=0.5$}{$v$}{$t=0.3$}{$-3.1e-6$}{$3.1e-6$}    
   \end{scope}

\end{tikzpicture}
\end{center}
\caption{Sphere eigenmodes.
 Left: exploded view of the surface patches of the overset grid for the interior of a sphere.
 Middle and right: Computed solution and errors for the fourth-order accurate SPIE-UW-PC scheme on grid $\Gcs^{(4)}$,
 eigenmode $(m_\phi,m_\theta,m_r)=(2,1,1)$.
 A coarsened version on the grid is shown.
    }
\label{fig:spherelEigs}
\end{figure}
}

Figure~\ref{fig:solidSphereConvergence} shows grid convergence results for
 the eigenmode $(m_\phi,m_\theta,m_r)=(2,1,1)$ at $t=0.4$.
 Results are shown for the second and fourth-order accurate IME-UW-PC and SPIE-UW-PC schemes.
 The graphs demonstrate that the solutions are converging at close to the expected rates.

\subsubsection{Long time simulations with random initial conditions} \label{sec:longTimeSimulations}

In this section we perform some very long-time simulations to confirm numerically that
the solutions computed using the IME and SPIE schemes with upwinding remain stable and bounded. 
Initial conditions are chosen with random grid values on $[0,1]$  so that all eigenmodes, 
including any possible unstable ones, would be seeded with an order one amount of energy.
The numerical schemes are integrated to very long times and the solutions are monitored for
any growth. Due to the upwinding, the magnitude of the computed solutions for a stable scheme is expected to decay slowly to
zero over time.

To assess the growth or decay of the solution we plot a discrete approximation to the energy given by
\ba
   \Ec(t) = \half \Big( \| \p_t u \|_{\Omega}^2 + c^2 \| \grad u \|_{\Omega}^2  \Big)  , \label{eq:energy}
\ea
where $\| \cdot \|_{\Omega}$ denotes the $L_2-$norm over the domain $\Omega$. We note that the energy defined in~\eqref{eq:energy}
remains constant in time for exact solutions of the wave equation on $\Omega$ assuming homogeneous Dirichlet or Neumann conditions specified
on the boundary of $\Omega$.
For purposes of this study, a first order accurate approximation to~\eqref{eq:energy} is sufficient, denoted by $\Ec_h$.
$\Dmt U_\iv^n$ is used to approximate the time derivative in~\eqref{eq:energy} and first order accurate backward differences
are used to approximate the spatial derivatives, for example on a Cartesian grid $\p_x u \approx \Dmx U_\iv^n$.
Note that the discrete energy $\Ec_h$ would remain approximately constant if the scheme is stable, but with upwind dissipation included the discrete energy
is expected to decay over time.

{
\newcommand{\figw}{5.4cm}
\newcommand{\figh}{5.4cm}
\begin{figure}[htb]
\begin{center}
\begin{tikzpicture}
  \useasboundingbox (0,.75) rectangle (16.5,5);  

  \begin{scope}[xshift=0cm]
     \figByWidth{0}{0}{randomDiskNorms}{\figw}[0][0][0][0]
     \figByWidth{5.5}{0}{randomSphereNorms}{\figw}[0][0][0][0]
     \figByWidth{11}{0}{randomDiskCFLNorms}{\figw}[0][0][0][0]
  \end{scope} 

\end{tikzpicture}
\end{center}
\caption{Long time simulations. Left: disk, SPIE-UW-PC. Middle: sphere, SPIE-UW-PC. Right: disk, IME-UW-PC.
 }
\label{fig:longTimeSimulations}
\end{figure}
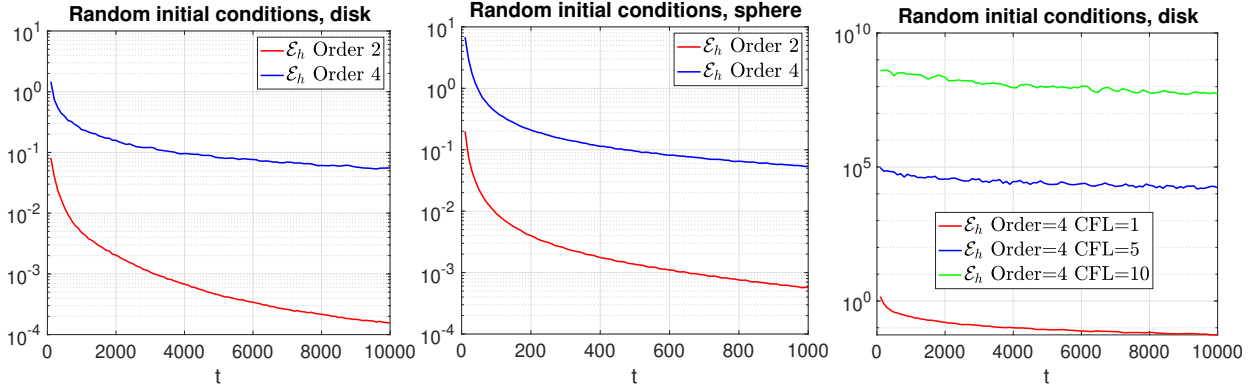
}

Figure~\ref{fig:longTimeSimulations} shows results from some long-time simulations for both the SPIE and IME schemes.
In all cases the schemes remained stable.
The left plot shows the discrete energy $\Ec_h(t)$ over time for a computation
on the disk grid $\Gcd^{(4)}$ as described in Section~\ref{sec:diskEigenModes}.
The final time is $t=10^4$ for the simulation and approximately $6\times 10^5$ time-steps are used.
Results are shown for the second and fourth-order accurate SPIE-UW-PC schemes.
The discrete energy is seen to decay rapidly at first as the high-frequency components of the solution are damped by the high-order upwind dissipation.
As time progresses the solution becomes smoother and the energy decays more slowly. 
The discrete energy for the fourth-order accurate scheme decays more slowly than the second-order accurate scheme since
its dissipation scales as $\Oc(h^5)$ compared to $\Oc(h^3)$ for the second-order accurate scheme.
The middle plot shows results for the three-dimensional solid sphere
grid $\Gcs^{(4)}$ as described in Section~\ref{sec:sphereEigenModes}.
In this case the final time is $t=10^3$ and the calculation requires approximately $10^5$ time-steps.
The results show that the discrete energy decays and schemes remain stable for the spherical case in qualitative agreement with the results for the disk case.

The right-most plot of Figure~\ref{fig:longTimeSimulations} compares the energy decay 
for the (fully implicit) IME4-UW-PC scheme on the disk grid $\Gcd^{(4)}$ for three different values 
of the CFL number, $1$, $5$ and $10$. In each case the scheme remains stable and
the discrete energy $\Ec_h$ decays. The dissipation parameter $\nu_p$ is the same
for each case. Note that the CFL=$10$ run takes $10$ times fewer time-steps than the CFL=$1$
run, and thus the dissipation has fewer time-steps to act.

\subsection{Performance of the SPIE scheme}

We now turn our attention to a set of examples that posses some geometric stiffness. For such problems
it is demonstrated that the SPIE scheme can be much faster than the fully explicit schemes.
Importantly, it is also shown that the accuracy of the computed solutions from the SPIE scheme are, in general,
quite similar to the accuracy of the explicit ME solutions. Thus, at least for the cases shown here,
taking a large CFL time-step using an
implicit ME method in small parts of the domain
where geometric stiffness occurs does not appear to have a significant effect on the overall accuracy.

\newcommand{\Gcsh}{\Gc_{\rm sh}}
\subsubsection{Scattering from a small hole} \label{sec:cylScatPerformance}

This section studies the accuracy and performance of the SPIE scheme for the scattering of a plane wave from 
a small cylinder in two space dimensions. The incident wave, exact solution, and overset grid topology
were described previously in Section~\ref{sec:cylinderScattering}. 
The results presented in this section show that the SPIE scheme can compute solutions much faster than
the EME scheme but with similar errors. 

{
\newcommand{\figw}{5.25cm}
\newcommand{\figh}{5.cm}
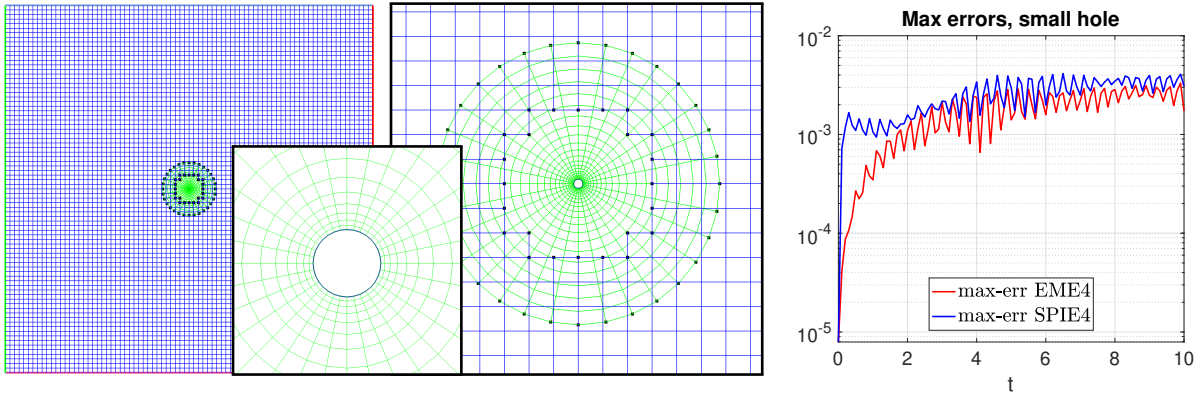
\begin{figure}[htb]
\begin{center}
\begin{tikzpicture}
  \useasboundingbox (0,.75) rectangle (16,5.1);  

  \figByWidth{-.2}{0.1}{holeGridRad0p01G2}{5.2cm}[0.][0.][0.][0.]
  \figByWidthb{5.1}{0.25}{holeGridRad0p01G2Zoom}{4.9cm}[0.05][0.05][0.05][0.05]

  \figByWidthb{3}{0.25}{holeGridRad0p01G2Zoom2}{3cm}[0.05][0.05][0.05][0.05]

  \begin{scope}[xshift=10.5cm]
     \figByWidth{0}{0}{smallHoleErrors}{\figw}[0][0][0][0]
  \end{scope} 

\end{tikzpicture}
\end{center}
\caption{Scattering from a small hole. Left: overset grid $\Gcscat^{(2)}$ and magnified views.
   Right: max-norm errors over time for the EME4 and SPIE4 schemes on grid $\Gcscat^{(4)}$. The SPIE4 scheme achieves similar errors to the 
   EME4 scheme but at a factor $21$ reduced CPU cost.
 }
\label{fig:smallHoleGridAndErrorsFig}
\end{figure}
}

{
\newcommand{\drawContour}[7]{%
\begin{scope}[#1]
\draw(0.0,0) node[anchor=south west,xshift=-4pt,yshift=+0pt] {\trimfiga{#2}{\figWidtha}};
  \draw(.75,.5) node[draw,fill=white,anchor=west,xshift=2pt,yshift=1pt,inner sep=2pt] {\scriptsize #3};
\begin{scope}[xshift=-.475cm,yshift=-2pt]
  \draw (\xcb,\ycb) node[anchor=south west,xshift=0.25cm,yshift=.5cm,rotate=-90] {\trimfigcb{colourBarLines}{\cbWidth}{\cbHeight}};
  \draw (.8,0) node[anchor=north,xshift=+3pt,yshift=+2pt] {\scriptsize $#6$};
  \draw (4.8,0) node[anchor=north,xshift=+0pt,yshift=+2pt] {\scriptsize $#7$};
\end{scope}
\end{scope}
}
\newcommand{\drawContourB}[7]{%
\begin{scope}[#1]
\draw(0.0,0) node[anchor=south west,xshift=-4pt,yshift=+0pt] {\trimfiga{#2}{\figWidtha}};
  \draw(.75,.5) node[draw,fill=white,anchor=west,xshift=2pt,yshift=1pt,inner sep=2pt] {\scriptsize #3};
\begin{scope}[xshift=+.2cm,yshift=-2pt]
  \draw (\xcb,\ycb) node[anchor=south west,xshift=0.25cm,yshift=.5cm,rotate=-90] {\trimfigcb{colourBarLines}{\cbWidth}{\cbHeight}};
  \draw (.8,0) node[anchor=north,xshift=+3pt,yshift=+2pt] {\scriptsize $#6$};
  \draw (4.8,0) node[anchor=north,xshift=+0pt,yshift=+2pt] {\scriptsize $#7$};
\end{scope}
\end{scope}
}
\newcommand{\cbWidth}{.2cm}
\newcommand{\cbHeight}{4cm}
\newcommand{\xcb}{.5cm}
\newcommand{\ycb}{-.2cm}
\setlength{\ycbTop}{\ycb+\cbHeight}
\setlength{\ycbMid}{\ycb+\cbHeight*\real{.5}}
\newcommand{\trimfigcb}[3]{\includegraphics[width=#2, height=#3, clip, trim=17cm 2.35cm 1.65cm 2.35cm]{#1}}
\newcommand{\figWidtha}{4cm}
\newcommand{\trimfiga}[2]{\trimw{#1}{#2}{.12}{.117}{.1}{.12}}
\begin{figure}[htb]
\begin{center}
\begin{tikzpicture}
   \useasboundingbox (0,.3) rectangle (15,4.25);  

   \begin{scope}[yshift=0cm]
    \drawContour{xshift= 0.cm,yshift=0.00cm}{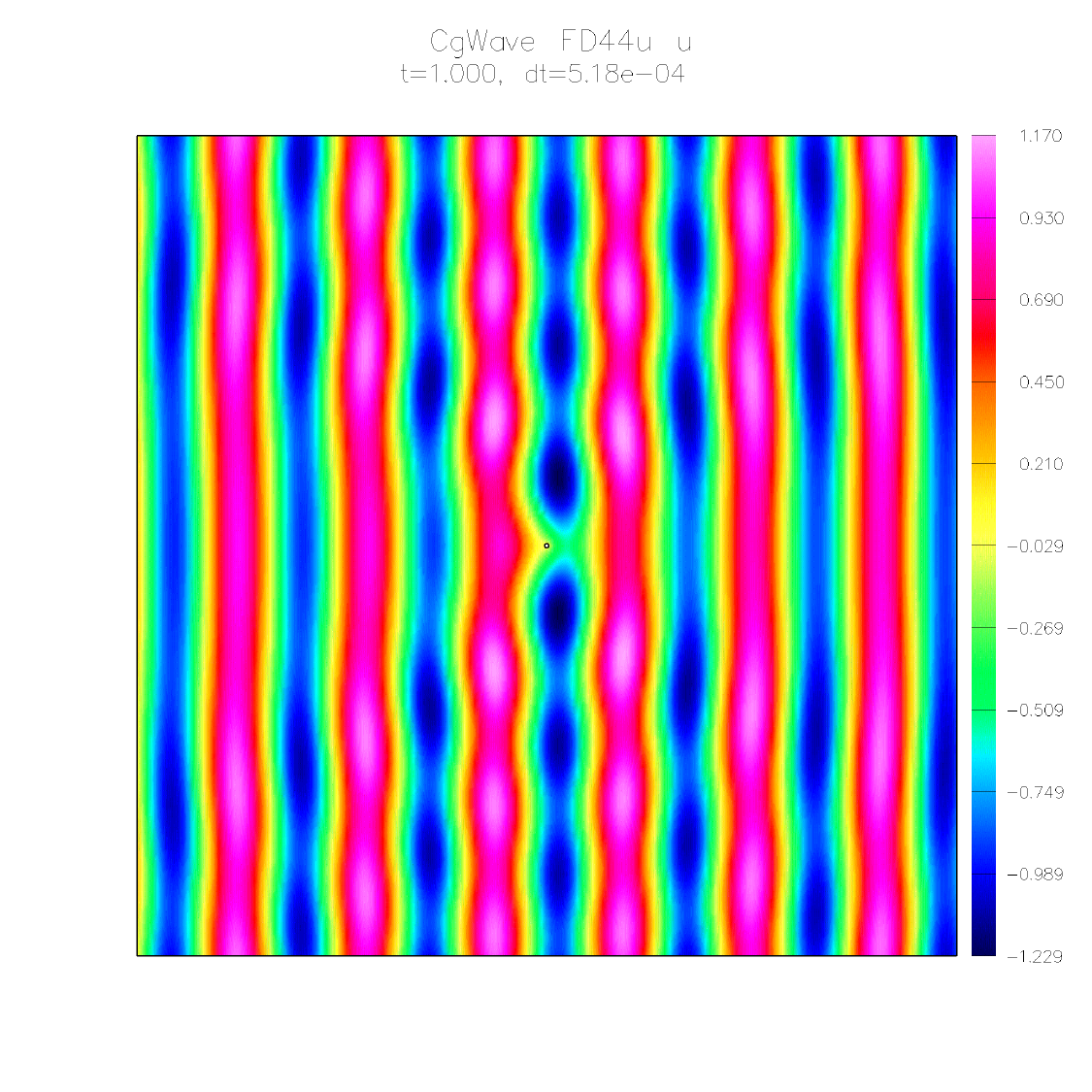}{$u$, $t=1.0$}{$v$}{$t=0.3$}{$-1.23$}{$1.17$} 
  
    \renewcommand{\figWidtha}{5.5cm}   
    \renewcommand{\trimfiga}[2]{\trimw{#1}{#2}{.0}{.117}{.1}{.3}} 
    \drawContourB{xshift=4.5cm,yshift=0.00cm}{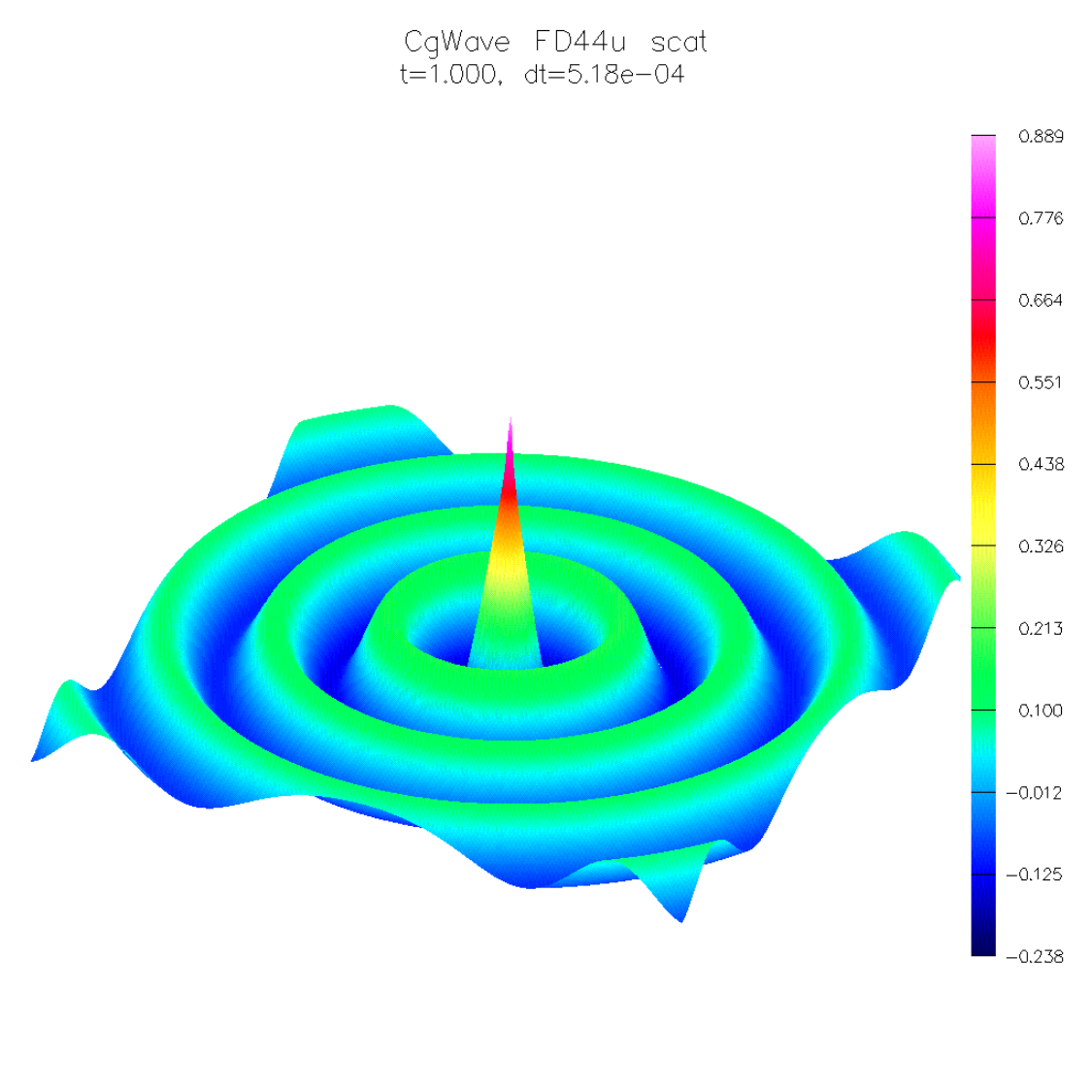}{$u_s$, $t=1.0$}{$v$}{$t=0.3$}{$-.24$}{$.89$} 

    \renewcommand{\figWidtha}{4cm}
    \renewcommand{\trimfiga}[2]{\trimw{#1}{#2}{.12}{.117}{.1}{.12}}    
    \drawContour{xshift=10.5cm,yshift=0.00cm}{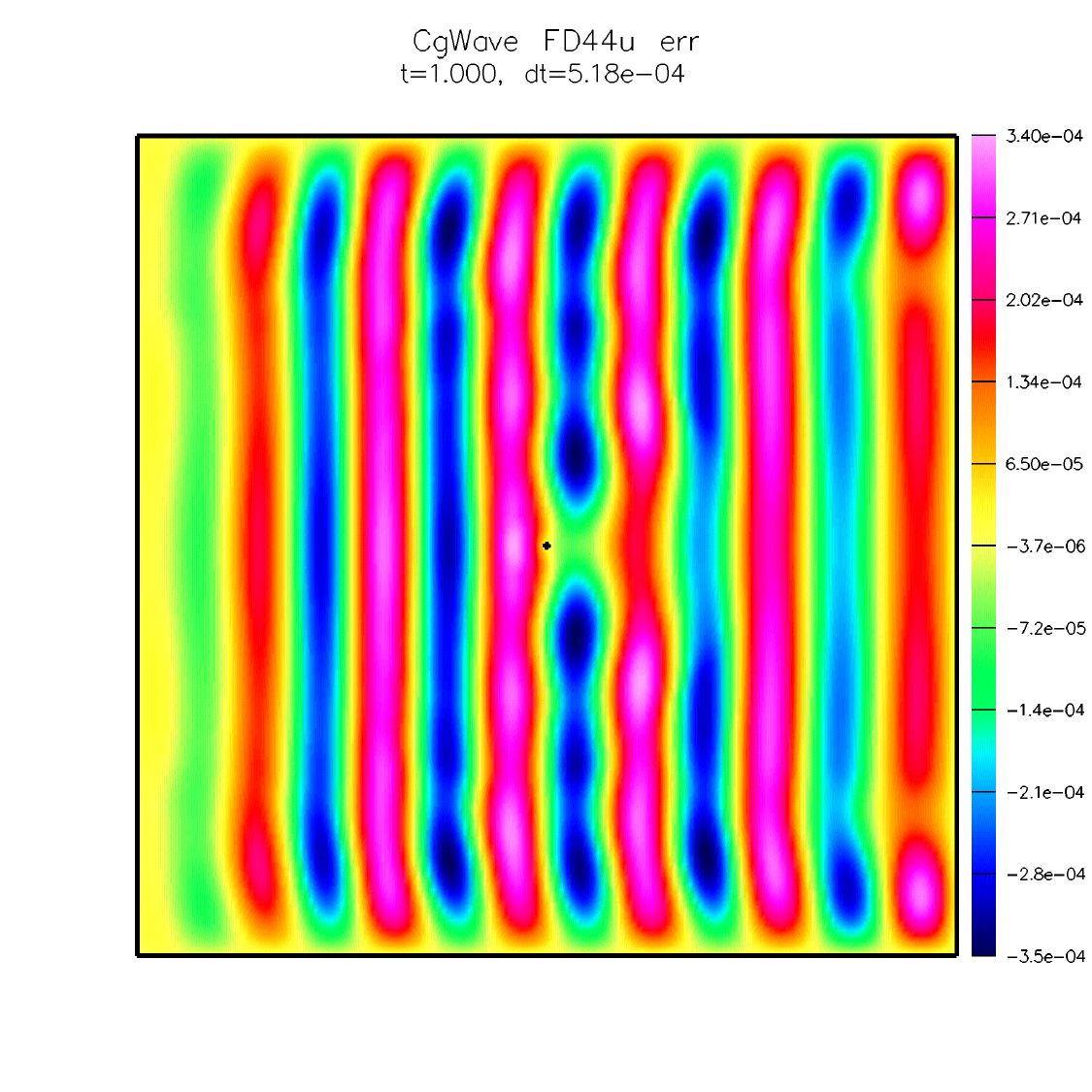}{err, $t=1.0$}{$v$}{$t=0.3$}{$-3.5e-4$}{$3.4e-4$}
    \draw(4.5,3.5) node[draw,fill=white,anchor=west,xshift=2pt,yshift=1pt,inner sep=3pt] {\scriptsize Small hole};   
    \draw[thick,black,->] (4.5,3.5) -- (2.5,2.3);            
   \end{scope}

\end{tikzpicture}
\end{center}
\caption{Scattering from a small hole. Left: contour plot of the total field $u$. Middle: surface plot of the
scattered field $u_s$. Right: errors in $u$. Solution at $t=1$ computed on grid $\Gcscat^{(4)}$ using the EME4-UW-PC scheme.
    }
\label{fig:smallHoleContours}
\end{figure}
}

A very small cylindrical hole of radius $a=0.01$ sits at the center of a square domain $[-2,2]^2$. The overset grid
is shown in Figure~\ref{fig:smallHoleGridAndErrorsFig}.
An incident field with wave-number $k=10$ impinges on the hole where a homogeneous Dirichlet boundary conditions is applied.
The exact solution is imposed on the outer boundaries of the square.
The solution is computed to a final time of $t=10$ using the EME4-UW-PC scheme and the SPIE4-UW-PC scheme 
with $\Nuc=2$ upwind corrections.  
For the SPIE scheme, implicit time-stepping is used for the boundary-fitted annular grid with radial stretching near the small hole, while explicit time-stepping is used for the Cartesian background grid.
Figure~\ref{fig:smallHoleContours} shows the computed solution at $t=1$ for the total field, $u$, 
the scattered field, $u_s$, and the error in $u$. Note that there is a significant scattered field for this case even though
the radius of the cylinder $a=0.01$ is fairly small compared to the wavelength, $2\pi/k \approx 0.63$, of the incident field.
The error is seem to be smooth with the largest errors distributed throughout the domain; there are no particularly large errors
in the vicinity of the hole.

The right graph in Figure~\ref{fig:smallHoleGridAndErrorsFig} compares the max-norm errors over time for 
the EME4 and SPIE4 schemes. The error in the EME scheme starts out smaller but then becomes similar in magnitude
to the errors in the SPIE4 scheme.
The time-step for the EME4 scheme is approximately $30$ times smaller than that for the SPIE4 scheme, and the  CPU time
required to compute the solution at $t=10$ using the SPIE4 scheme is approximately 20 times smaller than that needed for the EME4 scheme.

\newcommand{\Gcha}{\Gc_{\rm h,a}}
\newcommand{\Gcho}{\Gc_{\rm h,o}}
\subsubsection{Scattering of a modulated Gaussian plane wave by a collection of small holes} \label{sec:holes}

We consider the scattering of a modulated Gaussian plane wave from two different arrays of small holes.
This example demonstrates an interesting scattering problem for a geometry with small geometric features for which
the SPIE scheme gives a good speedup over the explicit scheme.
The incident field consists of a modulated Gaussian plane wave traveling from left to right and given
by the formula
\ba
   u(\xv,t) = e^{-\beta ( x-x_0 - ct )^2 } \, \cos( 2\pi\, k_0 \, (x-x_0-c t )) , \label{eq:modulatedGaussian}
\ea
where the Gaussian shape parameter is $\beta=20$, the modulation wave-number is $k_0=8$, and the center of the pulse is at $x_0=-2$ initially.

Two configurations of holes are considered for the region $[-3,2]\times[-2,2]$.
The first, called the \textit{aligned-hole} configuration, contains an array of $M_x \times M_y$ holes, each of radius $a=0.01$, with $M_x=7$ and $M_y=26$.
The centers of the holes are located at
\ba
   \xv_{m_x,m_y} = 
              \begin{bmatrix} m_x \, s_x   \\ m_y \, s_y \end{bmatrix} 
      - \half \begin{bmatrix} (M_x-1)s_x \\ (M_y-1) s_y \end{bmatrix}, 
     \quad m_x=0,1,\ldots, M_x, \quad m_y=0,1,\ldots, M_y,
\ea
where $s_x=0.15$ and $s_y=0.15$ denote the hole separations in the $x$ and $y$ directions, respectively.
The second configuration, called the \textit{offset-hole} configuration, also contains $M_x=7$ columns of holes, 
but every second column is shifted vertically by $s_x/2$ and contains $27$ holes instead of $26$.

{
\newcommand{\figw}{5.75cm}
\newcommand{\figh}{5.75cm}
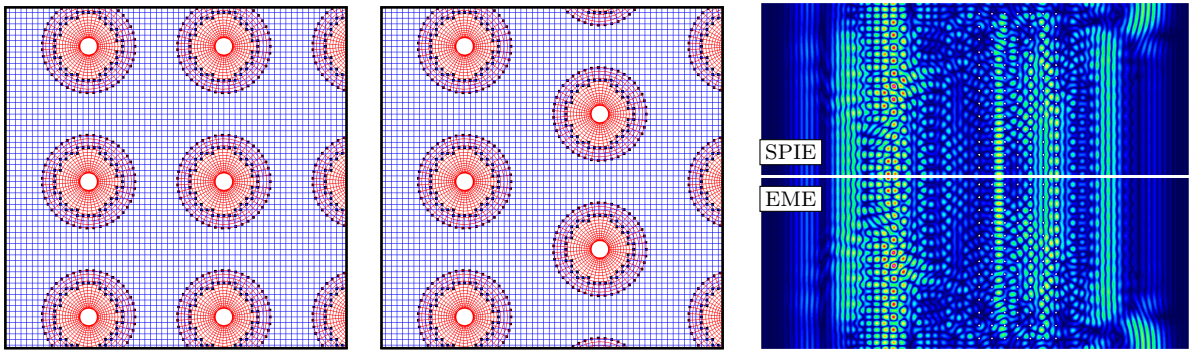
\begin{figure}[htb]
\begin{center}
\begin{tikzpicture}
  \useasboundingbox (0,.7) rectangle (16,4.75);  

\figByWidthb{0.0}{0.1}{holesGridG16}{4.5cm}[0.][0.1][0.05][0.05]
  \figByWidthb{5}{0.1}{offsetHolesGridG16}{4.5cm}[0.][0.1][0.05][0.05]

  \begin{scope}[xshift=10cm,yshift=2.4cm]
     \figByWidth{0}{0}{holesRad01G16N7M26SPIEabsu36}{\figw}[0.12][.12][.5][.19]
     \draw (0,.25) node[draw,fill=white,anchor=west,xshift=0pt,yshift=2pt,inner sep=2pt] {\scriptsize SPIE};
  \end{scope}
  \begin{scope}[xshift=10cm]
     \figByWidth{0}{0}{holesRad01G16N7M26EMEabsu36}{\figw}[0.12][.12][.19][.5]
     \draw (0,2.1) node[draw,fill=white,anchor=west,xshift=0pt,yshift=0pt,inner sep=2pt] {\scriptsize EME};
  \end{scope}        

\end{tikzpicture}
\end{center}
\caption{Left: Closeup of the aligned hole grid.
    Middle: Closeup of the offset hole grid.
    The white dots on the right plots are small holes with a grid around each as shown on the left.
    Right: comparison of the implicit-explicit SPIE solution (top half of computation) 
    to the explicit EME scheme (bottom half of computation) on the aligned hole grid. The results are nearly
    indistinguishable.
 }
\label{fig:holesGridFig}
\end{figure}
}

The overset grid for the aligned-hole configuration, denoted by $\Gcha^{(j)}$,
consists of a background Cartesian grid for the region $[-3,2]\times[-2,2]$,
together with small annular grids around the holes as shown in Figure~\ref{fig:holesGridFig}.
The nominal grid spacing is $\ds^{(j)} = 1/(10 j)$ with the grid lines on the annulii slightly smaller
and clustered near the boundary as shown in the figure.
The overset grid for the offset-hole configuration, denoted by $\Gcho^{(j)}$, has a similar construction to that
of the aligned-hole grid following its placement of holes.
The boundary conditions are taken as Dirichlet on the holes, Dirichlet on the left and right ends of the outer rectangle and
periodic in the y-direction of the outer rectangle.

{
\newcommand{\drawContour}[7]{%
\begin{scope}[#1]
\draw(0.0,0) node[anchor=south west,xshift=-4pt,yshift=+0pt] {\trimfiga{#2}{\figWidtha}};
  \draw(.5,.5) node[draw,fill=white,anchor=west,xshift=2pt,yshift=1pt,inner sep=2pt] {\scriptsize #3};
\begin{scope}[xshift=0cm,yshift=0pt]
  \draw (\xcb,\ycb) node[anchor=south west,xshift=0.25cm,yshift=.5cm,rotate=-90] {\trimfigcb{colourBarLines}{\cbWidth}{\cbHeight}};
  \draw (.8,0) node[anchor=north,xshift=+3pt,yshift=+2pt] {\scriptsize $#6$};
  \draw (4.8,0) node[anchor=north,xshift=+0pt,yshift=+2pt] {\scriptsize $#7$};
\end{scope}
\end{scope}
}
\newcommand{\cbWidth}{.2cm}
\newcommand{\cbHeight}{4cm}
\newcommand{\xcb}{.5cm}
\newcommand{\ycb}{-.2cm}
\setlength{\ycbTop}{\ycb+\cbHeight}
\setlength{\ycbMid}{\ycb+\cbHeight*\real{.5}}
\newcommand{\trimfigcb}[3]{\includegraphics[width=#2, height=#3, clip, trim=17cm 2.35cm 1.65cm 2.35cm]{#1}}
\newcommand{\figWidtha}{5cm}
\newcommand{\trimfiga}[2]{\trimw{#1}{#2}{.12}{.12}{.18}{.19}}
\begin{figure}[htb]
\begin{center}
\begin{tikzpicture}
   \useasboundingbox (0.75,.3) rectangle (16,9.5);  

   \begin{scope}[yshift=5cm]
      \drawContour{xshift=0.cm,yshift=0.00cm}{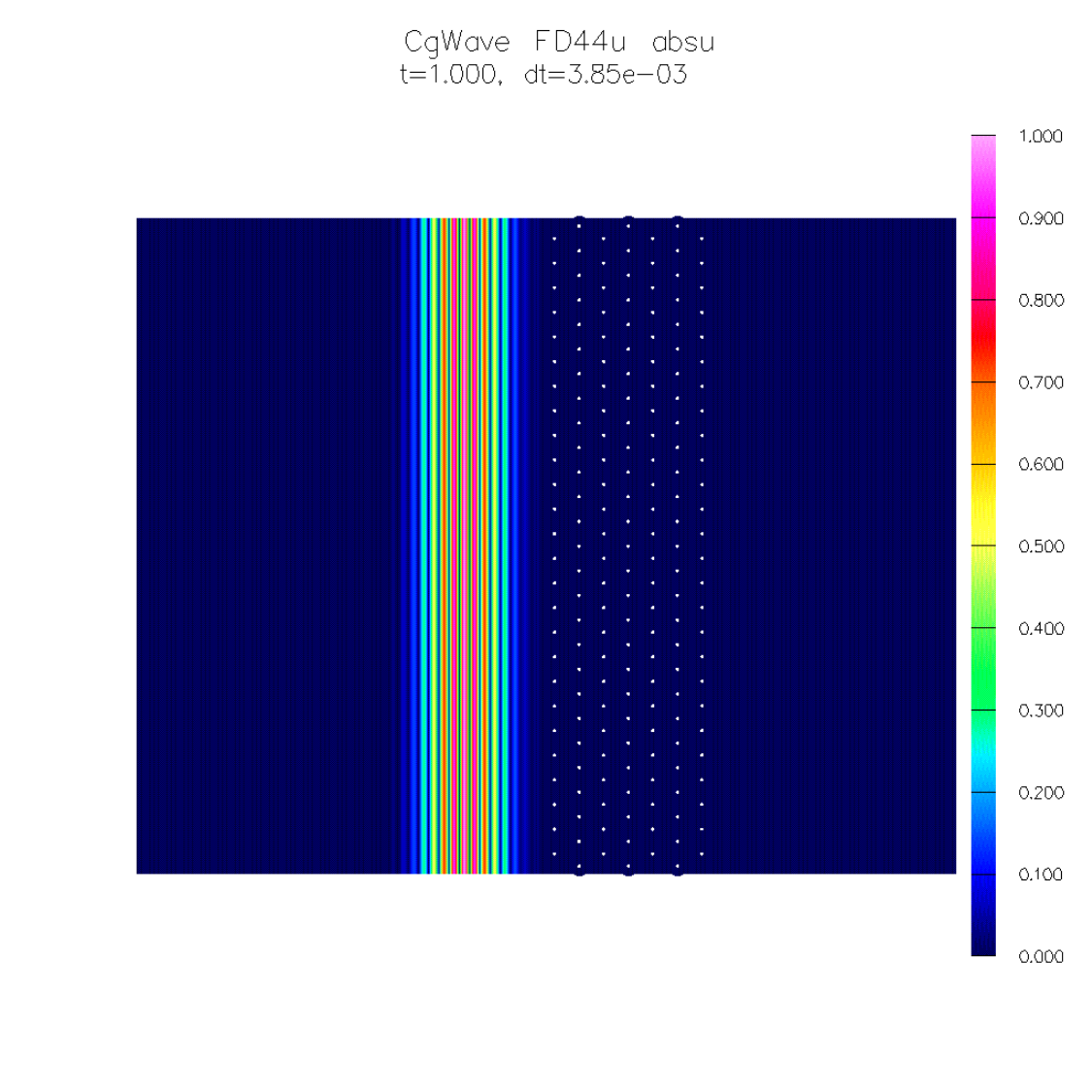}{$|u|$, $t=1.0$}{$v$}{$t=1.0$}{$0.0$}{$1.0$}     
     \drawContour{xshift=5.5cm,yshift=0.00cm}{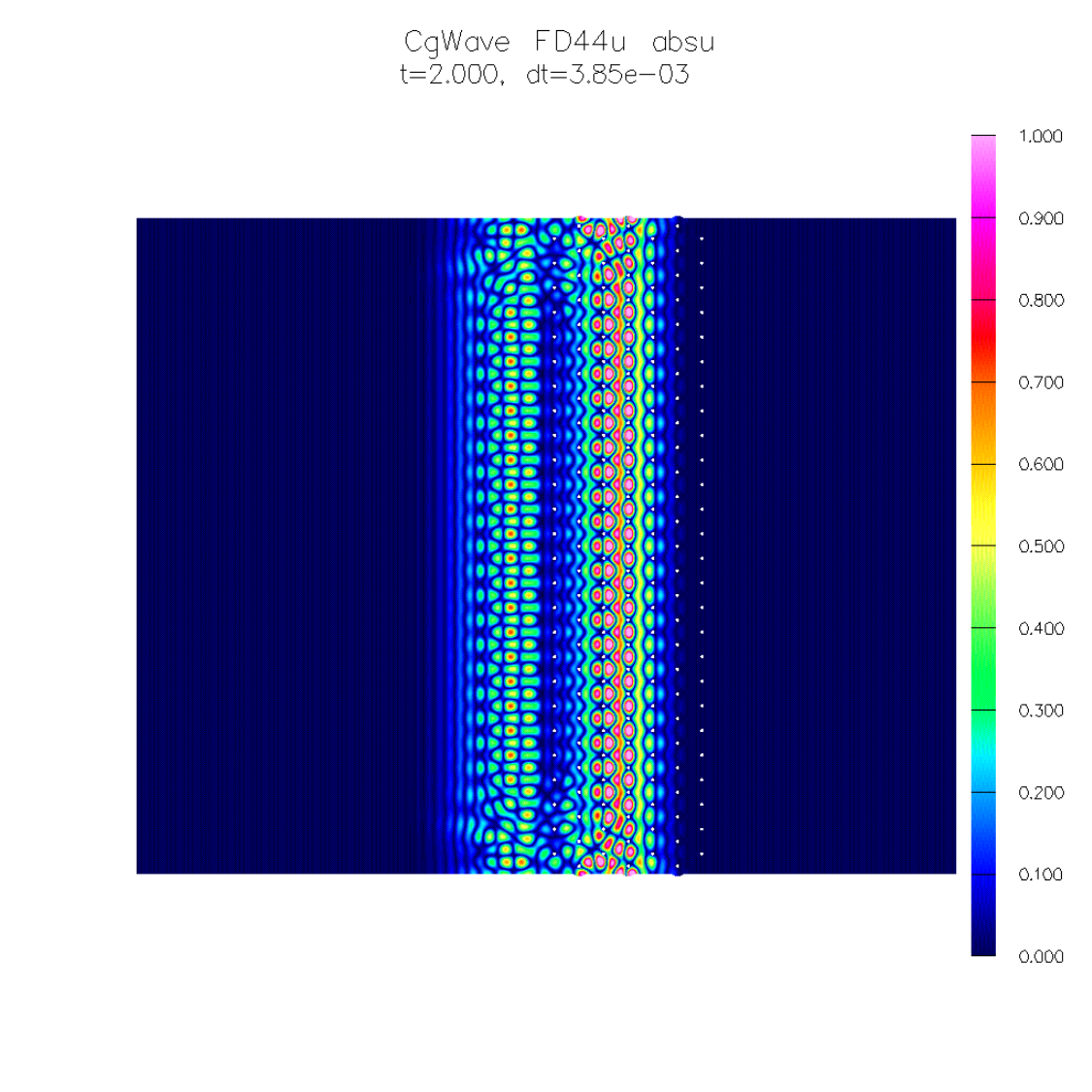}{$|u|$, $t=2.0$}{$v$}{$t=1.0$}{$0.0$}{$1.0$}   
      \drawContour{xshift=11cm,yshift=0.00cm}{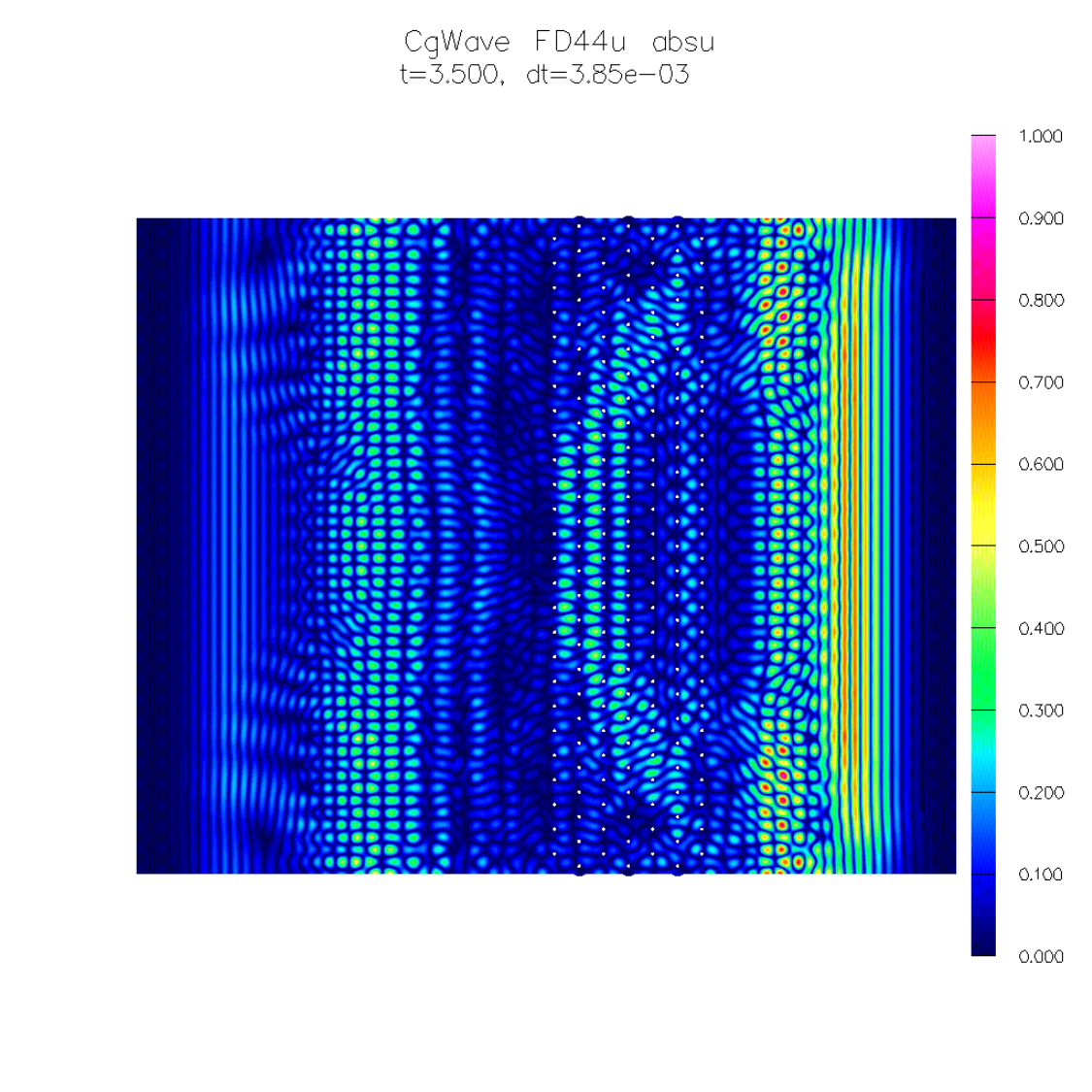}{$|u|$, $t=3.5$}{$v$}{$t=1.0$}{$0.0$}{$1.0$}
      \draw(5,3.5) node[draw,fill=white,anchor=west,xshift=2pt,yshift=1pt,inner sep=3pt] {\scriptsize Offset holes};   
      \draw[thick,white,->] (5,3.5) -- (4,3.);   
      \draw[thick,white,->] (1.75,2.5) -- (2.9,2.5);   
   \end{scope}
   \begin{scope}[yshift=0cm]
     \drawContour{xshift=0.0cm,yshift=0.cm}{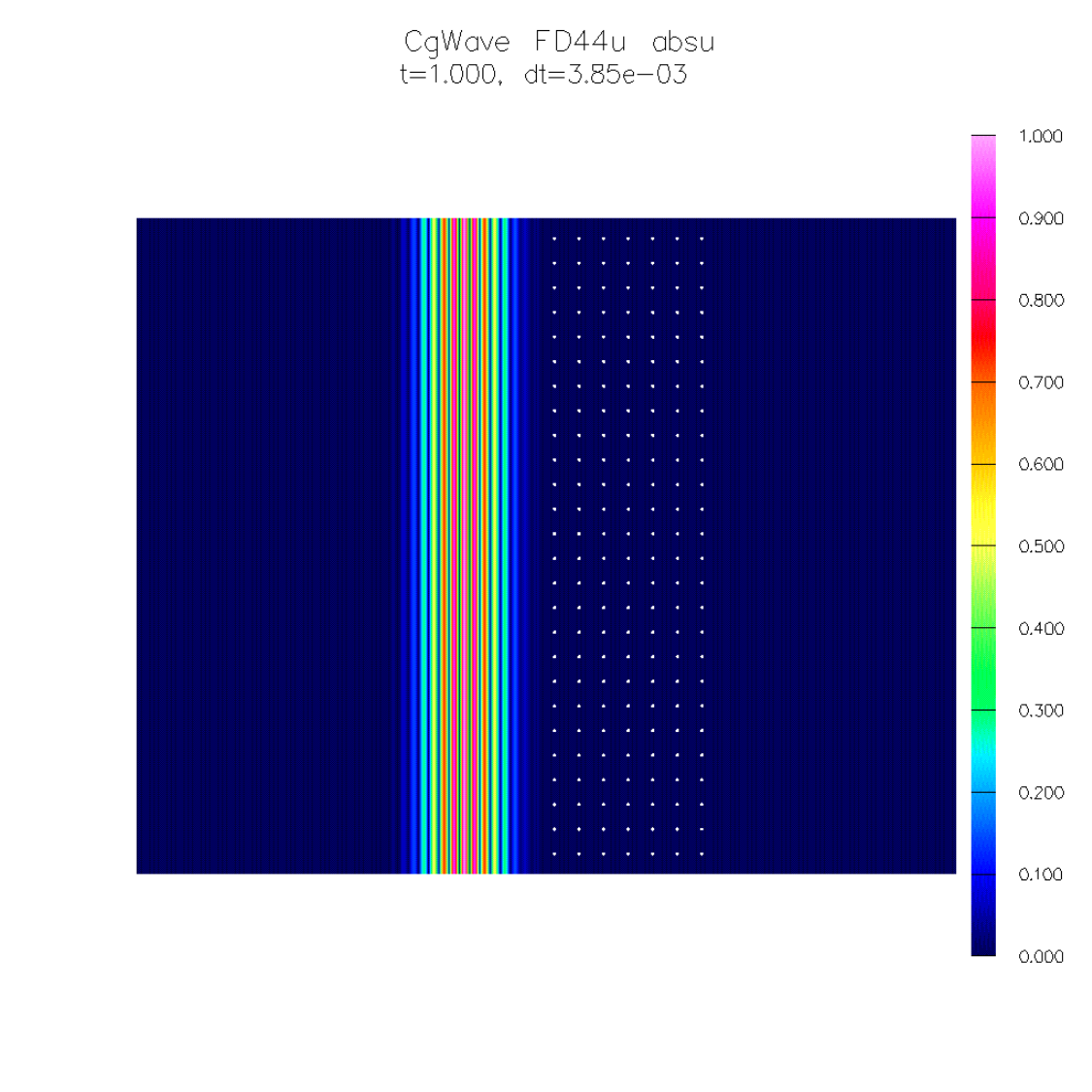}{$|u|$, $t=1.0$}{$v$}{$t=1.0$}{$0.0$}{$1.0$}    
     \drawContour{xshift=5.5cm,yshift=0.cm}{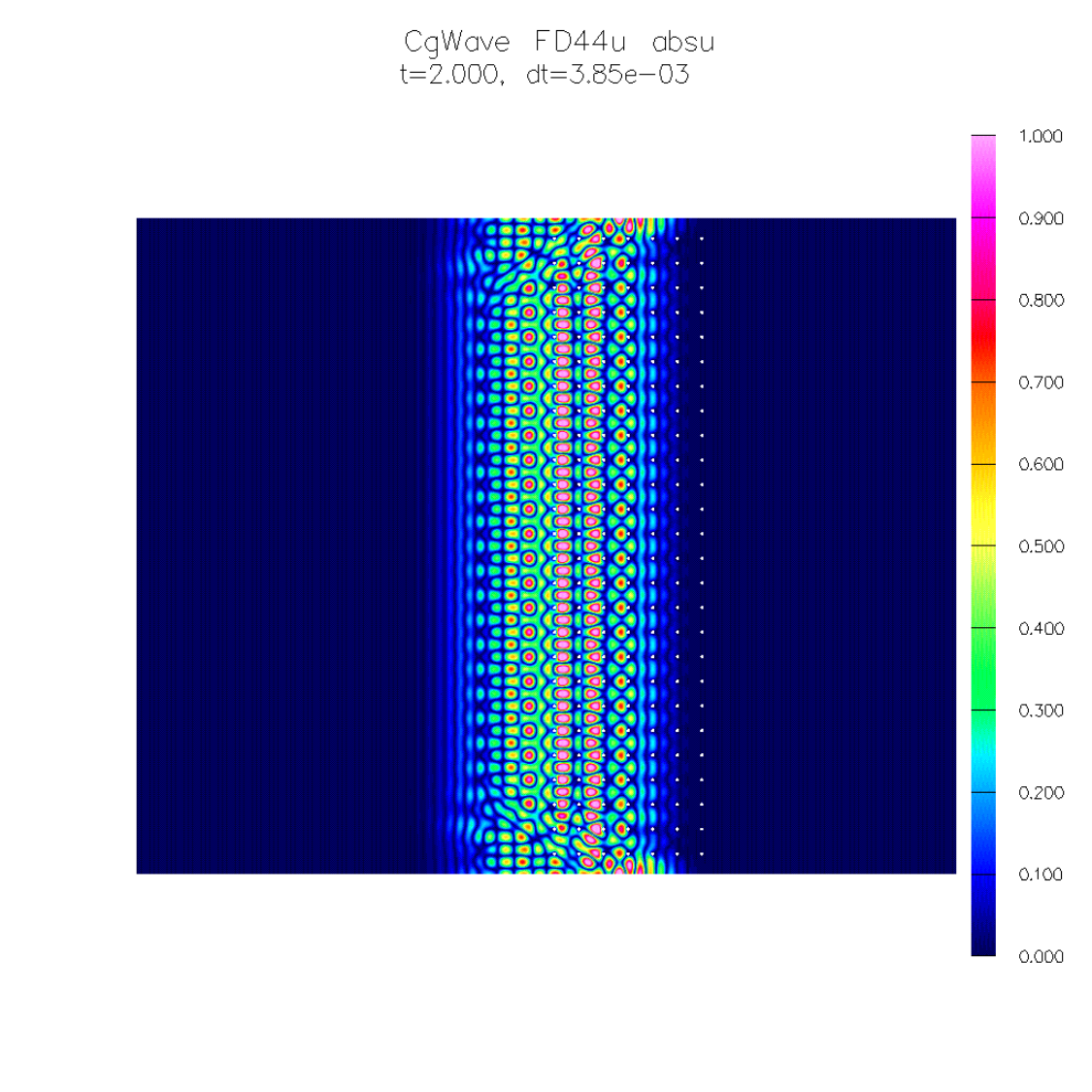}{$|u|$, $t=2.0$}{$v$}{$t=1.0$}{$0.0$}{$1.0$}    
     \drawContour{xshift=11.cm,yshift=0.cm}{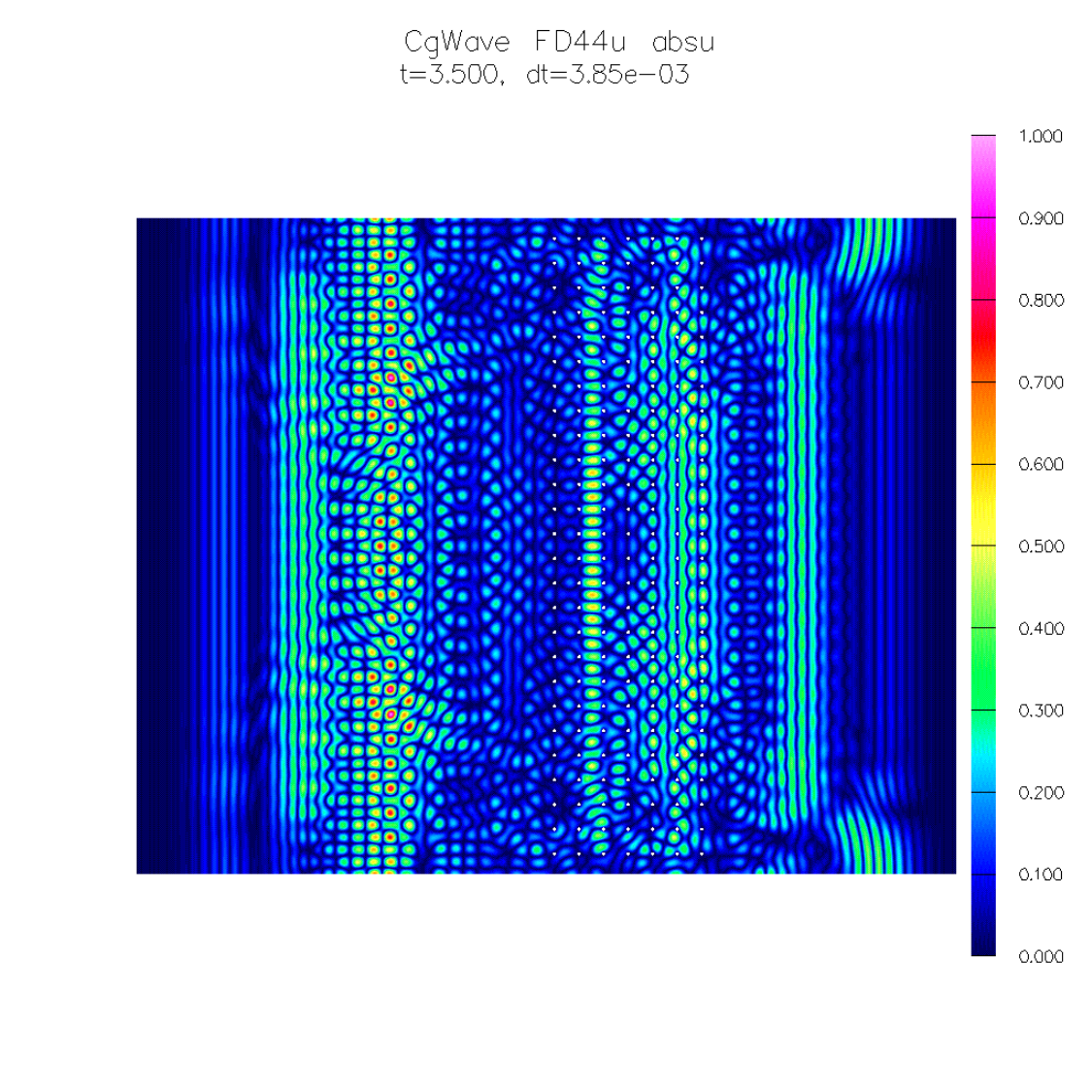}{$|u|$, $t=3.5$}{$v$}{$t=1.0$}{$0.0$}{$1.0$}    
    \draw(5,3.5) node[draw,fill=white,anchor=west,xshift=2pt,yshift=1pt,inner sep=3pt] {\scriptsize Aligned holes};
    \draw[thick,white,->] (5,3.5) -- (4,3.);   
    \draw[thick,white,->] (1.75,2.5) -- (2.9,2.5);    
   \end{scope}   

\end{tikzpicture}
\end{center}
\caption{Scattering of a modulated Gaussian plane wave by small holes (the white dots are small holes with a grid around each as shown in Figure~\ref{fig:holesGridFig}).
   Top: Offset holes. Bottom: Aligned holes.
    }
\label{fig:holes}
\end{figure}

}

Figure~\ref{fig:holes} shows the solution at three times for the two grid configurations of aligned and offset holes.
The solution is computed with the SPIE4-UW-PC scheme on grids $\Gcha^{(16)}$ and $\Gcho^{(16)}$.
Note that there are some edge effects in the solutions near the top and bottom periodic boundaries of the domain, 
due to the arrangement of the hole grids near these boundaries.
The solution at $t=1$ shows the incident Gaussian plane wave just starting to impact the first column of holes.
At $t=2$ the wave has travelled through most of the holes and a reflected wave is beginning to appear.
By $t=3.5$ most of the incident wave has been reflected or transmitted, although some residual wave motion 
resides within the array of holes.
Perhaps surprisingly, the transmitted wave is much stronger for the offset arrangement of holes.

Returning to Figure~\ref{fig:holesGridFig}, the right plot compares contours of the solutions computed using the SPIE and EME schemes.
The top half of the plot shows the SPIE4 solution, while the bottom half shows the EME4 solution.
After accounting for the reflection symmetry about the horizontal centerline, the results are nearly indistinguishable.
The speedup of the SPIE scheme over the EME scheme was about a factor of~$2$ for this case.
The SPIE time-step is about $4$ times that for the EME scheme. 
A better implementation of the implicit solvers should show an even bigger speedup of perhaps a factor of~$3$ or more (see the comments in Section~\ref{sec:implicitSolvers}).

\newcommand{\Gcke}{\Gc_{\rm ke}}
\subsubsection{Scattering of a modulated Gaussian plane wave from a knife edge} \label{sec:knifeEdge}

In this example, a modulated Gaussian plane wave given by~\eqref{eq:modulatedGaussian}
travels from left to right and diffracts off a thin \textsl{knife edge}
as shown in Figures~\ref{fig:knifeEdge} and~\ref{fig:knifeEdgeContours}.
This example demonstrates a problem that is geometrically stiff due to a sharp corner in the domain geometry, 
and one for which only a small portion of the overset grid is treated implicitly.

The overset grid for the geometry, denoted by $\Gcke^{(j)}$, is shown in Figure~\ref{fig:knifeEdge},
and consists of four component grids. A background Cartesian grid covers the domain $[-1.25,1]\times[0,1]$.
Two other Cartesian grids lie adjacent to the lower sides of the knife edge which has a total height of $0.5$.
A curvilinear grid is used
over the tip of the knife edge.
The nominal grid spacing was $\ds^{(j)}=1/(10 j)$, although the tip grid used a finer mesh with stretching to resolve
the sharp tip of the knife edge.

{
\newcommand{\drawContour}[7]{%
\begin{scope}[#1]
\draw(0.0,0) node[anchor=south west,xshift=-4pt,yshift=+0pt] {\trimfiga{#2}{\figWidtha}};
  \draw(.5,4.8) node[draw,fill=white,anchor=west,xshift=2pt,yshift=1pt,inner sep=2pt] {\scriptsize #3};
\begin{scope}[xshift=0cm,yshift=-2pt]
  \draw (\xcb,\ycb) node[anchor=south west,xshift=0.25cm,yshift=.5cm,rotate=-90] {\trimfigcb{colourBarLines}{\cbWidth}{\cbHeight}};
  \draw (.8,0) node[anchor=north,xshift=+3pt,yshift=+2pt] {\scriptsize $#6$};
  \draw (4.8,0) node[anchor=north,xshift=+0pt,yshift=+2pt] {\scriptsize $#7$};
\end{scope}
\end{scope}
}
\newcommand{\cbWidth}{.2cm}
\newcommand{\cbHeight}{4cm}
\newcommand{\xcb}{.5cm}
\newcommand{\ycb}{-.2cm}
\setlength{\ycbTop}{\ycb+\cbHeight}
\setlength{\ycbMid}{\ycb+\cbHeight*\real{.5}}
\newcommand{\trimfigcb}[3]{\includegraphics[width=#2, height=#3, clip, trim=17cm 2.35cm 1.65cm 2.35cm]{#1}}
\newcommand{\figWidtha}{5.15cm}
\newcommand{\trimfiga}[2]{\trimw{#1}{#2}{.05}{.12}{.07}{.1}}
\begin{figure}[htb]
\begin{center}
\begin{tikzpicture}
   \useasboundingbox (0.5,.2) rectangle (16,5.2);  

   \begin{scope}[xshift=0cm,yshift=-5pt]
      \drawContour{xshift=0.0cm,yshift=0.cm}{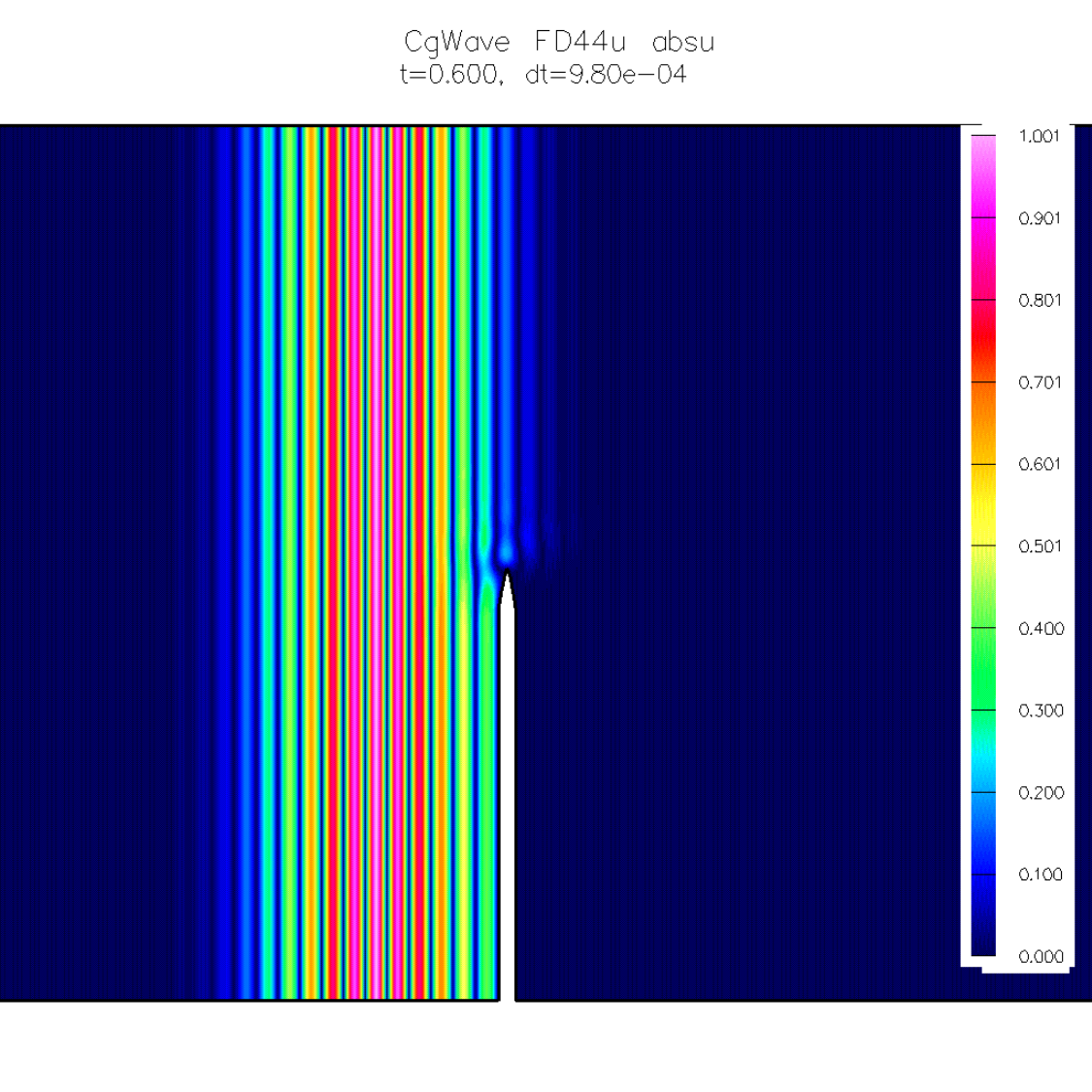}{$|u|$, $t=.6$}{$v$}{$t=1.0$}{$0.$}{$1.$}     
     \drawContour{xshift=5.25cm,yshift=0.cm}{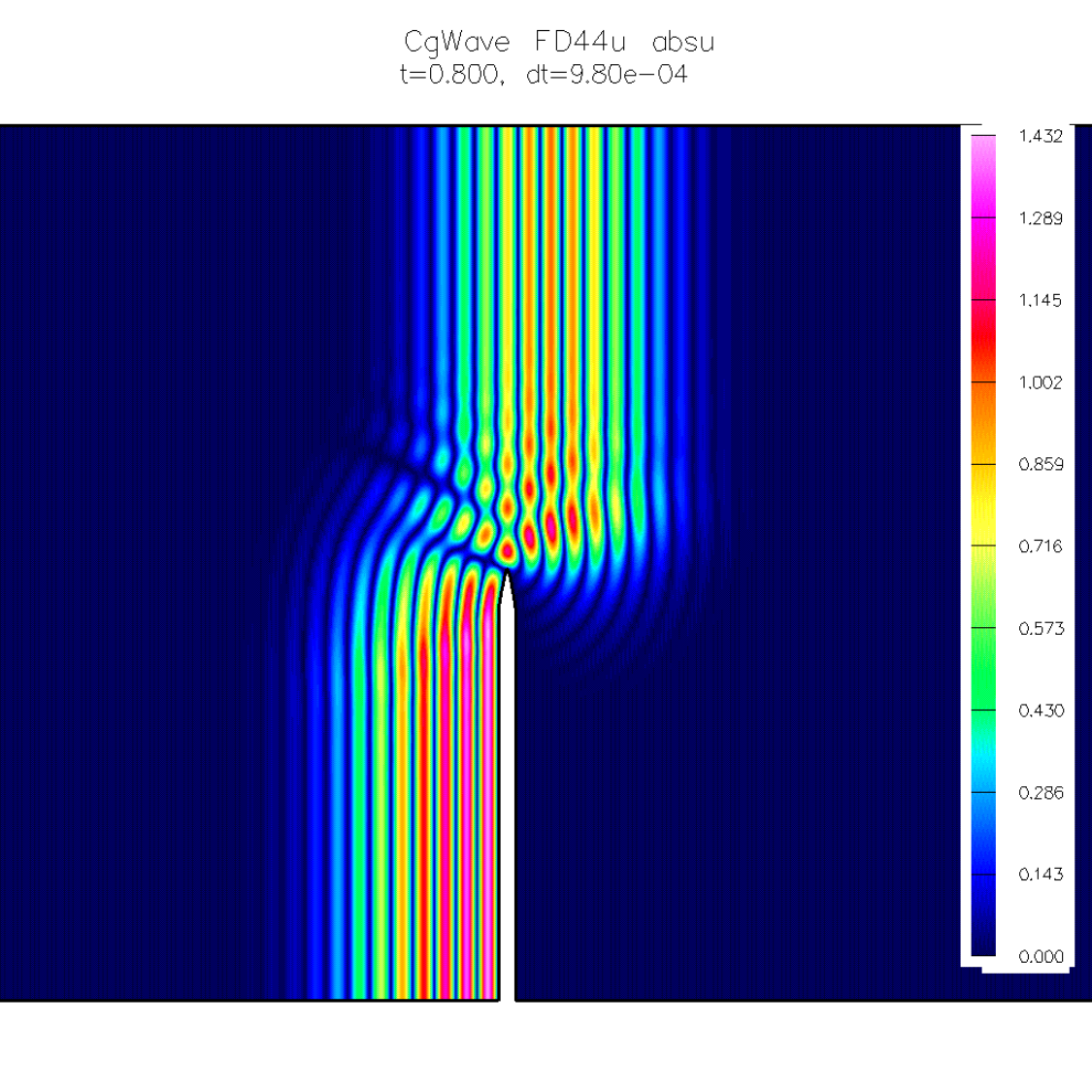}{$|u|$, $t=.8$}{$v$}{$t=1.0$}{$0$}{$1.4$}     
     \drawContour{xshift=10.5cm,yshift=0.cm}{tipGridG64O5SPIEk20_showabsu21}{$|u|$, $t=1$}{$v$}{$t=1.0$}{$0$}{$1.2$}   
   \end{scope}   

\end{tikzpicture}
\end{center}
\caption{Scattering of a modulated Gaussian plane wave from a knife edge. 
Contours of the solution at times $t=.6,.8,1$ for the modulation wave number $k_0=20$ using the fourth-order accurate
 implicit-explicit scheme SPIE4-UW-PC. The time-step was about $20$ times larger compared to the corresponding explicit scheme.
    }
\label{fig:knifeEdgeContours}
\end{figure}

}

Figure~\ref{fig:knifeEdgeContours} shows contours of the solution for three times 
computed on grid $\Gcke^{(16)}$ using the fourth-order accurate SPIE4-UW-PC scheme.
The Gaussian is centered at $x_0=-.75$ initially and the Gaussian shape parameter is taken as $\beta=80$.
Neumann boundary conditions are used on the outer boundaries of the domain and a Dirichlet condition is used on the knife edge.
The SPIE scheme is used with only the curvilinear tip grid treated implicitly.  As a result, the scheme is able to use a time-step that is about $20$ times larger than that required by the fully explicit EME-UW-PC scheme.
The speedup factor over the fully explicit scheme is found to be about~$11$ for both the second and fourth-order accurate SPIE schemes. 
Note that the tip grid has just $1,760$ grid points out of a total of $928,765$, or $0.2\%$ of the points.
Obviously a more efficient implementation of the implicit solver should lead to speedups closer to a factor of $20$, a task for future work.

\section{Conclusions} \label{sec:conclusions}

We have described and analyzed a class of new implicit and implicit-explicit time-stepping methods for the numerical solution
of the wave equation in second-order form.
These single-step, three time-level, schemes are based on the modified equation (ME) approach.  
Second and fourth-order accurate schemes are developed, although the approach supports higher-order accurate schemes.
The coefficient matrix implied by the implicit scheme is definite and well suited for solution by modern Krylov methods or multigrid.
Conditions for accuracy and unconditional stability of the implicit ME (IME) schemes are derived.
Several approaches for incorporating upwind dissipation into the IME schemes are discussed. 
A predictor-corrector approach that adds the upwinding in a separate explicit step appears to be quite useful.
For problems on overset grids that are geometrically stiff due to locally small cells, we have developed a spatially
partitioned implicit-explicit (SPIE) scheme whereby component grids with small cells are integrated with
the IME scheme while others grids use an explicit ME (EME) scheme. We have shown that for geometrically stiff problems
the resulting SPIE scheme can be many times faster and more accurate that using the EME scheme everywhere.
Although developed for the wave equation, the new schemes can be extended to other wave propagation problems written
in second-order form such as
Maxwell's equations of electromagnetics, elasticity, and acoustics.


\appendix

\newcommand{\LpHat}{\widehat{\Lc_{p,h}}}
\newcommand{\LtwoHat}{\widehat{\Lc_{2,h}}}
\newcommand{\LfourHat}{\widehat{\Lc_{4,h}}}
\section{Stability proofs} \label{sec:stabilityProofs}

For the stability analyses we consider a Cartesian grid on a
$2\pi$ periodic domain $\Omega=[0,2\pi]^{\nd}$.
Von Neumann analysis expands the solution in a discrete Fourier series in space.
The stability condition is enforced by ensuring each Fourier mode satisfies the condition.

\subsection{Stability of the second-order accurate implicit ME scheme (IME2)} \label{sec:stabilityProof2IME2}

Here is the proof of Theorem~\ref{thm:stabilityIME2}, the statement of which is repeated here for clarity.
\begin{theoremNoNumber}[IME2 Stability]
\label{thm:stabilityIME2A}
 The IME2 scheme~\eqref{eq:ME2discrete} is unconditionally stable on a periodic domain provided
 \ba
     \alpha_2 \ge \f{1}{4}.       
\ea
\end{theoremNoNumber}

\begin{proof}

We look for solutions consisting of a single Fourier mode,
\ba
  U_{\jv}^n  = a^n \, e^{i\kv\cdot\xv_{\jv}} , \label{eq:vonNemannAnsatz}
\ea
where $a$ is the amplification factor and $\kv=[k_1,k_2,k_3]^T$ is the vector of wave numbers,
with $k_d=-N_d/2,-N_d/2+1,\ldots,N_d/2-1$, assuming $N_d$ is even.
Substituting the anstaz~\eqref{eq:vonNemannAnsatz} into~\eqref{eq:ME2discrete}
leads to a quadratic equation for $a$
\bse
\ba  
   & a^2 - 2 b \, a + 1 = 0 , 
\intertext{where}
 & b \eqdef  \f{1 + (\alpha_2- \half) \,\lamHat_2^2 \, z}{ 1+ \alpha_2 \lamHat_2^2 \, z }, \\
 &  \lamHat_2^2 \eqdef - \LtwoHat, \quad z\eqdef \dt^2, 
\ea
\ese
and where $\LtwoHat$ is the Fourier symbol of $\Lc_{2,h}$, 
\ba
   & \LtwoHat \eqdef   \, c^2 \sum_{d=0}^{\nd} \, \f{-4\sin^2(k_d h_d)}{h_d^2} .
\ea
Note that $\lamHat_2^2 \ge 0$, with strict inequality $\lamHat_2^2 >0$ when $\kv\ne 0$.
It is not hard to show that for stability (Definition~\eqref{def:stability}) we require $b\in\Real$ and $|b|<1$
(the end cases $b=\pm 1$ lead to double roots $a=\pm 1$ and linearly growing modes).
Thus, when $\kv\ne\zerov$, we require
\ba
    -1 < \f{1 + (\alpha_2- \half) \,\lamHat_2^2 \, z}{ 1+ \alpha_2 \lamHat_2^2 \, z } < 1, 
        \label{eq:stabilityFractionIME2}
\ea
for all $z>0$. 
The right inequality in~\eqref{eq:stabilityFractionIME2} gives
\bse
\ba
   & 1 + (\alpha_2- \half) \,\lamHat_2^2 \, z <  1+ \alpha_2 \lamHat_2^2 \, z , \\
 \implies  & - \half \,\lamHat_2^2 \, z <  0 , 
\ea
\ese
which is always true.
The left inequality implies
\bse
\ba
         & -( 1+ \alpha_2 \lamHat_2^2 \, z) <   1 + (\alpha_2- \half) \,\lamHat_2^2\, z, \\
\implies &  -2  + \half \lamHat_2^2 z <  2 \alpha_2 \lamHat_2^2 \, z , \\
\implies & \alpha_2 > \f{1}{4} - \f{1}{\lamHat_2^2 \, z}.
\ea
\ese
Therefore we require
\ba
   \alpha_2 \ge \f{1}{4},
\ea
and this completes the proof.
\end{proof}

\subsection{Stability of the fourth-order accurate implicit ME scheme (IME4)} \label{sec:stabilityProof2IME4}

Here is the proof of Theorem~\ref{thm:stabilityIME4}, the statement of which is repeated here for clarity.
\begin{theoremNoNumber}[IME4 Stability]
\label{thm:stabilityIME4A}
 The IME4 scheme~\eqref{eq:ME2discrete} is unconditionally stable on a periodic domain provided
\bse
\label{eq:stabilityConditionIME4A}
\ba
    & \alpha_2 \ge \f{1}{12}, \\
    & \alpha_4 \ge 
       \begin{cases}
           \f{1}{4} \alpha_2 - \f{1}{48},                                  &   \text{when $\alpha_2\ge \f{1}{4}$}, \\
           \f{1}{4} \alpha_2 - \f{1}{48} + \f{8}{9} (\f{1}{4}-\alpha_2)^2, &  \text{when $\f{1}{12} \le \alpha_2 \le \f{1}{4}$ }.
       \end{cases} 
\ea
\ese
\end{theoremNoNumber}

\begin{proof}
Using the anstaz~\eqref{eq:vonNemannAnsatz} in~\eqref{eq:ME2discrete} leads to 
following quadratic for the time-stepping amplification factor $a$ ,
\ba
  & a^2 - 2 b \, a + 1 = 0 , 
\intertext{where} 
& b \eqdef \f{1 - \half \beta_2 \,\lamHat_4^2 \, z - \half \beta_4 \, \lamHat_2^4 \, z^2}{ 1+ \alpha_2 \lamHat_4^2 \, z + \alpha_4  \lamHat_2^4\, z^2}, \\
&  \lamHat_4^2 \eqdef - \LfourHat, \quad
  \lamHat_2^2 \eqdef   -\LtwoHat, \quad z\eqdef \dt^2 ,
\ea
and where
\ba  
   \beta_2 = 1-2\alpha_2, \quad \beta_4 = \alpha_2 -2 \alpha_4 - \f{1}{12}. \label{eq:beta4alpha4}
\ea
For stability we require that $|b|<1$ for $\kv\ne\zerov$, 
\ba
  & \Big|  \f{1 - \half \beta_2 \lamHat_4^2 \, z - \half \beta_4 \lamHat_2^4 \, z^2}{ 1+ \alpha_2 \lamHat_4^2 \, z + \alpha_4 \lamHat_2^4 \, z^2} 
     \Big| < 1,     \label{eq:IME4StabFormula}
\ea
which will give constraints on $\alpha_2$ and $\alpha_4$.
Requiring~\eqref{eq:IME4StabFormula} leads to two conditions,
\bse
\label{eq:IME4stabConditions}
\ba
  &  1 - \half \beta_2 \lamHat_4^2 \, z - \half \beta_4 \lamHat_2^4 \, z^2 < 1+ \alpha_2 \lamHat_4^2 \, z + \alpha_4 \lamHat_2^4 \, z^2, \\
  & 
  -(1+ \alpha_2 \lamHat_4^2 \, z + \alpha_4 \lamHat_2^4 \, z^2) <  1 - \half \beta_2 \lamHat_4^2 \, z - \half \beta_4 \lamHat_2^4 \, z^2 .
\ea
\ese
These can be simplified to 
\bse
\ba
  &  \half  \lamHat_4^2 \, z  + ( \half \alpha_2 - \f{1}{24} )\lamHat_2^4 \, z^2 > 0 ,  \label{eq:IME4stab1} \\
  &  (\alpha_4 -\half\beta_4) \lamHat_2^4 \, z^2 + (2\alpha_2-\half)  \lamHat_4^2 \, z + 2 >0 .\label{eq:IME4stab2}
\ea
\ese
Inequality~\eqref{eq:IME4stab1} must hold for all $z>0$ which implies 
\ba    
   \alpha_2 \ge \f{1}{12}. 
\ea
Inequality~\eqref{eq:IME4stab2} is a quadratic inequality in $z=\dt^2$,
\ba  
   &  A z^2 + B z + C > 0, \\
   & A \eqdef (\alpha_4 -\half\beta_4) \lamHat_2^4, \quad B\eqdef (2\alpha_2-\half)  \lamHat_4^2, \quad C\eqdef 2, 
\ea
which must hold for all $z>0$.
This quadratic must be flat or concave upward 
which implies $A = \alpha_4 -\half\beta_4 \ge 0$ or 
\ba
  \alpha_4 \ge \f{1}{4} \alpha_2 - \f{1}{48} \ge 0 ,  \text{~~when~~} \f{1}{12} \le \alpha_2 .
\ea
The minimum of the quadratic with $z\ge 0$ occurs when $z_m=-B/(2A) \ge 0 $, which implies $B \le 0$ or
$\alpha_2\le \f{1}{4}$. The minimum value of the quadratic is $C - B^2/(4A) $ and this should be greater than or equal to zero which implies
$B^2 \le 4 A C $ or 
\ba
   (2\alpha_2-\half)^2  \lamHat_4^4 \le 8 (\alpha_4 -\half\beta_4) \lamHat_2^4, \quad \text{~~when~~} \f{1}{12} \le \alpha_2\le \f{1}{4}. 
\ea
This last inequality is re-arranged as a condition on $\alpha_4$ in terms of $\alpha_2$, (using~\eqref{eq:beta4alpha4}),
\ba
   \alpha_4 \ge \f{1}{4} \alpha_2 - \f{1}{48} + \f{1}{2} (\alpha_2-\f{1}{4})^2 \, \f{\lamHat_4^4}{\lamHat_2^4} , \quad \text{~~when~~} \f{1}{12} \le \alpha_2\le \f{1}{4}. 
\ea
Using $\f{\lamHat_4^4}{\lamHat_2^4} \le (4/3)^2$ in the last term gives 
\ba
   \alpha_4 \ge \f{1}{4} \alpha_2 - \f{1}{48} + \f{8}{9} (\f{1}{4}-\alpha_2)^2 , \quad \text{~~when~~} \f{1}{12} \le \alpha_2\le \f{1}{4}. 
\ea
In summary $\alpha_2$ must satisfy
\ba
   & \alpha_2 \ge \f{1}{12}, 
\ea
while $\alpha_4$ is constrained by 
\ba
    & \alpha_4 \ge 
       \begin{cases}
           \f{1}{4} \alpha_2 - \f{1}{48}, &  \text{when $\alpha_2\ge \f{1}{4}$}, \\
           \f{1}{4} \alpha_2 - \f{1}{48} + \f{8}{9} (\f{1}{4}-\alpha_2)^2 ,&  \text{when $\f{1}{12} \le \alpha_2 \le \f{1}{4}$ },
       \end{cases}
\ea
which completes the proof.
\end{proof}

\newcommand{\qHat}{\hat{q}}
\newcommand{\QvHat}{\hat{\Qv}}
\newcommand{\LambdaHat}{\hat{\Lambda}}
\subsection{IME stability with single stage upwind dissipation (IME-UW)} \label{sec:stabilityIME-UW}

In this section we prove Theorem~\ref{th:IME-UW}, which is repeated here or clarity.
\begin{theoremNoNumber} \label{th:IME-UWA}
  The IME-UW schemes~\eqref{eq:IME-UW} for $p=2,4$ on a periodic or infinite domain Cartesian grid 
  are unconditionally stable for any $\nu_p>0$ provided $\alpha_2$ satisfies the conditions of
  Theorem~\ref{thm:stabilityIME2}, for $p=2$, or $\alpha_2$ and $\alpha_4$ satisfy the conditions for 
  Theorem~\ref{thm:stabilityIME4} for $p=4$.
\end{theoremNoNumber}

\begin{proof}
We show the proof for $p=4$, the case $p=2$ is similar.
Using the anstaz~\eqref{eq:vonNemannAnsatz} in~\eqref{eq:IME-UW} leads to 
following quadratic for the time-stepping amplification factor $a$ ,
\ba
  & a^2 - 2 b \, a + c = 0 , \label{eq:genericQuadratic}
\intertext{where} 
& b \eqdef \f{1 - \half \beta_2 \,\lamHat_4^2 \, \dt^2 - \half \beta_4 \, \lamHat_2^4 \, \dt^4}
             {1+ \LambdaHat + \f{\nu_p}{2} \dt\, \qHat_p^2 }, \\
& c \eqdef \f{1+ \LambdaHat - \f{\nu_p}{2} \dt\, \qHat_p^2}
             {1+ \LambdaHat + \f{\nu_p}{2} \dt\, \qHat_p^2} ,  \label{eq:productOfRoots} \\
&  \LambdaHat \eqdef \alpha_2 \lamHat_4^2 \, \dt^2 + \alpha_4  \lamHat_2^4\, \dt^4 , \\
&  \qHat_p^2 \eqdef \QvHat_p,  
\ea
and where $\QvHat_p >0 $ (for $\kv\ne0$) is the symbol of the dissipation operator $\Qv_p$ in~\eqref{eq:Qdef}, 
\ba  
   & \QvHat_p \eqdef \sum_{d=1}^{\nd} \f{c}{h_d} \left[ 4 \sin^2(h_d/2)  \right]^{p/2+1} .
\ea
The conditions for stability come from the theory of Schur and von Neumann polynomials~\cite{Miller1971,Strikwerda89}
which for the quadratic~\eqref{eq:genericQuadratic} are
\bse
\label{eq:IME-UWvonNeumannPoly}
\ba
   &  |c| < 1 ,   \label{eq:IME-UWvonNeumannPolyA} \\
   & |b| \le \half | 1 + c |  . \label{eq:IME-UWvonNeumannPolyB} 
\ea
\ese
Note that the magnitude of the product of the roots $c$ in~\eqref{eq:productOfRoots} is now less than $1$, $|c|<1$,  when $\nu_p>0$, since we have assumed 
$\LambdaHat > 0$ for $\kv\ne\zerov$. Thus the first condition~\eqref{eq:IME-UWvonNeumannPolyA} is true.
Note that 
\ba
  \half(1 + c) &= \f{ 1+ \LambdaHat + \f{\nu_p}{2} \dt\, \qHat_p^2 + 1+ \LambdaHat - \f{\nu_p}{2} \dt\, \qHat_p^2} 
              {1+ \LambdaHat + \f{\nu_p}{2} \dt\, \qHat_p^2}  
        = \f{ 1+ \LambdaHat }{1+ \LambdaHat + \f{\nu_p}{2} \dt\, \qHat_p^2} ,
\ea
and thus $1+c>0$.
The inequality~\eqref{eq:IME-UWvonNeumannPolyB} thus
requires the two conditions
\bse
\ba
 & \f{1 - \half \beta_2 \,\lamHat_4^2 \, \dt^2 - \half \beta_4 \, \lamHat_2^4 \, \dt^4}
             {1+ \LambdaHat + \f{\nu_p}{2} \dt\, \qHat_p^2 }
     \le 
         \f{1+ \LambdaHat}
          {1+ \LambdaHat + \f{\nu_p}{2} \dt\, \qHat_p^2}, \\
& -  \f{1+ \LambdaHat}
          {1+ \LambdaHat + \f{\nu_p}{2} \dt\, \qHat_p^2} \le
\f{1 - \half \beta_2 \,\lamHat_4^2 \, \dt^2 - \half \beta_4 \, \lamHat_2^4 \, \dt^4}
             {1+ \LambdaHat + \f{\nu_p}{2} \dt\, \qHat_p^2 }
\ea
\ese
or upon multiplying through by the denominator, 
\bse
\ba           
 & 1 - \half \beta_2 \,\lamHat_4^2 \, \dt^2 - \half \beta_4 \, \lamHat_2^4 \, \dt^4
    \le 1+ \LambdaHat , \\
 & -( 1+ \LambdaHat ) \le 1 - \half \beta_2 \,\lamHat_4^2 \, \dt^2 - \half \beta_4 \, \lamHat_2^4 \, \dt^4 .
\ea
\ese
These last two conditions (note that $\nu_p$ has dropped out) are satisfied since these are essentially the same
inequalities~\eqref{eq:IME4stabConditions}
hold from Theorem~\ref{thm:stabilityIME4A} (the only difference is that $\le$ is replaced by $<$ in~\eqref{eq:IME4stabConditions}).
This proves the theorem.
\end{proof}

\subsection{IME Stability with predictor-corrector upwind dissipation (IME-UW-PC)} \label{sec:stabilityIME-UW-PC}

Here is the proof of theorem~\ref{th:IME-UW-PC-MUSTA}.
\begin{proof}
We prove the result for $p=2$, the proof for $p=4$ follows in a similar fashion.
The second-order accurate IME-UW-PC scheme with multiple stages is
\bse
\label{eq:IME-UW-PC-MUSTAa}
\ba
 & \f{U_\jv^{(0)} -2 U_\jv^n + U_\jv^{n-1}}{\dt^2 } = \Lv_{\alpha p}\, (U_\jv^{(0)}, U_\jv^n, U_\jv^{n-1}) , \\
 & U_\jv^{n+1} =  \Rv_p^{\Nuc} U_\jv^{(0)} + (I- \Rv_p^{\Nuc} ) U_\jv^{n-1}, 
\ea
\ese
where
\ba
    \Rv_p \eqdef I - \f{\nu_p \dt}{2} \Qv_p.
\ea
Substituting the ansatz $U_\jv^{(0)} = \Uhat^{(0)} \, e^{i\kv\cdot\xv_{\jv}}$ and 
$U_\jv^n = \Uhat^n\,  e^{i\kv\cdot\xv_{\jv}}$ leads to 
\bse
\ba
  & \f{\Uhat^{(0)} -2 \Uhat^n + \Uhat^{n-1}}{\dt^2 } 
     =  - \lamHat_2^2 \Big( \alpha_2 \, \Uhat^{(0)} + (1-2\alpha_2) \,\Uhat^{n}  + \alpha_2 \, \Uhat^{n-1} \Big), \label{eq:IME-UW-PCa2} \\
 & \Uhat^{n+1} = \Rhat_p^{\Nuc} \Uhat_\jv^{(0)} + (I- \Rhat_p^{\Nuc} ) \Uhat_\jv^{n-1} ,        \label{eq:IME-UW-PCb2} 
\ea
\ese
where
\ba
    \Rhat_p = 1 - \f{\nu_p \dt}{2} \qHat_p.
\ea
Solving~\eqref{eq:IME-UW-PCa2} for $\Uhat^{(0)} $ ,
\ba
 & \Uhat^{(0)} = 
        2 \, \f{1 - \lamHat_2^2 \dt^2 (\half-\alpha_2) }{1 + \alpha_2 \lamHat_2^2 \dt^2} \, \Uhat^n 
       -  \Uhat^{n-1}               
\ea
and substituting into~\eqref{eq:IME-UW-PCb2}  gives 
\bse
\ba
    \Uhat^{n+1} & = \f{\Rhat_p^{\Nuc}}{1 +\alpha_2 \dt^2\, \lamHat_2^2  }
        \left[ 
         2 \Uhat^n - \Uhat^{n-1} - \dt^2 \lamHat_2^2\, \big( (1-2\alpha_2) \,\Uhat^{n}  + \alpha_2 \, \Uhat^{n-1} \big) 
          \right] 
         + (1 - \Rhat_p^{\Nuc} )\, \Uhat^{n-1}, \\
   & =  2  \f{ 1 - \dt^2 \lamHat_2^2\, (\half -\alpha_2)}{1 +\alpha_2 \dt^2\, \lamHat_2^2  } \Rhat_p^{\Nuc} \,\Uhat^{n}
          +  (1 - 2 \Rhat_p^{\Nuc} ) \Uhat^{n-1} .
\ea
\ese
Now looking for solutions of the form $\Uhat^n = c_0\, a^n$ for some constant $c_0$ leads to a quadratic equation for a,
\ba
   &  a^2 - 2 b a + c = 0, 
\ea
where
\ba
   &  b \eqdef  \Rhat_p^{\Nuc}\,  \f{ 1 - \dt^2 \lamHat_2^2\, (\half -\alpha_2)}{1 +\alpha_2 \dt^2\, \lamHat_2^2  }  , \\
   &  c \eqdef -1  + 2 \Rhat_p^{\Nuc} 
\ea
For stability we require the two conditions~\eqref{eq:IME-UWvonNeumannPoly} from Section~\ref{sec:stabilityIME-UW},
\ba
   | c | <1                 & \implies | 1 - 2 \Rhat_p^{\Nuc} | < 1, \label{eq:MUSTAcondition1}\\
   |b | \le \half |1 + c |  & \implies | \Rhat_p |^{\Nuc} 
    \left| 
       \f{ 1 - \dt^2 \lamHat_2^2\, (\half -\alpha_2)}{1 +\alpha_2 \dt^2\, \lamHat_2^2  }
    \right| \le | \Rhat_p |^{\Nuc}   \label{eq:MUSTAcondition2}
\ea
If we assume the parameters $\alpha_2$ and $\lamHat_2$ are chosen to make the scheme without dissipation stable then
\ba
      \left| 
       \f{ 1 - \dt^2 \lamHat_2^2\, (\half -\alpha_2)}{1 +\alpha_2 \dt^2\, \lamHat_2^2  }
    \right| \le 1
\ea
and~\eqref{eq:MUSTAcondition2} is satisfied. Condition~\eqref{eq:MUSTAcondition1} implies $0 < \Rhat_p^\Nuc <1 $ or 
\ba
    0 < \big( 1 - \f{\nu_p \dt}{2} \qHat_p \big)^\Nuc < 1 
\ea
which implies (ignoring the special case when $\nu_p \qHat_2=0$ )
\ba
  \f{\nu_p \dt}{2} \qHat_p \, < \,
  \begin{cases}
      2 & \text{if $\Nuc$ is even}, \\
      1 & \text{if $\Nuc$ is odd} .
  \end{cases} 
\ea
The conclusions of the proof now follow.
\end{proof}

\subsection{Proof of a lemma} \label{sec:GKSlemma}

Here is the proof of Lemma~\ref{thm:gksLemmaIME2}.
\begin{proof}
   If $|\kappa|=1$ then it can be written as $\kappa=e^{i\theta}$ for $\theta\in \Real$.
   Then, using $\kappa-2+\kappa^{-1} = -4\sin^2(\theta/2)$ gives 
   \ba
    & b = \frac{1 - 4 (\half-\alpha_2)\, \lambda^2 \, \sin^2(\theta/2)}
                                          {1 +4 \alpha_2 \, \lambda^2 \, \sin^2(\theta/2)}.    
   \ea
   Note that $b\in\Real$ and $|b| \le 1$ since $b\le 1$ implies
   \ba
      & 1 - 4 (\half-\alpha_2)\, \lambda^2 \, \sin^2(\theta/2) \le 1 +4 \alpha_2 \, \lambda^2 \, \sin^2(\theta/2), \\
   \implies
      & -2 \lambda^2 \, \sin^2(\theta/2)  \le 0 ,
   \ea
   which is true, while $b \ge -1$ implies
   \ba
      & -(1 +4 \alpha_2 \, \lambda^2 \, \sin^2(\theta/2)) \le 1 - 4 (\half-\alpha_2)\, \lambda^2 \, \sin^2(\theta/2) , \\
      & \implies   (1 - 4 \alpha_2) \lambda^2 \, \sin^2(\theta/2) \le 1,
   \ea
   which holds when $\lambda<1$ and $\alpha_2\ge0$, or for any $\lambda>0$ when $\alpha_2 \ge 1/4$.
   Now, when $b\in \Real$ and $|b|\le 1$ then the magnitude of the roots $a$ of~\eqref{eq:lemmaIME2quadratic} satisfy
   \ba
      | a | = |  b \pm \sqrt{b^2-1} |= |  b \pm i \, \sqrt{1 - b^2} |  = \sqrt{ b^2 + (1-b^2)} =1.
   \ea
   This proves the lemma.
\end{proof}

\bibliographystyle{elsart-num}
\bibliography{journal-ISI,henshaw,henshawPapers,wave}

\begin{thebibliography}{10}
\expandafter\ifx\csname url\endcsname\relax
  \def\url#1{\texttt{#1}}\fi
\expandafter\ifx\csname urlprefix\endcsname\relax\def\urlprefix{URL }\fi

\bibitem{TaubeDumbserMunzSchneider2009}
A.~Taube, M.~Dumbser, C.-D. Munz, R.~Schneider, A high-order discontinuous
  {Galerkin} method with time-accurate local time stepping for the {Maxwell}
  equations, International Journal of Numerical Modelling: Electronic Networks,
  Devices and Fields 22~(1) (2009) 77--103.
\newline\urlprefix\url{https://onlinelibrary.wiley.com/doi/abs/10.1002/jnm.700}

\bibitem{appelo2020waveholtz}
D.~Appelo, F.~Garcia, O.~Runborg, {WaveHoltz}: Iterative solution of the
  {Helmholtz} equation via the wave equation, SIAM J. Sci. Comput. 42~(4)
  (2020) A1950--A1983.

\bibitem{EmWaveHoltzPengAppelo2022}
T.~Rylander, A.~Bondeson, {EM-WaveHoltz}: A flexible frequency-domain method
  built from time-domain solvers, IEEE Transactions on Antennas and Propogation
  70~(7) (2022) 5659 -- 5671.

\bibitem{appeloElWaveHoltz2022}
D.~Appelö, F.~Garcia, A.~{Alvarez Loya}, O.~Runborg, {El-WaveHoltz}: A
  time-domain iterative solver for time-harmonic elastic waves, Computer
  Methods in Applied Mechanics and Engineering 401 (2022) 115603.
\newline\urlprefix\url{https://www.sciencedirect.com/science/article/pii/S0045782522005655}

\bibitem{GroteMehlinMitkova2015}
M.~J. Grote, M.~Mehlin, T.~Mitkova, {Runge}--{Kutta}-based explicit local
  time-stepping methods for wave propagation, SIAM Journal on Scientific
  Computing 37~(2) (2015) A747--A775.
\newline\urlprefix\url{https://doi.org/10.1137/140958293}

\bibitem{AlmquistMehlin2017}
M.~Almquist, M.~Mehlin, Multilevel local time-stepping methods of
  {Runge--Kutta}-type for wave equations, SIAM Journal on Scientific Computing
  39~(5) (2017) A2020--A2048.
\newline\urlprefix\url{https://doi.org/10.1137/16M1084407}

\bibitem{Dablain1986}
M.~Dablain, High order differencing for the scalar wave equation, Geophysics 51
  (1986) 54--66.

\bibitem{ShubinBell1987}
G.~R. Shubin, J.~B. Bell, A modified equation approach to constructing fourth
  order schemes for acoustic wave propagation, SIAM J. Sci. Stat. Comput. 8~(2)
  (1987) 135--151.

\bibitem{Berger3}
M.~J. Berger, J.~Oliger, Adaptive mesh refinement for hyperbolic partial
  differential equations, J. Comput. Phys. 53 (1984) 484--512.

\bibitem{LiuLiHu2014}
L.~Liu, X.~Li, F.~Q. Hu, Nonuniform-time-step explicit {Runge–Kutta} scheme
  for high-order finite difference method, Computers \& Fluids 105 (2014)
  166--178.
\newline\urlprefix\url{https://www.sciencedirect.com/science/article/pii/S0045793014003454}

\bibitem{DiazGrote2009}
J.~Diaz, M.~J. Grote, Energy conserving explicit local time stepping for
  second-order wave equations, SIAM Journal on Scientific Computing 31~(3)
  (2009) 1985--2014.
\newline\urlprefix\url{https://doi.org/10.1137/070709414}

\bibitem{BeznosovAppelo2021}
O.~Beznosov, D.~Appel\"o, Hermite-discontinuous {Galerkin} overset grid methods
  for the scalar wave equation, Communications on Applied Mathematics and
  Computation 3~(3) (2021) 391--418.
\newline\urlprefix\url{https://doi.org/10.1007/s42967-020-00075-5}

\bibitem{BrittTurkelTsynkov2018}
S.~Britt, E.~Turkel, S.~Tsynkov, A high order compact time/space finite
  difference scheme for the wave equation with variable speed of sound, Journal
  of Scientific Computing 76 (2018) 777--811.
\newline\urlprefix\url{https://api.semanticscholar.org/CorpusID:207198477}

\bibitem{LiLiaoLin2019}
K.~Li, W.~Liao, Y.~Lin, A compact high order {Alternating Direction Implicit}
  method for three-dimensional acoustic wave equation with variable
  coefficient, Journal of Computational and Applied Mathematics 361 (2019)
  113--129.
\newline\urlprefix\url{https://www.sciencedirect.com/science/article/pii/S0377042719301992}

\bibitem{KahanaSmithTurkelTsynkov2022}
A.~Kahana, F.~Smith, E.~Turkel, S.~Tsynkov, A high order compact time/space
  finite difference scheme for the 2d and 3d wave equation with a damping
  layer, Journal of Computational Physics 460 (2022) 111161.
\newline\urlprefix\url{https://www.sciencedirect.com/science/article/pii/S0021999122002236}

\bibitem{LimKimDouglas2007}
H.~Lim, S.~Kim, J.~Douglas, Numerical methods for viscous and nonviscous wave
  equations, Appl. Numer. Math. 57~(2) (2007) 194–212.
\newline\urlprefix\url{https://doi.org/10.1016/j.apnum.2006.02.004}

\bibitem{KimLim2007}
S.~Kim, H.~Lim, High-order schemes for acoustic waveform simulation, Applied
  Numerical Mathematics 57~(4) (2007) 402--414.
\newline\urlprefix\url{https://www.sciencedirect.com/science/article/pii/S0168927406001012}

\bibitem{Piperno2006}
{Piperno, Serge}, Symplectic local time-stepping in non-dissipative {DGTD}
  methods applied to wave propagation problems, ESAIM: M2AN 40~(5) (2006)
  815--841.
\newline\urlprefix\url{https://doi.org/10.1051/m2an:2006035}

\bibitem{Verwer2011}
J.~Verwer, Component splitting for semi-discrete {Maxwell} equations, BIT
  Numerical Mathematics 51~(2) (2011) 427--445.
\newline\urlprefix\url{https://doi.org/10.1007/s10543-010-0296-y}

\bibitem{ChabassierImperiale2016}
J.~Chabassier, S.~Imperiale, Fourth-order energy-preserving locally implicit
  time discretization for linear wave equations, International Journal for
  Numerical Methods in Engineering 106~(8) (2016) 593--622.
\newline\urlprefix\url{https://onlinelibrary.wiley.com/doi/abs/10.1002/nme.5130}

\bibitem{mxsosup2018}
J.~Angel, J.~W. Banks, W.~D. Henshaw, High-order upwind schemes for the wave
  equation on overlapping grids: {Maxwell}'s equations in second-order form, J.
  Comput. Phys. 352 (2018) 534--567.

\bibitem{adegdm2019}
J.~Angel, J.~W. Banks, W.~D. Henshaw, M.~J. Jenkinson, A.~V. Kildishev,
  G.~{Kova\v ci\v c}, L.~J. Prokopeva, D.~W. Schwendeman, A high-order accurate
  scheme for {M}axwell's equations with a generalized dispersion model, J.
  Comput. Phys. 378 (2019) 411--444.

\bibitem{adegdmi2020}
J.~W. Banks, B.~Buckner, W.~D. Henshaw, M.~J. Jenkinson, A.~V. Kildishev,
  G.~{Kova\v ci\v c}, L.~J. Prokopeva, D.~W. Schwendeman, A high-order accurate
  scheme for {M}axwell's equations with a generalized dispersive material
  ({GDM}) model and material interfaces, J. Comput. Phys. 412 (2020) 109424.

\bibitem{maxwellMLA2022}
Q.~Xia, J.~W. Banks, W.~D. Henshaw, A.~V. Kildishev, G.~Kovačič, L.~J.
  Prokopeva, D.~W. Schwendeman, High-order accurate schemes for {M}axwell's
  equations with nonlinear active media and material interfaces, J. Comput.
  Phys. 456 (2022) 111051.

\bibitem{ssmx2023}
J.~B. Angel, J.~W. Banks, A.~Carson, W.~D. Henshaw, Efficient upwind
  finite-difference schemes for wave equations on overset grids, J. Comput.
  Phys. 45~(5) (2023) A2703--A2724.

\bibitem{smog2012}
D.~Appel\"o, J.~W. Banks, W.~D. Henshaw, D.~W. Schwendeman, Numerical methods
  for solid mechanics on overlapping grids: Linear elasticity, J. Comput. Phys.
  231~(18) (2012) 6012--6050\citeCount{24}.

\bibitem{flunsi2016}
J.~W. Banks, W.~D. Henshaw, A.~Kapila, D.~W. Schwendeman, An added-mass
  partitioned algorithm for fluid-structure interactions of compressible fluids
  and nonlinear solids, J. Comput. Phys. 305 (2016) 1037--1064\citeCount{8}.

\bibitem{ism2023}
J.~W. Banks, W.~D. Henshaw, A.~Newell, D.~W. Schwendeman, Fractional-step
  finite difference schemes for incompressible elasticity on overset grids, J.
  Comput. Phys. 488.

\bibitem{max2006b}
W.~D. Henshaw, A high-order accurate parallel solver for {Maxwell}'s equations
  on overlapping grids, SIAM J. Sci. Comput. 28~(5) (2006)
  1730--1765\citeCount{55}.

\bibitem{lcbc2022}
N.~G. {Al Hassanieh}, J.~W. Banks, W.~D. Henshaw, D.~W. Schwendeman, Local
  compatibility boundary conditions for high-order accurate finite-difference
  approximations of {PDE}s, SIAM J. Sci. Comput. 44 (2022) A3645--A3672.

\bibitem{Strand1994}
B.~Strand, Summation by parts for finite difference approximations for d/dx,
  J.\ Comput.\ Phys. 110 (1994) 47--67.

\bibitem{Olsson1995}
P.~Olsson, Summation by parts, projections, and stability. {II}, Math.\ Comput.
  64 (1995) 1473--1493.

\bibitem{Mattsson-Nordstrom-2004}
K.~Mattsson, J.~Nordstr{\"o}m, Summation by parts operators for finite
  difference approximations of second derivatives, J. Comput. Phys. 199 (2004)
  503--540.

\bibitem{AppeloPetersson2009}
D.~Appel\"o, N.~A. Petersson, A stable finite difference method for the elastic
  wave equation on complex geometries with free surfaces, Communications in
  Computational Physics 5 (2009) 84--107.

\bibitem{DuruKreissMattsson2014}
K.~Duru, G.~Kreiss, K.~Mattsson, Stable and high-order accurate boundary
  treatments for the elastic wave equation on second-order form, SIAM Journal
  on Scientific Computing 36~(6) (2014) A2787--A2818.

\bibitem{sosup2012}
J.~W. Banks, W.~D. Henshaw, Upwind schemes for the wave equation in
  second-order form, J. Comput. Phys. 231~(17) (2012) 5854--5889\citeCount{3}.

\bibitem{AllisonCarson2023}
A.~M. Carson, High-order accurate implicit-explicit time-stepping schemes for
  wave equations on overset grids, Ph.D. thesis, Dept. of Mathematical
  Sciences, Rensselaer Polytechnic Institute (2023).

\bibitem{GKSI}
H.-O. Kreiss, Stability theory of difference approximations of mixed initial
  boundary value problems. {I}, Mathematics of Computation 22 (1968) 703--714.

\bibitem{GKSII}
B.~Gustafsson, H.-O. Kreiss, A.~Sundstr\"om, Stability theory of difference
  approximations for mixed initial boundary value problems. {II}, Mathematics
  of Computation 26~(119) (1972) 649--686.

\bibitem{CGNS}
G.~S. Chesshire, W.~D. Henshaw, Composite overlapping meshes for the solution
  of partial differential equations, J. Comput. Phys. 90~(1) (1990)
  1--64\citeCount{527}.

\bibitem{ogen}
W.~D. Henshaw, Ogen: An overlapping grid generator for {O}verture, Research
  Report UCRL-MA-132237, Lawrence Livermore National Laboratory (1998).

\bibitem{Miller1971}
J.~J.~H. Miller, On the location of zeros of certain classes of polynomials
  with applications to numerical analysis, IMA J. Appl. Math. 8~(3) (1971)
  397--406.

\bibitem{Strikwerda89}
J.~C. Strikwerda, Finite Difference Schemes and Partial Differential Equations,
  Wadsworth and Brooks/Cole, 1989.

\end{thebibliography}
 
\end{document}